\documentstyle[amscd]{amsart}
\newsymbol\varkappa 207B
\newsymbol\restriction 1316
\newsymbol\between 1347
\def\dj{d\kern-.30em\raise1.25ex\vbox{\hrule width .3em height .03em}}
\def\Dj{D\kern-.75em\raise0.75ex\vbox{\hrule width .3em height .03em}
\kern.03em}
\def\pod{\it\mediumseries Extended Version}
\makeatletter
\renewcommand{\subsection}{\@startsection{subsection}{2}{\z@}%
{\baselineskip}{0.5\baselineskip}{\defaultfont\bf}}
\makeatother
\newcommand{\rig}{\wp}
\newcommand{\restr}{{\restriction}}
\newcommand{\grten}{\mathbin{\widehat{\otimes}}}
\newcommand{\Lf}{\frak{l}(P)}
\newcommand{\J}{\frak{T}}
\newcommand{\JJ}{\Upsilon(P)}
\newcommand{\lie}{\mbox{\shape{n}\selectfont lie}}
\newcommand{\Sum}{\displaystyle{\sum}}
\newcommand{\hor}{\mbox{\family{euf}\shape{n}\selectfont hor}}
\newcommand{\ver}{\mbox{\family{euf}\shape{n}\selectfont ver}}
\newcommand{\con}{\mbox{\family{euf}\shape{n}\selectfont con}}
\newcommand{\rc}{\mbox{\family{euf}\shape{n}\selectfont r}}
\newcommand{\vh}{\mbox{\family{euf}\shape{n}\selectfont vh}}
\newcommand{\ad}{\mbox{\shape{n}\selectfont ad}}
\newcommand{\id}{\mbox{\shape{n}\selectfont id}}
\newcommand{\adj}{\varpi}
\newcommand{\k}{\kappa}
\newcommand{\e}{\epsilon}
\newcommand{\Mor}{\mbox{\shape{n}\selectfont Mor}}
\def\bla#1{$(${\it #1\/{}}$)$}
\def\lh{\ell}
\newcommand{\gen}{\mbox{\shape{n}\selectfont gen}}
\newcommand{\proj}{\varkappa}\newcommand{\Proj}{\rho_\perp}
\newcommand{\kG}{\sharp^\kappa}
\newcommand{\fG}{\sharp^\phi}
\newcommand{\inv}{i\!\hspace{0.8pt}n\!\hspace{0.6pt}v}
\newcommand{\tG}{\sharp}
\newcommand{\hr}{h\!or}
\newcommand{\Plft}{q_\star}
\newcommand{\Jlft}{\sharp_\star}
\newcommand{\Prig}{q^\star}
\newcommand{\Jrig}{\sharp^\sstar}
\newcommand{\im}{\mbox{\shape{n}\selectfont im}}
\def\sstar{{\raise0.2ex\hbox{$\scriptstyle\star$}}}
\def\L{\cal{L}}
\newtheorem{thm}{Theorem}[section]
\newtheorem{lem}[thm]{Lemma}
\newtheorem{pro}[thm]{Proposition}
\theoremstyle{definition}
\newtheorem{defn}{Definition}[section]
\numberwithin{equation}{section}
\begin{document}
\title[QUANTUM PRINCIPAL BUNDLES]{GEOMETRY OF QUANTUM PRINCIPAL BUNDLES II\\
\hfill \pod}
\author{Mi\'co {\Dj}ur{\Dj}evi\'c}
\address{Instituto de Matematicas, UNAM,
Area de la Investigacion Cientifica,
Circuito Exterior,
Ciudad Universitaria,
M\'exico DF, CP 04510, M\'EXICO\newline
\mbox{}\newline
\indent {\it Written In}\newline
\indent Faculty of Physics, University of Belgrade, pBox 550,
Studentski Trg 12, 11001 Beograd, SERBIA\newline
\indent {\it Extended Version}\newline
\indent Centro de Investigaciones Teoricas, UNAM, Facultad
de Estudios Superiores Cuautitlan, Cuautitlan Izcalli, M\'EXICO}
\maketitle
\begin{abstract}
A general noncommutative-geometric theory of
principal bundles  is
presented. Quantum groups  play  the  role  of  structure  groups.
General quantum  spaces  play  the  role  of  base  manifolds.  A
differential calculus on quantum principal bundles is studied.  In
particular, algebras of horizontal and  verticalized  differential
forms on the bundle are introduced and investigated. The formalism
of connections is developed.
Operators  of  horizontal  projection,  covariant
derivative and curvature are constructed and analyzed.  A  quantum
generalization  of  classical  Weil's  theory  of   characteristic
classes is  sketched.  Quantum  analogs  of  infinitesimal  gauge
transformations are studied. Illustrative examples and
constructions are presented.
\end{abstract}
\tableofcontents
\filbreak
\section{Introduction}
\renewcommand{\thepage}{}
     In this study we continue the presentation of the  theory  of
quantum principal bundles.

     The theory developed in the previous paper \cite{D}
was ``semiclassical'': structure groups were
considered as quantum  objects, however
base  spaces were classical smooth manifolds. Algebraic formalism
developed in the previous paper will be now  generalized  and  incorporated
into a completely quantum framework, following general  philosophy
of  non-commutative   differential   geometry \cite{C}.
Base   manifolds,
structure groups  and  corresponding  principal  bundles  will  be
considered as quantum objects.

     The paper is organized as follows.

     \smallskip
     Exposition of the theory begins in  Section 3,  with  a
general definition of quantum principal bundles.  This  definition
will translate into a noncommutative-geometric  context classical
idea that a principal bundle is a space on  which  the
structure  group acts  freely  on  the  right,  such  that
the  base  manifold  is diffeomorphic to the corresponding orbit space.

     After the main definition we pass  to  questions  related  to
differential calculus on quantum principal bundles. As first, we
introduce and  analyze  a  differential  *-algebra  consisting  of
``verticalized'' differential forms on the bundle. This algebra will
be introduced independently  of  a  specification  of  a  complete
differential calculus on the bundle.
\renewcommand{\thepage}{\arabic{page}}

     The calculus  on  the  bundle  is based  on  a
graded-differential *-algebra (representing  differential  forms)
possessing two important properties. As first, we  require
that this differential algebra is generated by ``functions'' on the
bundle. This condition  ensures  uniqueness  of  various  entities
naturally  appearing  in  the  study  of  differential   calculus.
Secondly, we postulate that the group action on ``the  functions''  on
the bundle is extendible to an  appropriate  differential  algebra
homomorphism (imitating the corresponding ``pull back''
of differential forms).

     Quantum counterparts of various important entities associated
to  differential  calculus  in  the  classical  theory   will   be
introduced in a constructive manner, starting from the algebra  of
differential forms on the bundle (and from  a  given  differential
calculus on the  structure  quantum group).  In  particular,  a graded
*-algebra representing horizontal forms
will  be  introduced  and  analyzed.
Also, graded-differential *-alge\-bras representing differential
forms  on  the base manifold and ``verticalized'' differential
forms on the bundle will be described.

     It is important to point out a conceptual difference  between
this approach to differential calculus and the approach presented
in  the  previous  paper, where a differential
calculus on the bundle was constructed starting from the  standard
calculus on the base manifold, and an appropriate calculus on
the structure quantum group. The main property of the  calculus
was a variant of {\it local triviality}, in
the  sense  that  all  local
trivializations of the bundle  locally  trivialize  the  calculus,
too. This property implies certain restrictions on  a  possible
differential calculus on the structure quantum group, as
discussed in details in \cite{D}. On the other hand, in this  paper  we
start from a fixed calculus on the group (based on the universal
differential envelope of a given first-order differential
structure) and the calculus on  the
base manifold  is  determined  by  the  calculus  on  the  bundle.
However, the calculus on the bundle is not uniquely determined  by
mentioned initial two conditions.

     In  Section 4  the  formalism  of  connections   will   be
presented. All  corresponding  basic  ``global'' constructions  and
results of the previous paper will be translated  into  the  general
quantum context. In particular, operators of horizontal  projection,
covariant  derivative  and  curvature  will  be constructed   and
investigated. Further, two  particularly  interesting  classes  of
connections will be  introduced  and  analyzed.  The  first  class
consists of connections possessing certain multiplicativity property.
This  is  a  trivial  generalization   of   multiplicative
connections of the previous paper. The second  class  consists  of
connections  that  are  counterparts  of   classical   connections
introduced in the  previous  paper.  Here,  these
connections will be called regular. Intuitively speaking,  regular
connections  are  ``maximally   compatible''   with the internal
geometrical structure of the bundle.

     In Section 5 a generalization of classical Weil's theory
of characteristic classes will be presented.

     Finally, in Section 6 some examples, remarks  and  additional
constructions are included. In particular,  we  shall  present  a
general {\it re}-construction of differential calculus on
the  bundle,
starting from  a  given  algebra  of  horizontal  forms,  and  two
operators imitating the covariant derivative and the curvature  of
a regular connection. Further, quantum
analogs of infinitesimal gauge transformations will be studied,
from two different viewpoints.

We shall also briefly discuss
interrelations with a theory of quantum principal
bundles presented in \cite{BM}.

Concerning  concrete
examples of quantum principal bundles, we shall  consider  trivial
bundles and principal bundles based on quantum homogeneous
spaces. The main structural elements of differential
calculus and the formalism of connections will be
illustrated on these examples.

The paper ends with two appendices. The first appendix is  devoted
to the analysis of the calculus on the bundle in the case when the
higher-order calculus on the structure quantum group is  described
by the corresponding bicovariant exterior  algebra  \cite{W3}.  In
particular, it will be shown that if the first-order  calculus  on
the group is
compatible with all ``transition functions'' (in  the  context  of
the previous paper) then the higher-order  calculus  based  on  the
exterior algebra possesses this property  too.  In  fact  this  is
equivalent to a possibility of constructing the calculus on  the
bundle such that all local trivializations of the  bundle  locally
trivialize  the  calculus.  Further,  we  shall  prove  that
bicovariant  exterior  algebras  describe,  in  a  certain  sense,
{\it the minimal\/} higher-order calculus on the group such that the
corresponding calculus  on  the  bundle  possesses  the  mentioned
trivializability property (universal envelopes always describe {\it the
maximal\/} higher-order calculus).  We  shall  also  analyze  similar
questions in the context of general theory.

In  the  second  appendix  the   structure   of   the    *-algebra
representing ``functions'' on the bundle is analyzed,
in the light of the decomposition  of  the  right  action  of  the
structure quantum group into multiple irreducible components.

\section{Preparatory Material}

     Before passing to quantum principal bundles we shall fix  the
notation, and introduce in the game relevant quantum group entities.
We  shall use  the  symbol $\grten$ for  a  graded  tensor  product
of graded (differential *-) algebras.

     Here, as in the previous paper, we shall  deal  with  compact
matrix quantum groups \cite{W2} only (however
the compactness assumption is not essential for
a large part of the formalism). Let $G$ be such a group. The algebra
of ``polynomial functions'' on $G$ will be denoted  by $\cal{A}$.
The  group
structure  on  $G$   is   determined   by   the   comultiplication
$\phi\colon \cal{A}\rightarrow\cal{A}\otimes \cal{A}$, the
countit    $\e\colon \cal{A}\rightarrow \Bbb{C}$
and     the     antipode
$\k\colon \cal{A}\rightarrow\cal{A}$. The result of an $(n-1)$-fold
comultiplication of an element $a\in\cal{A}$ will be symbolically
written  as  $a^{(1)}\otimes \cdots\otimes   a^{(n)}$.  The  adjoint
action of $G$ on itself will be  denoted
by $\ad\colon \cal{A}\rightarrow\cal{A}\otimes \cal{A}$. Explicitly, this
map is given by
\begin{equation}\label{ad}
\ad(a)=a^{(2)}\otimes  \k(a^{(1)})a^{(3)}.
\end{equation}

     Let $\Gamma$ be a first-order differential calculus \cite{W3}
over $G$ and let
    $$\Gamma^{\wedge}=\sideset{}{^\oplus}\sum_{k\geq 0}\Gamma^{\wedge k}$$
be the universal differential envelope (\cite{D}--Appendix B) of $\Gamma$.
Here, the space $\Gamma^{\wedge k}$
consists of $k$-th order elements.

For each
$k\geq 0$ let $p_k \colon \Gamma^{\wedge}\rightarrow \Gamma^{\wedge k}$
be the corresponding projection (we shall use the same symbols for
projection operators associated to an arbitrary graded
algebra built over $\Gamma$).
Further, let
$$\Gamma^{\otimes }=\sideset{}{^\oplus}\sum_{k\geq 0}
\Gamma^{\otimes  k}$$
be the tensor bundle algebra over $\Gamma$. Here,
$$\Gamma^{\otimes k}=\Gamma
\otimes_{\cal A}\dotsb{}_{\cal A}
\otimes\Gamma$$
is the tensor product over $\cal{A}$ of $k$-copies of $\Gamma$.
The algebra $\Gamma^{\wedge}$ can be obtained from $\Gamma^{\otimes}$
by factorizing through the ideal $S^{\wedge}$ generated by elements $Q\in
\Gamma^{\otimes 2}$ of the form
$$Q=\sum_i da_i\otimes_{\cal A}db_i$$
where $a_i,b_i\in\cal A$ satisfy
$$\sum_i a_idb_i=0.$$

     Let  us  assume  that  $\Gamma$  is  left-covariant.  Let
$\ell_\Gamma\colon \Gamma\rightarrow\cal{A}\otimes \Gamma$ be
the corresponding left action of $G$ on $\Gamma$. We shall  denote
by  $\Gamma_{\inv}$  the  space  of  left-invariant   elements   of
$\Gamma$. In other words
$$\Gamma_{\inv}=\bigl\{\,\vartheta\in\Gamma\colon \ell_\Gamma(\vartheta)
=1\otimes\vartheta\,\bigr\}.$$
Further, $\cal{R}\subseteq \ker(\e)$  will  be  the  right
$\cal{A}$-ideal which canonically, in the sense of \cite{W3},
corresponds to $\Gamma$. The map
$\pi\colon \cal{A}\rightarrow\Gamma_{\inv}$ given by
\begin{equation}\label{pi}
\pi(a)=\k(a^{(1)})da^{(2)}
\end{equation}
is surjective and $\ker(\pi)=\Bbb{C}\oplus{\cal R}$. Because  of  this
there exists a natural isomorphism
          $$\Gamma_{\inv}\leftrightarrow\ker(\e)/\cal{R}.$$
The above isomorphism induces a right $\cal{A}$-module structure on
$\Gamma_{\inv}$, which will be denoted by $\circ $. Explicitly,
\begin{equation}\label{mod}
\pi (a)\circ  b=\pi \bigl(ab-\e(a)b\bigr),
\end{equation}
for each $a,b\in\cal{A}$.
The maps $\phi$ and $\ell_\Gamma$ admit common extensions to
homomorphisms $\ell_\Gamma^{\wedge}\colon \Gamma^{\wedge}
\rightarrow\cal{A}\otimes \Gamma^{\wedge}$ and
$\ell_\Gamma^{\otimes}\colon \Gamma^{\otimes
}\rightarrow\cal{A}\otimes
\Gamma^{\otimes }$ (left actions of $G$ on $\Gamma^{\wedge}$ and
$\Gamma^{\otimes }$).

     The tensor product of $k$-copies of $\Gamma_{\inv}$ will  be
denoted  by $\Gamma^{\otimes  k}_{\inv}$. The tensor algebra over
$\Gamma_{\inv}$ will be denoted by $\Gamma_{\inv}^{\otimes }$. It is
naturally isomorphic to the space of left-invariant elements of
$\Gamma^{\otimes }$.

     The subalgebra of  left-invariant  elements  of
$\Gamma^{\wedge}$ will be denoted by $\Gamma^{\wedge}_{\inv}$. This
algebra is naturally graded. We shall denote by
$\Gamma^{\wedge k}_{\inv}$ the space of left-invariant
$k$-th order elements. Let
$\pi_{\inv}\colon \Gamma^{\wedge}\rightarrow\Gamma^{\wedge}_{\inv}$
be the canonical projection map \cite{W3} onto left-invariant elements.
In the framework of the canonical identification $\Gamma^\wedge
\leftrightarrow\cal A\otimes\Gamma^\wedge_{\inv}$ we have
$\pi_{\inv}\leftrightarrow \e\otimes \id $.

     The following natural isomorphism holds
$$\Gamma^{\wedge}_{\inv}=\Gamma^{\otimes }_{\inv}/S^{\wedge}_{\inv}.$$
Here $S^{\wedge}_{\inv}$ is the ideal in
$\Gamma^{\otimes }_{\inv}$, generated by elements
$q\in\Gamma_{\inv}^{\otimes 2}$ of the form
$$q=\pi(a^{(1)})\otimes \pi(a^{(2)}),$$ where $a\in\cal{R}$. This
space is in fact the left-invariant part of the ideal $S^{\wedge}$.

The right $\cal{A}$-module structure $\circ $ can  be
uniquely extended from $\Gamma_{\inv}$ to
$\Gamma_{\inv}^{\wedge,\otimes}$, such that
\begin{align}
1\circ a&=\e(a)1\label{1mod}\\
(\vartheta\eta)\circ a&=(\vartheta\circ a^{(1)})(\eta\circ
a^{(2)})\label{prodmod}
\end{align}
for each $\vartheta,\eta\in\Gamma_{\inv}^{\otimes ,\wedge}$ and
$a\in\cal{A}$. Explicitly, $\circ$ is given by
\begin{equation}\label{moddef}
\vartheta\circ  a=\k(a^{(1)})\vartheta a^{(2)}.
\end{equation}

     The algebra $\Gamma_{\inv}^{\wedge}\subseteq\Gamma^{\wedge}$
is $d$-invariant. The following identities hold
\begin{gather}
d(\vartheta\circ  a)=d(\vartheta)\circ  a-\pi(a^{(1)})(\vartheta\circ
a^{(2)})+(-1)^{\partial\vartheta}(\vartheta\circ
a^{(1)})\pi(a^{(2)})\label{dmod}\\
d\pi(a)=-\pi(a^{(1)})\pi(a^{(2)}).\label{dpi}
\end{gather}

     If $\Gamma$    is    *-covariant    then
the    *-involution
$*\colon \Gamma\rightarrow\Gamma$ is naturally extendible to
$\Gamma^{\wedge}$ and $\Gamma^{\otimes }$ (such  that $(\vartheta\eta)^*
=(-1)^{\partial\vartheta\partial\eta}\eta^*\vartheta^*$  for   each
$\vartheta,\eta\in\Gamma^{\wedge,\otimes }$ ). The maps
$\ell_\Gamma^{\wedge,\otimes}$
are hermitian, in a natural manner.
Algebras
$\Gamma_{\inv}^{\wedge,\otimes }\subseteq\Gamma^{\wedge,\otimes }$
are *-invariant. We have
\begin{equation}
(\vartheta\circ  a)^*=\vartheta^*\circ  \k(a)^*\label{ mod*}
\end{equation}
for each $a\in\cal{A}$ and
$\vartheta\in\Gamma^{\wedge,\otimes }_{\inv}$.

     Explicitly, the *-involution on $\Gamma_{\inv}$ is determined by
\begin{equation}
\pi(a)^* =-\pi\left[\k(a)^*\right].
\end{equation}

     Let us now assume that the calculus $\Gamma$ is bicovariant, and
let $\rig_\Gamma\colon \Gamma\rightarrow\Gamma\otimes \cal{A}$ be the
right action of $G$  on $\Gamma$. Maps  $\phi$  and
$\rig_\Gamma$   admit
common   extensions    to    homomorphisms
$\rig_\Gamma^{\wedge}$, $\rig_\Gamma^\otimes\colon \Gamma^{\wedge,\otimes
}\rightarrow\Gamma^{\wedge,\otimes}\otimes \cal{A}$
(right actions of $G$ on corresponding algebras).
Let
$\adj \colon \Gamma_{\inv}\rightarrow\Gamma_{\inv}\otimes \cal{A}$
be the adjoint action of $G$ on $\Gamma_{\inv}$. The space
$\Gamma_{\inv}$ is right-invariant, that is
$\rig_\Gamma(\Gamma_{\inv})\subseteq\Gamma_{\inv}\otimes \cal{A}$.
We have $\adj =\rig_\Gamma\restr\Gamma_{\inv}$.
Explicitly,
\begin{equation}\label{adpi}
\adj \pi=(\pi\otimes  \id) \ad.
\end{equation}
We shall denote by $\adj^{\otimes }$,
$\adj^{\otimes  k}$, $\adj^{\wedge}$ and
$\adj^{\wedge k}$ the adjoint actions of $G$ on
the  corresponding  spaces  (coinciding  with  the   corresponding
restrictions of $\rig_\Gamma^{\otimes }$ and
$\rig_\Gamma^{\wedge}$.

     The coproduct map
$\phi\colon \cal{A}\rightarrow\cal{A}\otimes \cal{A}$ admits
the  unique   extension   to   the
homomorphism
$\widehat{\phi}\colon \Gamma^{\wedge}\rightarrow\Gamma^{\wedge}
\grten\Gamma^{\wedge}$ of graded-differential algebras.
We have
\begin{equation}
\widehat{\phi}(\vartheta)=\ell_\Gamma(\vartheta)+\rig_\Gamma
(\vartheta),
\end{equation}
for each $\vartheta\in\Gamma$.  Further,  we  have
$\widehat{\phi}(\Gamma_{\inv}^{\wedge})\subseteq\Gamma_{\inv}^{
\wedge}\grten\Gamma^{\wedge}$.
Let $\widehat{\adj} \colon \Gamma_{\inv}^{\wedge}
\rightarrow\Gamma_{\inv}^{\wedge}\grten\Gamma^{\wedge}$
be the corresponding restriction. Explicitly,
\begin{equation}\label{adex}
\widehat{\adj} (\vartheta)=1\otimes \vartheta
+\adj (\vartheta),
\end{equation}
for each $\vartheta\in\Gamma_{\inv}$.

     If $\Gamma$ is a bicovariant *-calculus then the
maps $\widehat{\phi}$, $\widehat{\adj} $ and
all the introduced adjoint and right actions are hermitian,  in  a
natural manner.

\section{Quantum Principal Bundles $\&$
The Corresponding Differential Calculus}

     The aim of this section  is  twofold.   As  first,  we  shall
define  quantum  principal  bundles,  and   briefly   describe   a
geometrical background for this definition. Then,  an  appropriate
differential calculus  over  quantum  principal  bundles  will  be
introduced and analyzed. In particular, besides the  main  algebra
consisting  of  ``differential  forms''  on  the  bundle,  we  shall
introduce and analyze algebras of ``verticalized'' and ``horizontal''
differential forms. Finally, an algebra representing differential
forms on the base manifold will be defined.

     Let us consider a quantum space $M$, formally represented by  a
(unital) *-algebra $\cal{V}$. At the geometrical level, the elements of
$\cal{V}$ play the role of appropriate ``functions'' on this space.
\begin{defn}
     A {\it quantum principal $G$-bundle}
over  $M$  is  a triplet of the form  $P=(\cal{B},i,F)$,
where $\cal{B}$  is  a
(unital) *-algebra while $i\colon \cal{V}\rightarrow\cal{B}$ and
$F\colon \cal{B}\rightarrow\cal{B}\otimes \cal{A}$ are unital
*-homomorphisms such that

\smallskip
({\it qpb1}\/) The map $i\colon \cal{V}\rightarrow\cal{B}$ is injective and
$$ b\in i(\cal{V})\quad \Leftrightarrow\quad F(b)=b\otimes 1$$
for each $b\in\cal{B}$.

\smallskip
({\it qpb2}\/) The following identities hold
\begin{gather}
\id= (\id \otimes  \e)F\label{ide}\\
(\id \otimes \phi)F=(F\otimes  \id) F.\label{coas}
\end{gather}

\smallskip
({\it qpb3}\/) A linear map
$X\colon \cal{B}\otimes \cal{B}\rightarrow\cal{B}\otimes \cal{A}$ defined by
$$X(q\otimes  b)=qF(b)$$
is surjective.
\end{defn}

    The elements of $\cal{B}$ are interpretable as appropriate
``functions''
on the quantum space $P$. The map $F$ plays the role of  the  dualized
right  action  of  $G$  on  $P$.   Condition   ({\it qpb2\/})
justifies   this
interpretation. The map $i\colon \cal{V}\rightarrow\cal{B}$
plays the role of the dualized
``projection'' of $P$ on the base manifold $M$. Condition ({\it qpb1\/})
says  that
$M$ is identificable with the corresponding ``orbit  space''  for  the
right action. Accordingly, the elements of $\cal{V}$ will be identified
with their images in $i(\cal{V})$.

Finally, condition ({\it qpb3\/})  is  an  effective  quantum
counterpart of the classical requirement that the action of $G$
on $P$ is free. It is easy to see that this condition can be
equivalently formulated as

\smallskip
({\it qpb4}\/) For  each  $a\in\cal{A}$  there  exist
elements $b_k,q_k\in\cal{B}$ such that
\begin{equation}
1\otimes  a=\sum_k q_kF(b_k).\label{free}
\end{equation}

     We now pass to questions related to differential calculus on
quantum principal bundles.
Let $P=(\cal{B},i,F)$ be  a
quantum  principal
$G$-bundle over $M$. Let us fix a bicovariant  first-order  *-calculus
$\Gamma$  over $G$  and  let  us  consider  a
graded vector space $\ver(P)=\cal{B}\otimes \Gamma_{\inv}^{\wedge}$
(the grading is induced from $\Gamma_{\inv}^{\wedge}$).

\begin{lem}\bla{i} The formulas
\begin{gather}
(q\otimes   \eta)(b\otimes   \vartheta)=\sum_k    qb_k\otimes
(\eta\circ
c_k)\vartheta \label{verprod}\\
(b\otimes \vartheta)^* =\sum_k b_k^* \otimes (\vartheta^*\circ
c_k^* ) \label{ver*}\\
d_v(b\otimes \vartheta)=b\otimes  d\vartheta +
\sum_kb\otimes \pi(c_k)\vartheta \label{verd}
\end{gather}
where $F(b)=\Sum_k b_k \otimes  c_k$, determine the structure
of a graded-differential *-algebra on $\ver(P)$.

\smallskip
    \bla{ii} As a diferential algebra, $\ver(P)$ is generated by
$\cal{B}=\ver^0(P)$.
\end{lem}

\begin{pf}
Let us first check the associativity of the introduced
product. Applying \eqref{verprod} and \eqref{prodmod} and
elementary properties of $F$ we obtain
\begin{equation*}
\begin{split}
\bigl[(f\otimes \zeta)(b\otimes \vartheta)\bigr](q\otimes \eta)&=\sum_k
\bigl(fb_k\otimes (\zeta\circ  c_k)\vartheta\bigr)(q\otimes \eta)\\
&=\sum_{kl}fb_kq_l\otimes\bigl[\bigl((\zeta\circ  c_k)\vartheta\bigr)\circ
d_l\bigr]\eta\\
&= \sum_{kl}fb_kq_l\otimes  \bigl(\zeta\circ (c_kd_l^{(1)})\bigr)(
\vartheta\circ  d_l^{(2)})\eta\\
&=(f\otimes\zeta)\sum_{l}b
q_l\otimes (\vartheta\circ d_l)\eta\\
&=(f\otimes \zeta)\bigl[(b\otimes
\vartheta)(q\otimes \eta)\bigr],
\end{split}
\end{equation*}
where $F(q)=\Sum_lq_l\otimes  d_l$.

     Evidently, $\ver(P)$ is a unital algebra, with  the  unity
$1\otimes 1$.
Now we prove that \eqref{ver*} determines a *-algebra
structure on $\ver(P)$. We have
\begin{multline*}
\bigl[(b\otimes \vartheta)^*\bigr]^*=\sum_k\bigl(b_k^*\otimes
(\vartheta^*\circ  c_k^
*)\bigr)^*=\sum_kb_k\otimes(\vartheta^*\circ c_k^{(2)*})^*\circ c_k^{(1)}\\
=\sum_kb_k\otimes\bigl(\vartheta\circ
\k^{-1}(c_k^{(2)})\bigr)\circ  c_k^{(1)}=
\sum_kb_k\otimes \vartheta\circ\bigl(\k^{-1}(c_k^{(2)})c_k^{(1)}\bigr)=b
\otimes\vartheta.
\end{multline*}
Thus, $*$ is involutive. Further,
\begin{equation*}
\begin{split}
\bigl[(q\otimes \eta)(b\otimes \vartheta)\bigr]^*&=\sum_k\bigl
(qb_k\otimes(\eta\circ c_k)\vartheta\bigr)^*\\
&=\sum_{kl}(q_lb_k)^*\otimes \bigl((\eta\circ  c_k^{(2)})\vartheta\bigr)^*
\circ (d_lc_k^{(1)})^*\\
&=(-1)^{\partial\vartheta\partial\eta}\sum_{kl}b_k^*q_l^*\otimes
\bigl[\vartheta^*\circ  c_k^{(1)*}d_l^{(1)*}\bigr]\bigl[\eta^*\circ
\bigl(\k(c_k^{(3)})^*c_k^
{(2)*}d_l^{(2)*}\bigr)\bigr]\\
&=(-1)^{\partial\vartheta\partial\eta}\sum_{kl}b_k^*q_l^*\otimes
 (\vartheta^*\circ c_k^*d_l^{(1)*})(\eta^*\circ  d_l^{(2)*})\\
&=(-1)^{\partial\vartheta\partial\eta}\sum_{kl}\bigl[b_k^*\otimes
(\vartheta^*\circ  c_k^*)\bigr]\bigl[q_l^*\otimes(\eta^*\circ d_l^*)\bigr]\\
&=(-1)^{\partial\vartheta\partial\eta}(b\otimes \vartheta)^*(q
\otimes\eta)^*.
\end{split}
\end{equation*}

     Let us check that \eqref{verd} defines a
hermitian  differential on the *-algebra $\ver(P)$. We compute
\begin{equation*}
\begin{split}
d_v\bigl[(q\otimes \eta)(b\otimes\vartheta)\bigr]&=\sum_kd_v\bigl(qb_k
\otimes (\eta\circ
c_k)\vartheta\bigr)=\sum_kqb_k\otimes
d\bigl((\eta\circ c_k)\vartheta\bigr)\\
&\phantom{=}+\sum_{kl}q_lb_k\otimes \pi(d_lc_k^{(1)})(\eta\circ
c_k^{(2)})\vartheta\\
&=\sum_kqb_k\otimes (d(\eta)\circ
c_k)\vartheta-\sum_kqb_k\otimes \pi(c_k^{(1)})(\eta\circ c_k^{(2)})
\vartheta\\
&\phantom{=}+(-1)^{\partial\eta}\sum_k qb_k\otimes
\left((\eta\circ c_k^{(1)})\pi
(c_k^{(2)})\vartheta
+(\eta\circ  c_k)d\vartheta\right)\\
&\phantom{=}+\sum_{kl}q_lb_k\otimes \bigl(\pi(d_l)\circ
c_k^{(1)}\bigr)(\eta\circ c_k^{(2)})\vartheta\\
&\phantom{=}+\sum_kqb_k\otimes \pi(c_k^{(1)})(\eta\circ
c_k^{(2)})\vartheta\\
&=(q\otimes d\eta)(b\otimes
\vartheta)+(-1)^{\partial\eta}(q\otimes
\eta)\sum_kb_k\otimes \pi(c_k)\vartheta\\
&\phantom{=}+(-1)^{\partial\eta}(q\otimes \eta)(b\otimes d\vartheta)
+\sum_l\bigl(q_l\otimes \pi(d_l)\eta\bigr)(b\otimes
\vartheta)\\
&=\bigl[d_v(q\otimes
\eta)\bigr](b\otimes \vartheta)+(-1)^{\partial\eta}(q\otimes \eta)d_v(b
\otimes\vartheta).
\end{split}
\end{equation*}
Here, we have used \eqref{mod}, \eqref{prodmod}, \eqref{dmod} and the main
properties  of $F$.

     Further
\begin{multline*}
d_v^2(b\otimes\vartheta)=d_v\Bigl(b\otimes
d\vartheta+\sum_kb_k\otimes
\pi(c_k)\vartheta\Bigr)=\sum_kb_k\otimes \pi(c_k)d\vartheta\\
+\sum_kb_k\otimes \pi(c_k^{(1)})\pi(c_k^{(2)})\vartheta+\sum_kb_k
\otimes d\bigl(\pi(c_k)\vartheta\bigr)=0,
\end{multline*}
according to \eqref{dpi}. Finally,
\begin{equation*}
\begin{split}
d_v\bigl[(b\otimes \vartheta)^*\bigr]&=
\sum_kb_k^*\otimes d(\vartheta^*\circ c_k^*)
+\sum_kb_k^*\otimes \pi(c_k^{(1)*})(\vartheta^*\circ
c_k^{(2)*})\\
&=\sum_kb_k^*\otimes (d\vartheta^*)\circ c_k^*
-\sum_kb_k^*\otimes \pi(c_k^{(1)*})(\vartheta^*\circ c_k^{(2)*})\\
&\phantom{=}+(-1)^{\partial\vartheta}
\sum_k b_k^*\otimes (\vartheta^*\circ c_k^{(1)*})\pi(c_k^{(2)*})
+\sum_k b_k^*\otimes\pi(c_k^{(1)*})(\vartheta^*\circ
c_k^{(2)*}) \\
&=(b\otimes
d\vartheta)^*+(-1)^{\partial\vartheta}\sum_kb_k^*\otimes (\vartheta
^*\circ  c_k^{(1)*})\pi(c_k^{(2)*})\\
&=(b\otimes d\vartheta)^*-(-1)^{\partial\vartheta}
\sum_kb_k^*\otimes (\vartheta
^*\circ  c_k^{(1)*})\bigl[\pi\bigl(\k(c_k^{(3)})^*\bigr)\circ
c_k^{(2)*}\bigr]\\
&=(b\otimes d\vartheta)^*+
(-1)^{\partial\vartheta}\sum_kb_k^*\otimes\bigl(\vartheta
^*\pi(c_k^{(2)})^*\bigr)\circ c_k^{(1)*}\\
&=(b\otimes d\vartheta)^*+
\sum_kb_k^*\otimes \bigl(
\pi(c_k^{(2)})\vartheta\bigr)^*\circ c_k^{(1)*}=
\bigl[d_v(b\otimes \vartheta)\bigr]^*.
\end{split}
\end{equation*}

     Hence, $d_v$  is a hermitian differential. To prove  ({\it ii\/})
it  is
sufficient to check that elements of  the  form $qd_v(b)$  linearly
generate $\ver^1(P)$. However, this  directly  follows  from  property
({\it qpb4}\/) in the definition of quantum principal bundles.
\end{pf}

     In the following it will be assumed that $\ver(P)$ is
endowed  with
the  constructed  graded-differential *-algebra  structure.  The
elements of $\ver(P)$ are interpretable as verticalized  differential
forms on the bundle. In classical  geometry,  these  entities  are
obtained by restricting the domain of differential  forms (on  the
bundle) to the Lie algebra of vertical vector fields.
\begin{lem}
There exists the unique homomorphism
$$ \widehat{F}_v\colon \ver(P)\rightarrow
\ver(P)\grten\Gamma^{\wedge}$$
of (graded) differential algebras extending the map $F$. We have
\begin{equation}\label{vFcoas}
(\widehat{F}_v\otimes
\id) \widehat{F}_v=(\id \otimes \widehat{\phi})\widehat{F}_v.
\end{equation}
The map $\widehat{F}_v$  is hermitian, in the sense that
\begin{equation}\label{vF*}
\widehat{F}_v*=(*\otimes *)\widehat{F}_v.
\end{equation}
\end{lem}
\begin{pf}
According to ({\it ii\/}) of the previous lemma, the map
$\widehat{F}_v$ is unique, if exists.
Let us {\it define} a linear map
$\widehat{F}_v\colon \ver(P)\rightarrow
\ver(P)\grten\Gamma^{\wedge}$ by
$$\widehat{F}_v(b\otimes \vartheta)=\sum_{kl}b_k\otimes \vartheta_l
\otimes  c_kw_l$$
where $F(b)=\Sum_kb_k\otimes c_k$ and
$\widehat{\adj} (\vartheta)=\Sum_l\vartheta_l\otimes  w_l$.
It is easy to see that such
defined map $\widehat{F}_v$ is a differential algebra homomorphism.

    Identities \eqref{vFcoas} and \eqref{vF*}
directly follow form the fact  that
$\ver(P)$ is generated by $\cal{B}$, as well as from property \eqref{coas}
and the hermicity of $d_v$ respectively.
\end{pf}

    Let us consider a *-homomorphism
$F_v\colon \ver(P)\rightarrow \ver(P)\otimes \cal{A}$ given by
$$
 F_v=(\id \otimes  p_0)\widehat{F}_v. $$
This map  extends  the  action  $F$. It is interpretable as the
(dualized) right action of $G$ on verticalized forms. The  following
identities hold:
\begin{gather}
(F_v\otimes  \id) F_v =(\id \otimes \phi)F_v\\
(\id \otimes \e)F_v=\id \\
(d_v\otimes  \id) F_v=F_vd_v.
\end{gather}

    The first two identities justify the interpretation of $F_v$   as
an action of $G$. The last identity says that the differential
$d_v$ is right-covariant.

    So far about verticalized differential forms. We shall  assume
that a {\it complete} differential calculus over the bundle $P$  is  based
on a graded-differential *-algebra $\Omega(P)$ possessing  the following
properties

\smallskip
({\it diff1\/}) As a differential algebra, $\Omega(P)$ is generated by
$\cal{B}=\Omega^0(P)$.

\smallskip
({\it diff2\/}) The map $F\colon \cal{B}\rightarrow \cal{B}\otimes \cal{A}$
is extendible to a homomorphism
$$
\widehat{F}\colon \Omega(P)\rightarrow\Omega(P)\grten\Gamma^{
\wedge}$$
of graded-differential algebras.

    Let us fix a graded-differential  *-algebra $\Omega(P)$
such  that the above
properties hold. The elements of $\Omega(P)$ will  play  the  role
of differential forms on $P$. It is easy to see that  the
map $\widehat{F}$ is uniquely determined.

\begin{lem}We have
\begin{align}
(\widehat{F}\otimes
\id) \widehat{F}&=(\id \otimes \widehat{\phi})\widehat{F}
\label{dvFcoas}\\
\widehat{F}*&=(*\otimes *)\widehat{F}.\label{dvF*}
\end{align}
\end{lem}
\begin{pf}
 Both identities directly follow from similar properties
of $F$, and from properties ({\it diff1/2\/}).
\end{pf}

    The formula
\begin{equation}
F^{\wedge}=(\id \otimes  p_0)\widehat{F}\label{defwF}
\end{equation}
defines a *-homomorphism
$F^{\wedge}\colon \Omega(P)\rightarrow\Omega(P)\otimes \cal{A}$  extending
the action $F$. The following identities hold:
\begin{align}
(F^{\wedge}\otimes \id) F^{\wedge}&=(\id \otimes \phi)F^{\wedge}
\label{wFcoas}\\
(\id \otimes  \e)F^{\wedge}&=\id \label{ewF} \\
(d\otimes  \id) F^{\wedge}&=F^{\wedge}d.\label{dwF}
\end{align}

The map $F^{\wedge}$  is interpretable as the (dualized) right action of
$G$ on differential forms. As a homomorphism between algebras,
$F^{\wedge}$  is
completely determined by \eqref{dwF}, and by the fact that it extends $F$.

     Now, a very important algebra representing horizontal forms
will be introduced in the game. Intuitively  speaking,  horizontal
forms can be characterized as those elements  of $\Omega(P)$  possessing
``trivial'' differential properties along vertical fibers.
\begin{defn}
The elements of the graded *-subalgebra
$$\hor(P)=\widehat{F}^{-1}\bigl(\Omega(P)\otimes\cal{A}\bigr) $$
of $\Omega(P)$ are called {\it horizontal forms}.
\end{defn}

     Evidently, $\cal{B}=\hor^0(P)$.
\begin{lem}
The algebra $\hor(P)$ is
$F^{\wedge}$-invariant. In other words
\begin{equation}
F^{\wedge}\bigl(\hor(P)\bigr)\subseteq \hor(P)\otimes\cal{A}.\label{defhor}
\end{equation}
\end{lem}
\begin{pf}
If $\varphi\in \hor(P)$  then
$(\widehat{F}\otimes
\id) \widehat{F}(\varphi)=(\id \otimes \widehat{\phi})\widehat{F}(\varphi)
=(\id \otimes\phi)F^{\wedge}(\varphi)$
belongs to $\Omega(P)\otimes \cal{A}\otimes \cal{A}$. Hence,
$\widehat{F}(\varphi)=F^{\wedge}(\varphi)\in \hor(P)\otimes \cal{A}$.
\end{pf}

     The  following  technical  lemma  will  be  helpful  in various
considerations.
\begin{lem}\label{lem:tex}
Let us consider a homogeneous element
$w\in\Omega^n(P)$. Let $0\leq k\leq n$ be an integer such that
$(\id \otimes p_l)\widehat{F}(w)=0$ for  each $l>k$.  Then
there exist  horizontal  forms  $\varphi_1,\ldots,\varphi_m\in
\hor^{n-k}(P)$  and  elements
$\vartheta_1,\ldots,\vartheta_m\in\Gamma_{\inv}^{\wedge k}$ such that
\begin{equation}\label{tex1}
 (\id \otimes
p_k)\widehat{F}(w)=\sum_{i=1}^mF^{\wedge}(\varphi_i)\vartheta_i.
\end{equation}
\end{lem}
\begin{pf}
The statement is non-trivial only if $(\id \otimes
p_k)\widehat{F}(w)\neq0$.

     We have
\begin{equation}
(\id \otimes  p_k)\widehat{F}(w)=\sum_{ij}\xi_{ij}\otimes
a_{ij}\vartheta_i,\label{tex2}
\end{equation}
where $\vartheta_i\in\Gamma_{\inv}^{\wedge k}$ are some
some linearly independent  elements,
$\xi_{ij}\in\Omega^{n-k}(P)$ and $a_{ij}\in\cal{A}$. Applying
\eqref{adex}, \eqref{dvFcoas} and \eqref{tex2}, and the definition  of
$k$ we find
\begin{multline*}
\sum_{ij}\widehat{F}(\xi_{ij})\otimes
a_{ij}\vartheta_i=(\widehat{F}\otimes
p_k)\widehat{F}(w)=(\id^2\otimes
p_k)\sum_{ij}\xi_{ij}\otimes \widehat{\phi}(a_{ij}\vartheta_i)\\
=\sum_{ij}\xi_{ij}\otimes  a_{ij}^{(1)}\otimes
a_{ij}^{(2)}\vartheta_i.
\end{multline*}

Acting by $\id^2\otimes \pi_{\inv}$    on both sides of  the
above equality we obtain
\begin{equation}
\sum_i\widehat{F}(\varphi_i)\otimes \vartheta_i
=\sum_{ij}\xi_{ij}\otimes  a_{ij}\otimes \vartheta_i,
\label{tex3}
\end{equation}
where
$\varphi_i =\Sum_j\xi_{ij}\e(a_{ij})$. In other words
\begin{equation}\label{tex4}
\widehat{F}(\varphi_i)=\sum_j\xi_{ij}\otimes  a_{ij},
\end{equation}
and in particular $\varphi_i\in \hor^{n-k}(P).$
Finally, combining \eqref{tex2} and \eqref{tex4} we conclude  that
\eqref{tex1} holds.
\end{pf}

     We  are  going  to  construct  a  quantum  analog   for   the
``verticalizing'' homomorphism (in classical geometry,  induced  by
restricting the domain of differential forms  on  vertical  vector
fields on the bundle).
\begin{pro}
    There exists the unique (graded) differential
algebra homomorphism $\pi_v\colon
\Omega(P)\rightarrow \ver(P)$ reducing to the identity map on $\cal{B}$.
The map $\pi_v$  is surjective and hermitian. Moreover,
\begin{align}
\widehat{F}_v\pi_v&=(\pi_v\otimes  \id) \widehat{F}\label{pivdF}\\
F_v\pi_v&=(\pi_v\otimes \id) F^{\wedge}.\label{pivF}
\end{align}
\end{pro}
\begin{pf} Let  us  {\it define}  a   linear   (grade-preserving)   map
$\pi_v\colon \Omega(P)\rightarrow \ver(P)$ by requiring
$$\pi_v(w)=(\id \otimes \pi_{\inv}p_k)\widehat{F}(w)$$
for each $w\in\Omega^k(P)$. Obviously, $\pi_v$  is reduced to the
identity on $\cal{B}$.

    Let us prove that $\pi_v$ is a differential  algebra  homomorphism.
For given forms $w\in\Omega^k(P)$ and $u\in\Omega^l(P)$ let us choose
elements $b_i,q_j\in\cal{B}$ and
$\vartheta_i,\eta_j\in\Gamma^{\wedge k,l}_{\inv}$ such that
\begin{equation*}
(\id \otimes  p_k)\widehat{F}(w)=\sum_iF(b_i)\vartheta_i \qquad
(\id \otimes  p_l)\widehat{F}(u)=\sum_jF(q_j)\eta_j,
\end{equation*}
in accordance with the previous lemma.

    We have then
\begin{equation*}
\pi_v(w)=\sum_ib_i\otimes \vartheta_i\qquad
\pi_v(u)=\sum_jq_j\otimes \eta_j.
\end{equation*}

    A direct computation now gives
\begin{equation*}
\begin{split}
\pi_v(wu)&=(\id \otimes \pi_{\inv}p_{k+l})\widehat{F}(wu)\\
&=\sum_{ij}
(\id \otimes \pi_{\inv})\bigl[\bigl(F(b_i)\vartheta_i
\bigr)\bigl(F(q_j)\eta_j\bigr)\bigr]\\
&=\sum_{ij}(\id \otimes \pi_{\inv})\sum_{rs}\bigl[b_{ir}q_{js}\otimes
c_{ir}d_{js}^{(1)}(\vartheta_i\circ  d_{js}^{(2)})\eta_j\bigr]\\
&=\sum_{ijs}b_iq_{js}\otimes (\vartheta_i\circ
d_{js})\eta_j=\biggl[\sum_ib_i\otimes \vartheta_i
\biggr]\biggl[\sum_jq_j\otimes \eta_j\biggr]\\
&=\pi_v(w)\pi_v(u),
\end{split}
\end{equation*}
where $F(b_i)=\Sum_rb_{ir}\otimes  c_{ir}$   and
$F(q_j)=\Sum_sq_{js}\otimes  d_{js}$. Further,
\begin{equation*}
\begin{split}
\pi_vd(w)&=(\id \otimes \pi_{\inv}p_{k+1})d\widehat{F}(w)=(\id \otimes
\pi_{\inv}d)\biggl[\sum_iF(b_i)\vartheta_i\biggr]\\
&=(\id \otimes \pi_{\inv})\biggl[\sum_{ir}b_{ir}\otimes  c_{ir}d\vartheta_i
+b_{ir}\otimes d(c_{ir})\vartheta_i\biggr]\\
&=\sum_ib_i\otimes d\vartheta_i+\sum_{ir}b_{ir}\otimes \pi(c_{ir})
\vartheta_i\\
&=d_v\biggl[\sum_ib_i\otimes \vartheta_i\biggr]=d_v\pi_v(w).
\end{split}
\end{equation*}

     Consequently, $\pi_v$ is a homomorphism of  differential algebras.
The map $\pi_v$ is hermitian, because differentials on $\Omega(P)$
and $\ver(P)$
are hermitian, and the differential algebra $\Omega(P)$ is
generated by $\cal{B}$. To prove \eqref{pivdF} it is sufficient to
observe that  its  both  sides
are differential algebra homomorphisms coinciding with  $F$  on
$\cal{B}$. Finally, \eqref{pivF} follows from \eqref{pivdF},
and definitions of
$F^{\wedge}$  and $F_v$ .
\end{pf}

     Let us now consider a sequence
\begin{equation}\label{seq}
0\rightarrow \hor^1(P)\hookrightarrow \Omega^1(P) @>{\pi_v}>> \ver^1(P)
\rightarrow 0
\end{equation}
of natural homomorphisms of *-$\cal B$-bimodules.
\begin{lem}\label{lem:seq}
The above sequence is exact.
\end{lem}
\begin{pf}
Clearly, $\hor^1(P)\subseteq
\ker(\pi_v)\cap\Omega^1(P)$ and $\pi_v\bigl(\Omega^1(P)\bigr)=\ver^1(P)$.
For each $w\in\Omega^1(P)$ we have
$$
(\id \otimes p_1)\widehat{F}(w)=\sum_iF(b_i)\vartheta_i, $$
for some $b_i\in\cal{B}$ and $\vartheta_i\in\Gamma_{\inv}$,
according to Lemma~\ref{lem:tex}. This implies
$$
\pi_v(w)=\sum_ib_i\otimes \vartheta_i. $$
Consequently if $w\in \ker(\pi_v)$ then $(\id \otimes
p_1)\widehat{F}(w)=0$, and hence $w\in \hor^1(P)$.
\end{pf}

     If a differential calculus on the bundle $P$ is given, then  it
is possible to construct a natural differential  calculus  on  the
base space $M$. This calculus is  based  on  a  graded-differential
*-subalgebra $\Omega(M)\subseteq\Omega(P)$ consisting  of
right-invariant horizontal forms. Equivalently,

$$\Omega(M)=\left\{w\in\Omega(P)\colon \widehat{F}(w)=w\otimes
1\right\}.$$

In a special case when the group $G$ is ``connected'' in the sense that
only scalar elements of $\cal{A}$ are anihilated by the differential
map, the algebra $\Omega(M)$ can be described as
$$\Omega(M)=d^{-1}\bigl(\hor(P)\bigr)\cap\hor(P).$$

The differential algebra $\Omega(M)$ is generally not generated
by its $0$-order subalgebra
$\Omega^0(M)=\cal{V}$, in contrast to $\Omega(P)$. However, property
({\it diff1\/}) is not essential for developing a large part of the
formalism. Furthermore, it turns out that the algebra $\hor(P)$
is generally not generated by $\cal{B}$ and $\hor^1(P)$.

\section{The Formalism of Connections}
     In this section a general theory of connections on quantum  principal
bundles will be  presented. As first, quantum
analogs of pseudotensorial forms will be defined.

     Let  $V$  be  a  vector  space,   and   let   $v\colon V\rightarrow
V\otimes\cal{A}$ be a representation
of $G$ in this space. Let $\psi (v,P)$ be the space  of  all  linear  maps
$\zeta\colon V\rightarrow\Omega(P)$ such that the diagram
\begin{equation}\label{defpse}
\begin{CD}
V @>{\mbox{$\zeta$}}>> \Omega(P)\\
@V{\mbox{$v$}}VV @VV{\mbox{$F^{\wedge}$}}V\\
V\otimes\cal{A} @>>{\mbox{$\zeta\otimes \id $}}> \Omega(P)\otimes\cal{A}
\end{CD}
\end{equation}
is commutative (intertwiners between $v$ and $F^\wedge$).

The space $\psi (v,P)$ is naturally graded
\begin{equation}
\psi (v,P)=\sideset{}{^\oplus}\sum_{k\geq 0}\psi^k(v,P),
\label{gradpse}
\end{equation}
where $\psi^k(v,P)$ consists of maps with values in $\Omega^k(P)$.
The  elements of $\psi^k(v,P)$ can be interpreted as
pseudotensorial  $k$-forms
on $P$ with
values in the dual space $V^*$ .   Further,   $\psi  (v,P)$   is   closed
with   respect   to
compositions with $d\colon \Omega (P)\rightarrow\Omega (P)$.
It is also a module,
in natural manner,
over the subalgebra consisting of right-invariant forms.
Let
\begin{equation}\label{gradten}
\tau (v,P)=\sideset{}{^\oplus}\sum_{k\geq 0}\tau^k(v,P)
\end{equation}
be the graded subspace of $\psi (v,P)$ consisting of  pseudotensorial
forms having the values in $\hor(P)$. The elements of $\tau (v,P)$  are
interpretable as tensorial forms on $P$ with values in $V^*$. The space
$\tau (v,P)$ is a module over $\Omega(M)$.

If the space $\Gamma_{\inv}$ is infinite-dimensional and if
$\Omega(P)$ possesses infinitely non-zero components then, generally,
sums figuring in \eqref{gradpse} and \eqref{gradten}
will be ``larger'' then
standard direct sums of spaces (and should be interpreted
in the appropriate way).

Forms   $\varphi\in\tau
(v,P)$
can be also defined by the following equality
\begin{equation}
\widehat{F}\varphi (\vartheta  )=(\varphi  \otimes  \id) v(\vartheta
).\label{deften}
\end{equation}

 If  the  space  $V$   is  endowed  with  an  antilinear  involution
$*\colon V\rightarrow V$ such that $v*=(*\otimes *)v$ then the formula
\begin{equation}\label{*pse}
\varphi^* (\vartheta )=\varphi (\vartheta^*)^*
\end{equation}
determines natural *-involutions on $\psi (v,P)$
and $\tau (v,P)$.

 For the purposes of this paper, the most  important  is  the  case
$V=\Gamma_{\inv}$,  and  $v=\adj $.   In   this   case  we
shall write $\psi (P)=\psi (\adj ,P)$ and $\tau (P)=\tau
(\adj ,P)$.

     We pass to the definition of connection forms.
\begin{defn}
{\it A  connection}  on  $P$   is   a   hermitian   map
$\omega\in\psi^1 (P)$ satisfying
\begin{equation}\label{pivcon}
\pi_v \omega(\vartheta)=1\otimes\vartheta
\end{equation}
for each $\vartheta\in\Gamma_{\inv}$.
\end{defn}

     The above condition corresponds to the classical
requirement that connection forms (understood as $\lie(G)$-valued
pseudotensorial $\ad$-type 1-forms) map fundamental
vector fields into their generators.

\begin{thm}
Every  quantum  principal  bundle  $P$  admits  at
least one connection.
\end{thm}
\begin{pf}
Let us consider the space
$W=\pi_v^{-1}(\Gamma_{\inv})\cap\Omega^1(P)$. By definition
this space is invariant under $*$ and $F^\wedge$, and
$\pi_v(W)=\Gamma_{\inv}$. We can write
$$W=\sideset{}{^\oplus}\sum_{\alpha\in\cal{T}}W^\alpha$$
where $W^\alpha$ are corresponding multiple irreducible subspaces
(the notation is explained in Appendix B). Similarly
$$\Gamma_{\inv}=\sideset{}{^\oplus}\sum_{\alpha\in\cal{T}}
\Gamma_{\inv}^\alpha.$$
The following decompositions hold
$$
W^\alpha\leftrightarrow\mbox{Mor}(u^\alpha,F^\wedge)
\otimes\Bbb{C}^n\qquad
\Gamma_{\inv}^\alpha\leftrightarrow\mbox{Mor}(u^\alpha,
\adj )\otimes\Bbb{C}^n
$$
where $n$ is the dimension of $\alpha$.
Further, $\pi_v(W^\alpha)=\Gamma_{\inv}^\alpha$ for each
$\alpha\in\cal{T}$.

In terms of the above identifications the restriction map $\pi_v\colon
W^\alpha\rightarrow\Gamma_{\inv}^\alpha$ is given by
$$\pi_v\bigl\{\mu\otimes x\bigr\}=\pi_v\mu\otimes x.$$
This map is surjective. Let
$\tau^\alpha\colon \mbox{Mor}(u^\alpha,\adj )\rightarrow
\mbox{Mor}(u^\alpha,F^\wedge)$ be a left inverse of $\pi_v\restr
W^\alpha$.

Let $\tau\colon \Gamma_{\inv}
\rightarrow W$ be a map defined by
$$\tau=\sideset{}{^\oplus}\sum_{\alpha\in\cal{T}}\tau_\alpha$$
where $\tau_\alpha\colon \Gamma_{\inv}^\alpha
\rightarrow W^\alpha$ are given by
$\tau_\alpha=\tau^\alpha\otimes \id $. By construction, $\tau$ intertwines
$\adj $
and $F^{\wedge}$ and satisfies
$\pi_v\tau(\vartheta)=\vartheta$ for each
$\vartheta\in\Gamma_{\inv}$. Without a lack of generality we can assume
that $\tau$ is hermitian (if not,  we  can
consider  another  intertwiner $*\tau *\colon \Gamma_{\inv}\rightarrow W$
and redefine $\tau\mapsto(\tau +*\tau *)/2$). Finally,
composing $\tau$ and  the
inclusion map $W\hookrightarrow\Omega(P)$ we obtain  a  connection
on $P$.
\end{pf}

 Let $\con(P)$ be the set of all connections on $P$.  This  is  a  real
affine  subspace  of $\psi^1(P)$.  The  corresponding  vector  space
consists of hermitian tensorial $1$-forms.

     Connections can be described in a  different,  more  concise,
but equivalent manner, using an algebraic  condition  which  is  a
symbiosis of the  verticalization  condition  \eqref{pivcon}
and the pseudotensoriality property.
\begin{lem}
A first-order linear map $\omega\colon \Gamma_{\inv}
\rightarrow \Omega(P)$ is a  connection  on  $P$ iff
\begin{align}
\omega(\vartheta^*)=&\omega(\vartheta )^*\label{*con}\\
\widehat{F}\omega (\vartheta)=&(\omega\otimes
\id) \adj (\vartheta )+1\otimes\vartheta ,\label{dxFcon}
\end{align}
for each $\vartheta\in\Gamma_{\inv}$.
\end{lem}
\begin{pf}
It is clear  that  above  listed  properties  imply  that
$\omega$
is a connection on $P$. Conversely, let us consider an arbitrary
$\omega\in \con(P)$. Then the
pseudotensoriality property and Lemma~\ref{lem:tex} imply  that
for  each $\vartheta\in\Gamma_{\inv}$ we have
$$
\widehat{F}\omega(\vartheta)=(\omega\otimes
\id) \adj (\vartheta )+\sum_k F(b_k)\vartheta_k,$$
for some $b_k\in\cal{B}$ and $\vartheta_k \in\Gamma_{\inv}$.
Properties \eqref{ide}
and \eqref{pivcon},   and   the     definition     of
$\pi_v$    imply
$1\otimes\vartheta=\Sum_k  b_k  \otimes\vartheta_k$.  Hence
\eqref{dxFcon} holds.
\end{pf}

 Every connection $\omega$ canonically gives rise to a splitting  of
the  sequence \eqref{seq},   undestood   as   a   sequence   of   left
$\cal{B}$-modules. Indeed,  the  map  $\mu_{\omega}\colon \ver^1(P)
\rightarrow\Omega^1(P)$ defined by
$$\mu_{\omega}(b\otimes\vartheta )=b\omega(\vartheta )$$
splits the mentioned sequence. Moreover, the  map $\mu_{\omega}$
intertwines the corresponding right actions.

 For each $\omega\in \con(P)$ let
$\omega^{\otimes}\colon \Gamma_{\inv}^{\otimes}\rightarrow\Omega (P)$ be the
corresponding unital multiplicative extension.

      Two particularly interesting classes of connections  naturally
appear  in  deeper  considerations.  These  classes   consist   of
connections possessing some additional  properties  that  will  be
called multiplicativity and regularity.

\begin{defn}
A connection $\omega$ is called
{\it multiplicative} iff
\begin{equation}\label{defmul}
\omega\pi(a^{(1)})\omega\pi(a^{(2)})=0
\end{equation}
for each $a\in\cal{R}$. Equivalently, $\omega$ is multiplicative iff
$\omega^{\otimes}\restr S_{\inv}^{\wedge}=0$.
\end{defn}

If $\omega$ is multiplicative then there exists the unique unital
multiplicative extension
$\omega^\wedge\colon\Gamma_{\inv}^\wedge\rightarrow
\Omega(P)$. This map can be obtained by
by factorizing $\omega^{\otimes}$ through
$S_{\inv}^{\wedge}$ (both $\omega^{\wedge}$ and $\omega^\otimes$
are *-preserving).

     Condition \eqref{defmul} gives a quadratic  constraint  in
the  space $\con(P)$. It is worth  noticing  that  in  the
classical theory all connections are multiplicative.
\begin{defn}
A connection $\omega$ is called {\it
regular} iff
\begin{align}
\omega(\vartheta )\varphi&=(-1)^{\partial\varphi}\sum_k\varphi_k
\omega(\vartheta\circ c_k)\label{defreg1}\\
\intertext{
holds for each $\varphi\in \hor(P)$ and $\vartheta\in\Gamma_{\inv}$.
  Here, $F^{\wedge}(\varphi)=\Sum_k\varphi_k \otimes c_k$.
Equivalently,}
\varphi\omega(\vartheta )&=(-1)^{\partial\varphi}\sum_k
\omega\bigl(\vartheta \circ \k^{-1}(c_k)\bigr)\varphi_k.\label{defreg2}
\end{align}
In particular, regular connections graded-commute with forms
from $\Omega(M)$.
\end{defn}

 Regular connections, if exist, form an affine subspace
$\rc(P)\subseteq \con(P)$. The corresponding vector space
consists of hermitian forms $\zeta\in\tau^1(P)$ satisfying
\begin{align}
\varphi \zeta(\vartheta )&=(-1)^{\partial\varphi}\sum_k
\zeta\bigl(\vartheta
\circ \k^{-1}(c_k)\bigr)\varphi_k,\label{vreg2} \\
\intertext{or equivalently }
\zeta(\vartheta)\varphi&=(-1)^{\partial\varphi}\sum_k \varphi_k
\zeta(\vartheta\circ c_k)\label{vreg1}
\end{align}
for each $\varphi\in \hor(P)$ and $\vartheta\in\Gamma_{\inv}$.

It is important to mention that if $\omega$  is  regular  then  the
corresponding splitting $\mu_{\omega}$ is a splitting of
*-bimodules.

If  the  theory  is  confined  to  bundles   over   classical
manifolds and if $\Gamma$ is
the minimal admissible calculus then  regular
connections  are   precisely   {\it classical   connections} (in   the
terminology of the previous paper). In particular, in the
case of classical principal bundles all connections are regular.

Let us consider an expression
\begin{equation}\label{lacreg}
\lh^{\omega}(\vartheta,\varphi)=\omega(\vartheta)\varphi
-(-1)^{\partial\varphi}\varphi_k\omega(\vartheta\circ c_k)
\end{equation}
where $\varphi\in\hor(P)$ and $\vartheta\in\Gamma_{\inv}$. The map
$\lh^{\omega}$ measures the lack of regularity of $\omega$.
\begin{lem}\label{lem:lac}
We have
\begin{equation}
\widehat{F}\lh^\omega(\vartheta,\varphi)=\sum_{jk} \lh^\omega(
\vartheta_j,\varphi_k)\otimes d_jc_k,
\end{equation}
where $\adj(\vartheta)=\Sum_j\vartheta_j\otimes d_j$. In particular,
$\lh^\omega$ is a $\hor(P)$-valued map.
\end{lem}
\begin{pf} A direct computation gives
\begin{equation*}
\begin{split}
\widehat{F}\lh^\omega(\vartheta,\varphi)&=\sum_{jk}\omega(\vartheta_j)
\varphi_k\otimes d_j c_k+(-1)^{\partial\varphi}\sum_k\varphi_k
\otimes\vartheta c_k\\
&\phantom{=}-(-1)^{\partial\varphi}\sum_{kj}\varphi_k\omega(\vartheta_j\circ
c_k^{(3)})\otimes c_k^{(1)}\k(c_k^{(2)})d_jc_k^{(4)}\\
&\phantom{=}-(-1)^{\partial\varphi}\sum_k\varphi_k\otimes
c_k^{(1)}(\vartheta\circ c_k^{(2)})\\
&=\sum_{jk} \lh^\omega(
\vartheta_j,\varphi_k)\otimes d_jc_k. \qed
\end{split}
\end{equation*}
\renewcommand{\qed}{}
\end{pf}

Let $\sigma\colon \Gamma_{\inv}^{\otimes 2} \rightarrow
\Gamma_{\inv}^{\otimes 2}$ be the (left-invariant part
of the) canonical braid  operator
\cite{W3}. This map is explicitly \cite{D} given by
\begin{equation}\label{flip}
 \sigma(\eta\otimes\vartheta )=\sum_k \vartheta_k \otimes(\eta\circ
c_k)
\end{equation}
where $\adj (\vartheta )=\Sum_k \vartheta_k \otimes c_k$. Let
$m_{\Omega}$ be the multiplication map in $\Omega(P)$.

If $\omega\in\rc(P)$ then
\begin{equation}\label{regflip}
m_{\Omega}\Bigl\{\omega\otimes\varphi
\Bigr\}=(-1)^{\partial\varphi}m_{\Omega}\Bigl\{\varphi\otimes\omega
\Bigr\}\sigma,
\end{equation}
for each $\varphi\in\tau (P)$.
The above  equality  directly  follows  from \eqref{flip},
the tensoriality  of  $\varphi$,   and   the   definition   of
regular connections.

 Let  us  fix  a  linear  map  $\delta\colon \Gamma_{\inv}
\rightarrow\Gamma_{\inv}{\otimes}\Gamma_{\inv}$  such that
$$
\adj^{\otimes 2}\delta=(\delta\otimes
\id) \adj
$$
and such that if
$$\delta(\vartheta
)=\sum_k\vartheta^1_k \otimes\vartheta_k^2$$
then
$$
d(\vartheta)=\sum_k\vartheta_k^1 \vartheta_k^2\qquad
-\delta(\vartheta^*)=\sum_k\vartheta_k^{2*}\otimes
\vartheta^{1*}_k
$$
for   each   $\vartheta\in\Gamma_{\inv}$.
Such maps  will  be
called {\it embedded differentials}. For each
$\vartheta\in\Gamma_{\inv}$ there exists $a\in \ker(\e)$ such that
\begin{equation}\label{delpi}
\begin{split}
\delta(\vartheta )&=-\pi(a^{(1)})\otimes\pi (a^{(2)})\\
\pi(a)&=\vartheta.
\end{split}
\end{equation}

 For  given  linear  maps  $\varphi,\eta\colon \Gamma_{\inv}
\rightarrow\Omega(P)$  let us define (as in \cite{D}) new linear maps
$ \langle  \varphi,\eta\rangle ,[\varphi,\eta]\colon
\Gamma_{\inv}\rightarrow\Omega(P)$ by
\begin{align}
\langle \varphi,\eta\rangle &=m_{\Omega}(\varphi\otimes\eta
)\delta\label{bra1}\\
[\varphi,\eta]&=m_{\Omega}(\varphi\otimes\eta)c^\top,\label{bra2}
\end{align}
where $c^\top\colon
\Gamma_{\inv}\rightarrow\Gamma_{\inv}{\otimes}\Gamma_{\inv}$ is
the ``transposed commutator'' map \cite{W3}, explicitly given by
\begin{equation}\label{tcom}
c^\top=(\id \otimes\pi)\adj .
\end{equation}

 The same brackets can be used for linear maps defined on $\Gamma_{\inv}$,
with values in an arbitrary algebra.

 It is easy to see that  if $\varphi\in\psi^i (P)$  and
$\eta\in\psi^j (P)$ then
$\langle \varphi,\eta\rangle ,
[\varphi,\eta]\in\psi^{i+j}(P)$ (the same holds for
$\tau(P)$). Further,
\begin{equation}
\langle \varphi,\eta\rangle ^* =-(-1)^{ij}\langle \eta^*,
\varphi^*\rangle ,\label{bra1*}
\end{equation}
as directly follows from \eqref{*pse} and the hermicity of $\delta$.

     For an arbitrary $\omega\in \con(P)$ let us consider a map
$R_{\omega}\colon \Gamma_{\inv}\rightarrow\Omega(P)$ given by
\begin{equation}
R_{\omega}=d\omega-\langle \omega,\omega\rangle .\label{defR}
\end{equation}
Clearly, this is a pseudotensorial 2-form. Moreover,
\begin{lem}
We have
\begin{equation}
\widehat{F}R_{\omega}(\vartheta )=(R_{\omega}\otimes
\id) \adj (\vartheta )\label{dxFR}
\end{equation}
for each $\vartheta\in\Gamma_{\inv}$.
\end{lem}
\begin{pf}
A direct computation gives
\begin{equation*}
\begin{split}
\widehat{F}R_{\omega}(\vartheta )&=d\widehat{F}\omega
(\vartheta)+\widehat{F}\omega\pi(a^{(1)})\widehat{F}\omega\pi(a^{(
2)})
=d\bigl[\omega\pi(a^{(2)})\otimes \k(a^{(1)})a^{(3)}
+1\otimes\pi(a)\bigr]\\
&\phantom{=}+1\otimes\pi(a^{(1)})\pi(a^{(2)})+
\omega \pi(a^{(2)})\omega\pi(a^{(5)})\otimes
\k(a^{(1)})a^{(3)}\k(a^{(4)})a^{(6)}\\
&\phantom{=}+\omega\pi(a^{(2)})\otimes \k(a^{(1)})a^{(3)}\pi(a^{(4)})-
\omega\pi(a^{(3)})\otimes\pi (a^{(1)})\bigl(\k(a^{(2)})a^{(4)}\bigr)\\
&=\Bigl[d\omega-\langle \omega,\omega\rangle \Bigr]\pi(a^{(2)})\otimes
\k(a^{(1)})a^{(3)}+\omega\pi(a^{(2)})\otimes \k(a^{(1)})da^{(3)}\\
&\phantom{=}-\omega\pi(a^{(2)})\otimes d\bigl(\k(a^{(1)})a^{(3)}\bigr)
-\omega\pi(a^{(4)}
)\otimes \k(a^{(1)})da^{(2)}\bigl(\k(a^{(3)})a^{(5)}\bigr)\\
&=(R_{\omega}\otimes \id) \adj (\vartheta ).
\end{split}
\end{equation*}
Here, it is assumed that $a\in \ker(\e)$ satisfies \eqref{delpi} and
we have applied \eqref{pi}, \eqref{dpi}, \eqref{adpi} and
\eqref{dxFcon}.
\end{pf}
\begin{defn}
The constructed map $R_{\omega}$  is called
{\it the curvature} of $\omega$.
\end{defn}

 The curvature $R_{\omega}$ implicitly depends on  the  choice
of $\delta$. This dependence disappears if $\omega$ is
multiplicative.

 We are going to introduce the operator of  covariant  derivative.
This operator  will  be  first  defined  on  a  restricted  domain
consisting of horizontal forms. Later on,  after  constructing  the
horizontal projection operator, we shall extend the domain of  the
covariant derivative to the whole algebra $\Omega(P)$.

     For each $\omega\in \con(P)$ and $\varphi\in \hor(P)$ let us
define a new form
\begin{equation}
D_{\omega}(\varphi )=d\varphi -(-1)^{\partial\varphi}\sum_k
\varphi_k\omega\pi(c_k),\label{Dhor}
\end{equation}
where $F^{\wedge}(\varphi)=\Sum_k \varphi_k \otimes c_k$.
\begin{lem}
We have
\begin{equation}
\widehat{F}D_{\omega}(\varphi)=\sum_kD_{\omega}(\varphi_k)\otimes
c_k\label{dxFD}
\end{equation}
for each $\varphi\in \hor(P)$. In particular
$$D_{\omega}\hor(P)\subseteq \hor(P).$$
\end{lem}
\begin{pf}We compute
\begin{equation*}
\begin{split}
\widehat{F}D_{\omega}(\varphi)&=\widehat{F}d\varphi-(-1)^
{\partial\varphi}\sum_k\widehat{F}(\varphi_k )\widehat{F}\omega\pi
(c_k)\\
&=d\biggl(\sum_k\varphi_k \otimes c_k\biggr)-(-1)^{\partial\varphi}\sum_k
\varphi_k\otimes c_k^{(1)}\pi(c_k^{(2)}) \\
&\phantom{=}-(-1)^{\partial\varphi}\sum_k\varphi_k
\omega\pi(c_k^{(3)})\otimes c^{(1)}_k \k(c_k^{(2)})c_k^{(4)}\\
&=\sum_k d\varphi_k \otimes c_k
-(-1)^{\partial\varphi}\sum_k\varphi_k \omega\pi
(c_k^{(1)})\otimes c_k^{(2)}\\
&=\sum_k D_{\omega}(\varphi_k)\otimes c_k.
\end{split}
\end{equation*}
Hence, $\hor(P)$ is $D_\omega$-invariant.
\end{pf}
\begin{defn}
The constructed map $D_{\omega}
\colon \hor(P)\rightarrow \hor(P)$  is  called
{\it the covariant derivative} associated to $\omega$.
\end{defn}
\begin{pro}\label{pro:propD}
\bla{i} The diagram
\begin{equation}\label{wFD}
\begin{CD}
\hor(P) @>{\mbox{$F^\wedge$}}>> \hor(P)\otimes\cal{A}\\
@V{\mbox{$D_\omega$}}VV  @VV{\mbox{$D_\omega\otimes \id $}}V\\
\hor(P) @>>{\mbox{$F^\wedge$}}> \hor(P)\otimes\cal{A}
\end{CD}
\end{equation}
is commutative.

\smallskip
\bla{ii} We have
\begin{equation}
D_{\omega}^2 (\varphi)=-\sum_k \varphi_k
\left[d\omega\pi(c_k)+\omega\pi (c_k^{(1)})\omega\pi
(c_k^{(2)})\right]\label{DD1}
\end{equation}
for each $\varphi\in \hor(P)$. In particular, if $\omega$ is
multiplicative then
\begin{equation}
D_\omega^2(\varphi)=-\sum_k \varphi_k R_{\omega}\pi(c_k).
\label{DD2}
\end{equation}

\bla{iii} If $\omega$ is regular then
\begin{align}
D_{\omega}(\varphi\psi)&=D_{\omega}(\varphi)\psi +(-1)^{\partial
\varphi}\varphi D_{\omega}(\psi)\label{Dprod}\\
D_{\omega}(\varphi^*)&=D_{\omega}(\varphi)^*,\label{D*}
\end{align}
for each $\varphi,\psi\in \hor(P)$.

\smallskip
\bla{iv} We have
\begin{equation}\label{d-Dw}
(d-D_{\omega})\bigl(\Omega(M)\bigr)=\{0\}
\end{equation}
for each $\omega\in \con(P)$.
\end{pro}
\begin{pf}
Diagram  \eqref{wFD} follows  from  identity  \eqref{dxFD}.
Property ({\it iv}\/) follows from definitions of $\Omega(M)$ and
$D_{\omega}$. For each $\varphi\in \hor(P)$ we have
\begin{equation*}
\begin{split}
D_{\omega}^2(\varphi
)&=D_{\omega}\biggl(d\varphi-(-1)^{\partial\varphi}\sum_k\varphi_k
\omega\pi(c_k)\biggr)\\
&=-(-1)^{\partial\varphi}\sum_kd\varphi_k \omega\pi(c_k
)-\sum_k\varphi_kd\omega\pi(c_k )\\
&\phantom{=}-(-1)^{\partial\varphi}\sum_k\left(-d\varphi_k
\omega\pi(c_k )+(-1)
^{\partial\varphi}\varphi_k\omega\pi(c_k^{(3)})\omega\pi\bigl(c_k^{(1)}\k(
c_k^{(2)})c_k^{(4)}\bigr)\right)\\
&=-\sum_k\varphi_k\left(d\omega\pi(c_k)+\omega\pi(c_k^{(1)}
)\omega\pi(c_k^{(2)})\right).
\end{split}
\end{equation*}
If $\omega$ is multiplicative then \eqref{defR}
and  the  definition  of
$\delta$ imply that \eqref{DD2} holds.

 Finally, let us assume that $\omega$  is  regular.  Then
\begin{equation*}
\begin{split}
D_{\omega}(\varphi\psi)&=(d\varphi)\psi+(-1)^{\partial\varphi}
\varphi d\psi\\
&\phantom{=}-(-1)^{\partial\varphi+\partial\psi}\sum_{kl}
\Bigl(\varphi_k
\psi_l \omega \bigl(\pi(c_k )\circ d_l\bigr)+
\varphi_k \psi_l \e(c_k )\omega\pi(d_l)\Bigr)\\
&=\biggl(d\varphi-(-1)^{\partial\varphi}\sum_k \varphi_k \omega\pi
(c_k)\biggr)\psi+
(-1)^{\partial\varphi}\varphi\biggl(d\psi-(-1)^{\partial\psi}\sum_l
\psi_l \omega\pi(d_l)\biggr)\\
&=D_{\omega}(\varphi)\psi+(-1)^{\partial\varphi}\varphi D_{\omega}
(\psi),
\end{split}
\end{equation*}
where
$\Sum_l \psi_l\otimes
d_l=F^{\wedge} (\psi)$. Further,
\begin{equation*}
\begin{split}
    D_{\omega}(\varphi)^*& =d(\varphi^*)+\sum_k
\omega\pi\bigl(\k(c_k)^*\bigr)\varphi_k^* \\
&= d(\varphi^*)+(-1)^{\partial\varphi}\sum_k\varphi_k^*
\omega\bigl[\pi\bigl(\k(c_k^{(2)})^*\bigr)\circ c_k^{(1)*}\bigr]\\
&=d(\varphi^*)-(-1)^{\partial\varphi}\sum_k\varphi_k^* \omega\pi(c_k^*)=
D_{\omega}(\varphi^*).
\end{split}
\end{equation*}
Hence, properties ({\it iii}\/) hold.
\end{pf}

Actually, a connection $\omega$ is regular if and only if $D_\omega$
satisfies the graded Leibniz rule. Because of \eqref{wFD}, the
space  $\tau (P)$  is  closed  under taking
compositions with $D_{\omega}.$
This fact enables us to define the covariant
derivative (which will be  denoted  by  the  same  symbol)  as  an
operator acting in the space $\tau (P)$.
\begin{pro}\label{pro:Dten}
We have
\begin{equation}
D_{\omega}\varphi=d\varphi-(-1)^{\partial\varphi}[\varphi,\omega]
\label{Dten}
\end{equation}
for each $\omega\in \con(P)$ and $\varphi\in\tau (P)$.
\end{pro}
\begin{pf}
This directly  follows  from  the tensoriality of $\varphi$ and
from definitions  of
$D_{\omega}$ and brackets $[,]$.
\end{pf}

 For  a  given $\omega\in \con(P)$  let
$q_{\omega}\colon \psi (P)\rightarrow\psi (P)$ be a linear map
defined by
\begin{equation}
q_{\omega} (\varphi )=\langle \omega,\varphi \rangle -(-1)^{\partial\varphi}
\langle \varphi ,\omega \rangle -(-1)^{\partial\varphi}[\varphi ,\omega ].
\label{qdef}
\end{equation}
\begin{lem}\label{lem:dxFq}
\bla{i} We have
\begin{equation}
\widehat{F}q_{\omega}(\varphi)(\vartheta)=(q_{\omega}\otimes
\id) F^{\wedge} \varphi(\vartheta)\label{dxFq}
\end{equation}
for each $\varphi\in\tau (P)$ and $\vartheta\in\Gamma_{\inv}$.
In particular $q_{\omega}\tau (P)\subseteq \tau (P)$.

\smallskip
\bla{ii} If $\omega\in\rc(P)$ then
\begin{equation}\label{qreg}
\bigl(q_{\omega}\restr\tau (P)\bigr)=0.
\end{equation}
\end{lem}
\begin{pf}
We compute
\begin{equation*}
\begin{split}
\widehat{F}q_{\omega}(\varphi
)(\vartheta)&=-\widehat{F}\bigl[\omega\pi(a^{(1)})\varphi\pi(a^{(2)})\bigr]
+(-1)^{\partial\varphi}\widehat{F}\bigl[\varphi\pi(a^{(1)})
\omega\pi(a^{(2)})\bigr]\\
&\phantom{=}-
(-1)^{\partial\varphi}\sum_k\widehat{F}\bigl[\varphi(\vartheta_k
)\omega\pi(c_k)\bigr]\\
&=-\omega\pi(a^{(2)})\varphi\pi(a^{(3)})\otimes
\k(a^{(1)})a^{(4)}\\
&\phantom{=}-(-1)^{\partial\varphi}\varphi\pi(a^{(3)})
\otimes\pi(a^{(1)}   )\k(a^{(2)}   )a^{(4)}\\
&\phantom{=}+(-1)^{\partial\varphi}\varphi\pi(a^{(2)})
\omega\pi(a^{3})\otimes \k(a^{(1)})a^{(4)}\\
&\phantom{=}+(-1)^{\partial\varphi}\varphi\pi(a^{(2)})\otimes \k(a^{(1)})
a^{(3)}\pi(a^{(4)})\\
&\phantom{=}-(-1)^{\partial\varphi}
\sum_k\left(\varphi(\vartheta_k )\omega\pi(c_k^{(3)})\otimes
c_k^{(1)}\k(c_k^{(2)})c_k^{(4)}+\varphi(\vartheta_k )\otimes
c_k^{(1)}\pi(c_k^{(2)})\right)\\
&=\sum_k\Bigl\{\langle \omega ,\varphi\rangle (\vartheta_k )\otimes c_k
-(-1)^{\partial\varphi}\langle \varphi,\omega\rangle (\vartheta_k)\otimes
c_k\Bigr\}\\
&\phantom{=}-(-1)^{\partial\varphi}\sum_k[\varphi,\omega](\vartheta_k)
\otimes c_k -(-1)^{\partial\varphi}\sum_k\varphi(\vartheta_k )\otimes dc_k\\
&\phantom{=}-(-1)^{\partial\varphi}\varphi\pi(a^{(4)})\otimes
\k(a^{(1)})(da^{(2)})\k(a^{(3)})a^{(5)}\\
&\phantom{=}+(-1)^{\partial\varphi}\varphi\pi(a^{(2)})\otimes
\k(a^{(1)})da^{(3)}=\sum_kq_{\omega}(\varphi)(\vartheta_
k)\otimes c_k.
\end{split}
\end{equation*}
Here, $a\in \ker(\e)$ satisfies \eqref{delpi}, and $\Sum_k
\vartheta_k \otimes c_k =\adj (\vartheta)$.

 Let us assume that $\omega$ is regular. Applying \eqref{mod}, \eqref{adpi}
and \eqref{defreg1}, and the definition of brackets $\langle ,\rangle $
and $[,]$ we obtain
\begin{equation*}
\begin{split}
\langle \omega,\varphi\rangle (\vartheta)&=-\omega\pi(a^{(1)})\varphi
\pi(a^{(2)})\\
&=-(-1)^{\partial\varphi}\varphi\pi(a^{(3)})\omega\bigl[\pi(a^{(1)})
\circ \bigl(\k(a^{(2)})a^{(4)}\bigr)\bigr]\\
&=-(-1)^{\partial\varphi}
\varphi\pi(a^{(3)})\omega\pi\bigl[\bigl(a^{(1)}-\e(a^{(1)})1\bigr)
\k(a^{(2)})a^{(4)}\bigr]\\
&=-(-1)^{\partial\varphi}\varphi\pi(a^{(1)})\omega\pi(a^{(2)})+
(-1)^{\partial\varphi}\varphi\pi(a^{(2)})
\omega\pi\bigl(\k(a^{(1)})a^{(3)}\bigr)\\
&=(-1)^{\partial\varphi}
\langle \varphi,\omega\rangle (\vartheta)+
(-1)^{\partial\varphi}[\varphi,\omega](\vartheta).
\end{split}
\end{equation*}
Hence, \eqref{qreg} holds.
\end{pf}

Actually, the following equality holds
$$ q_\omega(\varphi)(\vartheta)=
\sum_k\lh^\omega\bigl(\vartheta_k^1,\varphi(\vartheta_k^2)\bigr)$$
for each $\vartheta\in\Gamma_{\inv}$ and $\varphi\in\tau(P)$.
The next proposition gives a
quantum  counterpart  for  the classical Bianchi identity.
\begin{pro}
We have
\begin{equation}
(D_{\omega} -q_{\omega})(R_{\omega})
=\langle \omega,\langle \omega,\omega\rangle
\rangle -\langle \langle \omega,\omega\rangle ,\omega\rangle \label{bianchi}
\end{equation}
for each $\omega\in \con(P)$.
\end{pro}
\begin{pf}
A direct computation gives
\begin{equation*}
\begin{split}
(D_{\omega}
-q_{\omega})(R_{\omega})&=dR_{\omega}-\langle \omega,R_{\omega}
\rangle +\langle R_{\omega},\omega\rangle
=-d\langle \omega,\omega\rangle \\
&\phantom{=}-\langle \omega,
d\omega -\langle \omega,\omega\rangle \rangle +\langle d\omega-
\langle \omega,\omega\rangle ,\omega\rangle \\
&=\langle \omega,\langle \omega,\omega\rangle \rangle -\langle
\langle \omega,\omega\rangle ,\omega\rangle .\qed
\end{split}
\end{equation*}
\renewcommand{\qed}{}
\end{pf}

 If the connection $\omega$ is multiplicative  then  the  right-hand
side of equality \eqref{bianchi} vanishes. On the other hand, if
$\omega$ is regular then the second summand on the left-hand side
vanishes.

     It is important to point out that regular connections  are  not
necessarily  multiplicative.  However,  there  exists   a   common
obstruction to multiplicativity for all  regular  connections,  so
that if one regular connection is multiplicative, then every regular
connection possesses the same property.
This obstruction  will  be  now  analyzed  in  more
details.

 In general, the lack of multiplicativity of a  connection
$\omega$
is measured by the map $r_{\omega}\colon \cal{R}\rightarrow\Omega(P)$
given  by $r_{\omega}(a)=\omega\pi(a^{(1)})\omega\pi(a^{(2)}).$
\begin{lem}
\bla{i} The following identities hold
\begin{gather}
 r_{\omega}\bigl(\k(a)^*\bigr)=-r_{\omega}(a)\label{r*k}\\
\pi_v r_{\omega}(a)=0\label{pivr}\\
\widehat{F}r_{\omega}(a)=F^{\wedge}r_{\omega}(a)=(r_{\omega}
\otimes \id) \ad(a).\label{dxFr}
\end{gather}
In particular $r_{\omega}(a)\in \hor^2 (P)$ for each $a\in\cal{R}$.

\smallskip
\bla{ii} Let us assume that $P$ admits regular connections.
The map $\omega\mapsto r_{\omega}$ is constant on equivalence
classes from the
space $\con(P)/\rc(P)$. If $\omega\in\rc(P)$ then
\begin{equation}
r_{\omega}(a)\varphi=\sum_k\varphi_k r_{\omega}(ac_k)
\label{regr}
\end{equation}
for   each   $a\in\cal{R}$    and $\varphi\in \hor(P)$. Furthermore,
\begin{equation}
dr_{\omega}(a)=\langle \omega,\omega\rangle \pi(a^{(1)})\omega
\pi(a^{(2)})-\omega
\pi(a^{(1)})\langle \omega,\omega\rangle \pi(a^{(2)}).\label{dr}
\end{equation}
\end{lem}
\begin{pf} A direct computation gives
\begin{equation*}
\begin{split}
r_{\omega}(a)^*&=-\omega\pi(a^{(2)})^*\omega\pi(a^{(1)})^*=-\omega
\pi\bigl(\k(a^{(2)})^*\bigr)\omega\pi\bigl(\k(a^{(1)})^*\bigr)
=-r_{\omega}\bigl(\k(a)^*\bigr)\\
\intertext{and similarly}
\widehat{F}r_{\omega} (a)&=\left[(\omega\otimes
\id) \adj \pi(a^{(1)})+1\otimes\pi(a^{(1)}
)\right]\left[(\omega\otimes \id) \adj \pi(a^{(2)})
+1\otimes\pi(a^{(2)})\right]\\
&=\omega\pi(a^{(2)})\omega\pi(a^{(3)})\otimes \k(a^{(1)})a^{(4)}+\omega
\pi(a^{(2)})\otimes \k(a^{(1)} )a^{(3)}\pi(a^{(4)})\\
&\phantom{=}-\omega\pi(a^{(3)})\otimes\pi(a^{(1)})\k(a^{(2)})a^{(4)}\\
&=r_{\omega}(a^{(2)})\otimes \k(a^{(1)})a^{(3)}+
\omega\pi(a^{(2)})\otimes d\bigl[\k(a^{(1)})a^{(3)}\bigr]
=(r_{\omega}\otimes \id) \ad(a).
\end{split}
\end{equation*}

 Let  us  assume  that  $\omega$   is   regular.   Then
\begin{equation*}
\begin{split}
r_{\omega+\varphi}(a)&=\omega\pi(a^{(1)})\omega\pi(a^{(2)})-\varphi
\pi(a^{(3)})\omega\bigl[\pi(a^{(1)})\circ
\bigl(\k(a^{(2)})a^{(4)}\bigr)\bigr]\\
&\phantom{=}+\varphi\pi(a^{(1)})\varphi\pi(a^{(2)})+\varphi\pi(a^{(1)})
\omega\pi(a^{(2)})\\
&=\omega\pi(a^{(1)})\omega\pi(a^{(2)})+
\varphi\pi(a^{(1)})\varphi\pi(a^{(2)})+\varphi\pi(a^{(2)}
)\omega\pi\bigl(\k(a^{(1)})a^{(3)}\bigr)\\
&=\omega\pi(a^{(1)})\omega\pi(a^{(2)})+\varphi\pi(a^{(1)})
\varphi\pi(a^{(2)})
\end{split}
\end{equation*}
for   each $\varphi=\varphi^*\in\tau^1(P)$.

    Now if $\zeta=\zeta^*\in \tau^1 (P)$ satisfies \eqref{vreg1} then
$\varphi\pi(a^{(1)})\zeta\pi(a^{(2)})+\zeta\pi(a^{(1)})\varphi\pi(a^{(2)}
)=0$,  and in particular $\zeta\pi(a^{(1)})\zeta\pi(a^{(2)})=0$  for
each $a\in\cal{R}$. Hence,
$$
r_{\omega+\varphi+\zeta}=r_{\omega+\varphi}.$$

    Further, applying \eqref{mod}, \eqref{defR} and
the tensoriality of the curvature  we obtain
\begin{equation*}
\begin{split}
dr_{\omega}(a)&=R_{\omega}\pi(a^{(1)})\omega\pi(a^{(2)}
)+\langle \omega,\omega\rangle \pi(a^{(1)})\omega\pi(a^{(2)})\\
&\phantom{=}-\omega\pi(a^{(1)})R_{\omega}\pi(a^{(2)})-\omega\pi
(a^{(1)})\langle \omega,\omega\rangle \pi(a^{(2)})\\
&=R_{\omega}\pi(a^{(1)}   )\omega\pi(a^{(2)})+
\langle \omega,\omega\rangle \pi(a^{(1)})\omega\pi(a^{(2)})\\
&\phantom{=}-R_{\omega}\pi(a^{(3)})\omega\bigl[\pi(a^{(1)})\circ\bigl(
\k(a^{(2)})a^{(4)}\bigr)\bigr]
-\omega\pi(a^{(1)})\langle \omega,\omega\rangle \pi(a^{(2)})\\
&=\langle \omega,\omega\rangle \pi(a^{(1)})
\omega\pi(a^{(2)})-\omega\pi(a^{(1)}
)\langle \omega,\omega\rangle \pi(a^{(2)})\\
&\phantom{=}+R_{\omega}\pi(a^{(2)})\omega\pi\bigl(\k(a^{(1)})a^{(3)}\bigr)\\
&=\langle \omega,\omega\rangle \pi(a^{(1)})\omega\pi(a^{(2)})
-\omega\pi(a^{(1)}
)\langle \omega,\omega\rangle \pi(a^{(2)}).
\end{split}
\end{equation*}
Finally,
\begin{equation*}
\begin{split}
    r_{\omega}(a)\varphi&=\sum\varphi_k \omega\bigl[\pi(a^{(1)})\circ
c_k^{(1)}\bigr]\omega\bigl[\pi(a^{(2)})\circ c_k^{(2)}\bigr]\\
&=\sum_k\varphi_k\omega\pi\bigl[\bigl(a^{(1)}-\e(a^{(1)})1\bigr)
c_k^{(1)}\bigr]\omega
\pi\bigl[\bigl(a^{(2)}-\e(a^{(2)})1\bigr)c_k^{(2)}\bigr]\\
&=\sum_k\varphi_k \omega\pi(a^{(1)}c_k^{(1)})\omega
\pi(a^{(2)}c_k^{(2)})-\sum_k
\varphi_k \omega\pi(ac_k^{(1)})\omega\pi(c_k^{(2)})\\
&\phantom{=}-\sum_k\varphi_k \omega\pi(c_k^{(1)})
\omega\pi(ac_k^{(2)})+\e(a)\sum_k \varphi_k\omega\pi(c_k^{(1)})
\omega\pi(c_k^{(2)})\\
&=\sum_k\varphi_k\omega\pi(a^{(1)}c_k^{(1)})\omega
\pi(a^{(2)}c_k^{(2)})=\sum_k\varphi_k r_{\omega}(ac_k)
\end{split}
\end{equation*}
for each $\varphi\in \hor(P)$.
Here, we have used \eqref{defreg1}, and the fact that
$\cal{R}\subseteq \ker(\e)$  is  a right ideal.
\end{pf}

 Let us assume that $P$ admits regular connections, and let
$\J(P)$  be
the left ideal in $\Omega(P)$ generated by  the  space
$r_{\omega}(\cal{R})$,
for some $\omega\in \rc(P)$.
\begin{lem}\bla{i} The following properties hold:
\begin{gather}
\J(P)^* =\J(P)\label{J*}\\
\widehat{F}\J(P)\subseteq \J(P)\otimes\Gamma^{\wedge}\label{dxFJ}\\
\pi_v \J(P)=\{0\}\label{pivJ}\\
d\J(P)\subseteq \J(P).\label{dJ}
\end{gather}

\bla{ii} The space $\J(P)$ is a two-sided ideal in $\Omega(P)$.
\end{lem}
\begin{pf}
Properties \eqref{dxFJ}--\eqref{pivJ} directly  follow
from  identities
\eqref{pivr}--\eqref{dxFr}. Concerning \eqref{dJ},
it  follows  from  \eqref{dr},  and  the
following  observations.   As   first, $\delta\pi(a)+\pi
(a^{(1)})\otimes\pi(a^{(2)})$ belongs to $S_{\inv}^{\wedge 2}$ for each
$a\in\cal{A}$. Secondly
$\omega^\otimes(S_{\inv}^{\wedge 2})=r_\omega(\cal R)$. Thirdly, because
of \eqref{defreg2} and \eqref{dxFr},
$$
 r_{\omega}(a)\omega(\vartheta)\in \J(P)$$
for each $a\in\cal{R}$ and $\vartheta\in\Gamma_{\inv}$.

 Let us prove ({\it ii}\/). Because of the above inclusion  and  the  fact
that $\omega(\Gamma_{\inv})$  and  $\hor(P)$  generate $\Omega(P)$
(as follows from the fact that $\mu_\omega$ splits the sequence
\eqref{seq}) it is sufficient to check that
$$
r_{\omega}(a)\varphi\in \J(P)
$$
for each $a\in\cal{R}$ and  $\varphi\in \hor(P)$.  However,
this follows
from \eqref{regr} in a straightforward way. Finally, \eqref{J*}
follows from ({\it ii}\/) and \eqref{r*k}.
\end{pf}

     The ideal $\J(P)$  measures  the  lack  of  multiplicativity  of
regular connections.
Let $\omega$ be a regular connection on $P$, and let  us  assume
that
$\J(P)=\{0\}$.  Then  $\omega$ is multiplicative and the  map
$\omega^{\wedge}\colon \Gamma_{\inv}^{\wedge}\rightarrow\Omega(P)$
possesses the  following commutation properties with horizontal
forms
\begin{align}
\omega^{\wedge}(\vartheta)\varphi&=(-1)^{\partial\varphi\partial
\vartheta}\sum_k\varphi_k\omega^{\wedge}(\vartheta\circ c_k)
\label{wreg1}\\
\varphi\omega^{\wedge}(\vartheta)&=(-1)^{\partial\varphi\partial
\vartheta}\sum_k\omega^{\wedge}
\bigl(\vartheta\circ \k^{-1}(c_k)\bigr)
\varphi_k,\label{wreg2}
\end{align}
as easily follows from \eqref{defreg1} and \eqref{defreg2}.

     The formulas
\begin{align}
 (\psi\otimes\eta)(\varphi\otimes\vartheta)&=
(-1)^{\partial\varphi\partial\eta}\sum_k
\psi\varphi_k \otimes(\eta\circ c_k)\vartheta\label{hGprod}\\
(\varphi\otimes\vartheta)^*&=\sum_k\varphi_k^* \otimes(\vartheta^*
\circ c_k^*), \label{hG*}
\end{align}
determine  the   structure   of a *-algebra   in   the
space $\vh(P)=\hor(P)\otimes\Gamma_{\inv}^{\wedge}$.
The proof is essentially the
same as the proof of the *-algebra properties for $\ver(P)$ (actually
$\ver(P)$ is a *-subalgebra of $\vh(P)$).
Here as usual, $F^{\wedge}(\varphi)=\Sum_k\varphi_k \otimes c_k$.
The elements of $\vh(P)$ are interpretable as ``vertically-horizontally''
decomposed differential forms on $P$. In the following considerations
the space $\vh(P)$ will be endowed  with  this (graded) *-algebra
structure.

 We shall now construct, starting  from  an  arbitrary  connection
$\omega$, an important isomorphism between the spaces $\Omega(P)$  and
$\vh(P)$,   extending   the   splitting $\mu_{\omega}$.   This
isomorphism will be the base for the construction and the analysis
of  horizontal  projection  operators.  Also,   some   interesting
questions related to the structure of horizontal forms
will be considered.

Let us fix a splitting of the form
$$\Gamma_{\inv}^{\otimes}=\Gamma_{\inv}^{\wedge}\oplus
S_{\inv}^{\wedge}$$
in which $\Gamma_{\inv}^{\wedge}$
is realized as a complement to the space $S_{\inv}^{\wedge}$,
with the help of a *-preserving section
$\iota\colon \Gamma_{\inv}^{\wedge}\rightarrow\Gamma_{\inv}^{\otimes}$
(of the factor-projection map) intertwining  the
adjoint  actions.  Further,  let  us  assume  that  the   embedded
differential $\delta$ is given by $$\delta(\vartheta)=
\iota d(\vartheta).$$
Finally, let
us consider  a
linear map $m_{\omega}\colon \vh(P)
\rightarrow\Omega(P)$ given by
\begin{align}
 m_{\omega}(\varphi\otimes\vartheta)&=\varphi\omega^{\wedge}
(\vartheta),\label{defm}\\
\intertext{where}
\omega^{\wedge}&=\omega^{\otimes}\iota.\label{defw}
\end{align}

 It is important to mention  that  the  above  definition  of  map
$\omega^{\wedge}$
reduces, for multiplicative connections, to the previously
introduced multiplicative extension. In particular, if $\omega$  is
multiplicative then $m_{\omega}$  is independent of the section
$\iota$. Further, if $\iota$ intertwines $\circ$-structures then
\eqref{wreg1}--\eqref{wreg2} hold for each $\omega\in\rc(P)$,
independently of the multiplicativity property.
By construction $\omega^\wedge$ is hermitian and
$$F^\wedge\omega^\wedge=(\omega^\wedge\otimes \id)
\adj^\wedge.$$
Let $F_{v\!h}\colon\vh(P)\rightarrow
\vh(P)\otimes\cal{A}$
be the product of actions $F^\wedge\colon \hor(P)\rightarrow\hor(P)\otimes
\cal A$ and
$\adj^{\wedge}$. This map  is  a  *-homomorphism.
\begin{thm}\label{pro:decW}
     \bla{i} The map $m_{\omega}$
is bijective.

\smallskip
\bla{ii} The diagram
\begin{equation}\label{mcov}
\begin{CD}
\vh(P) @>{\mbox{$m_\omega$}}>>
\Omega(P)\\
@V{\mbox{$F_{v\!h}$}}VV @VV{\mbox{$F^{\wedge}$}}V\\
\vh(P)\otimes\cal A
@>>{\mbox{$m_{\omega}\otimes \id $}}>
\Omega(P)\otimes\cal{A}
\end{CD}
\end{equation}
is commutative.

\smallskip
 \bla{iii} If $\omega$  is  regular  and  if  $\J(P)=\{0\}$  then
$m_{\omega}$ is an isomorphism of *-algebras.
\end{thm}
\begin{pf}
We shall first prove that $m_{\omega}$
is injective. Let us assume that
$w\in \ker(m_{\omega})$, and let
$w=\Sum_i\varphi_i \otimes\vartheta_i$
where
$\vartheta_i \in\Gamma_{\inv}^{\wedge}$ are homogeneous and  linearly
independent,
and $\varphi_i \in \hor(P)$ are non-zero elements. Let
$k$ be the maximum of degrees of elements $\vartheta_i$.
Then
$$
0=(\id \otimes p_k)\widehat{F}m_{\omega}(w)=\sideset{}{^\prime}
\sum_{i}F^{\wedge}
(\varphi_i)\vartheta_i,
$$
according to \eqref{dxFcon}, and definitions of $F^{\wedge}$  and
$\hor(P)$. The above summation   is   performed    over    the
indexes $i$ corresponding to $k$-th order elements.
This implies
${\Sum}_i^\prime\varphi_i\otimes\vartheta_i =0$,
which is a contradiction. Hence, $m_{\omega}$  is injective.

    We shall now prove that $m_{\omega}$  is surjective. For each
integer $k\geq 0$  let us consider the space
 $$
\Omega_k(P)=(\widehat{F})^{-1}\bigl(\Omega(P)\otimes\Gamma_k^{\wedge}
\bigr),
$$
where
$$\Gamma_k^{\wedge}=
\sideset{}{^\oplus}\sum_{i\leq k} \Gamma^{\wedge i}.$$
In other words, $w\in\Omega_k(P)$ iff $(\id\otimes p_l)\widehat{F}(w)=0$
for each $l>k$.
Clearly, $\Omega_0(P)=\hor(P)$ and
$\Omega_k(P)\subseteq\Omega_{k+1}(P)$ for each $k\geq 0$. Further
$$\displaystyle{\bigcup}_k\Omega_k(P)=\Omega(P).$$  We
are going to prove,  inductively,  that  the  spaces
$\Omega_k(P)$ are
contained in the  image  of  $m_{\omega}$.

Evidently,  the
statement is true for $k=0$. Let us assume that it
holds  for  some $k$.
Let us assume that there exists an element
$w\in\Omega_{k+1}(P)\big\backslash
\Omega_k(P)$. Then we have, according to Lemma~\ref{lem:tex}
$$
\widehat{F}(w)=\psi+\sum_iF^{\wedge}(\varphi_i)\vartheta_i,
$$
where  $\psi\in\Omega(P)\otimes\Gamma_k^{\wedge}$ and $\varphi_i
\in \hor(P)\setminus \{0\}$  while $\vartheta_i \in
(\Gamma^{\wedge}_{\inv})^{
1+k}$ are linearly independent elements.
On  the other hand,
$$
\widehat{F}m_{\omega}\biggl(\sum_i\varphi_i
\otimes\vartheta_i\biggr)=
\psi'+\sum_iF^{\wedge}(\varphi_i)\vartheta_i
$$
where $\psi'\in\Omega(P)\otimes\Gamma^{\wedge}_k$.
Hence,
$w-m_{\omega}\Bigl(\Sum_i\varphi_i\otimes\vartheta_i\Bigr)$
belongs to $\Omega_k(P)$ and therefore $w\in m_{\omega}\bigl(\vh(P)\bigr)$.
We conclude that $m_{\omega}$ is surjective. Thus, ({\it i}\/) holds.

 Intertwining property \eqref{mcov} directly follows from
the definition of $m_{\omega}$ and from the
right-covariance of $\omega^\wedge$.

 Finally, let us assume that $\J(P)=\{0\}$ and let
$\omega\in \rc(P)$. Then $\omega$ is multiplicative
and $\omega^\wedge$ is a *-homomorphism. Using
\eqref{wreg1} and \eqref{hGprod}--\eqref{hG*} we
obtain
\begin{equation*}
\begin{split}
m_{\omega}\bigl[(\psi\otimes \eta)(\varphi\otimes\vartheta)\bigr]&=
\sum_k (-1)^{\partial\varphi\partial\eta}\psi\varphi_k
\omega^{\wedge}\bigl[(\eta\circ c_k)\vartheta\bigr]\\
&=\sum_k (-1)^{\partial\varphi\partial\eta}\psi\varphi_k
\omega^{\wedge} (\eta\circ c_k)\omega^{\wedge}(\vartheta)\\
&=\psi\omega^{\wedge} (\eta)\varphi\omega^{\wedge}(\vartheta)=
m_{\omega}(\psi\otimes\eta)m_{\omega} (\varphi\otimes\vartheta),\\
\intertext{and similarly}
m_{\omega} \bigl[(\varphi\otimes\vartheta)^*\bigr]&=\sum_k\varphi_k^*
\omega^{\wedge}(\vartheta^*\circ c_k^*)=
(-1)^{\partial\varphi\partial\vartheta}\omega^{\wedge}(\vartheta^*
)\varphi^*\\
&=(-1)^{\partial\varphi\partial\vartheta}
\omega^{\wedge}(\vartheta)^* \varphi^* =
\bigl[\varphi\omega^{\wedge}(\vartheta)\bigr]^* =
\bigl[m_{\omega}(\varphi\otimes\vartheta)\bigr]^*.
\end{split}
\end{equation*}

     In other words, $m_{\omega}$  is a *-algebra isomorphism.
\end{pf}

The spaces $\Omega_k(P)$ introduced in the above proof
form a filtration of $\Omega(P)$, compatible
with the graded-differential *-algebra structure. In particular
$$\Omega_k(P)=\sideset{}{^\oplus}\sum_{j\geq 0}\Omega_k^j(P).$$
For each $k$ and $\omega\in\con(P)$ the space $\Omega_k(P)$
is linearly spanned
by elements of the form $\varphi\omega^\wedge(\vartheta)$, where
$\varphi\in\hor(P)$ and $\vartheta\in\Gamma_{\inv}^{\wedge i}$, with
$i\leq k$. Let $\mho(P)$ be the graded-differential *-algebra associated
to the filtered algebra $\Omega(P)$. The map
$\amalg\colon\vh(P)\rightarrow\mho(P)$ given by
$$\amalg(\varphi\otimes\vartheta)=\bigl[\varphi
\omega^\wedge(\vartheta)+\Omega_{k-1}(P)\bigr],$$
where $\vartheta\in\Gamma_{\inv}^{\wedge k}$, is bijective (and independent
of the choice of $\omega$ and $\iota$). Moreover, $\amalg$ is a
*-isomorphism, as follows from Lemma~\ref{lem:lac} and the horizontality
of $r_\omega$.
The introduced filtration
is compatible with the map $\widehat{F}$, in the sense that
$$\widehat{F}\bigl(\Omega_k(P)\bigr)\subseteq\Omega_k(P)
\otimes\Gamma^\wedge_k.$$
In other words, $\widehat{F}$ is factorizable
through the filtration.
The diagrams
\begin{equation}
\begin{CD}
\vh(P)@>{\mbox{$\amalg$}}>>\mho(P)\\
@V{\mbox{$d_{v\!h}$}}VV @VV{\mbox{$d$}}V\\
\vh(P)@>>{\mbox{$\amalg$}}>\mho(P)
\end{CD}\qquad\qquad
\begin{CD}
\vh(P)@>{\mbox{$\amalg$}}>>\mho(P)\\
@V{\mbox{$\widehat{F}_{v\!h}$}}VV @VV{\mbox{$\widehat{F}$}}V\\
\vh(P)\grten\Gamma^\wedge
@>>{\mbox{$\amalg\otimes\id$}}>\mho(P)
\grten\Gamma^\wedge
\end{CD}
\end{equation}
describe the corresponding factorized maps in terms of $\vh(P)$.
The differential
$d_{v\!h}\colon\vh(P)\rightarrow
\vh(P)$ is given by
$$ d_{v\!h}(\varphi\otimes\vartheta)=(-1)^{\partial\varphi}
\sum_k\varphi_k\otimes
\pi(c_k)\vartheta+(-1)^{\partial\varphi}\varphi\otimes d(\vartheta).$$
Similarly, $\widehat{F}_{v\!h}\colon\vh(P)
\rightarrow\vh(P)\grten\Gamma^\wedge$
is given by
$$ \widehat{F}_{v\!h}(\varphi\otimes\vartheta)=F^\wedge(\varphi)
\widehat{\adj}(\vartheta).$$
It is worth noticing that $d_{v\!h}\restr\ver(P)=d_v$ and
$\widehat{F}_{v\!h}\restr\ver(P)=\widehat{F}_v$.

The next lemma gives a more detailed  description  of
higher-order horizontal forms.
\begin{lem}
If the bundle
admits regular connections then the algebra
$$
\hor^+(P)= \sideset{}{^\oplus}\sum_{k \geq 1} \hor^k (P)
$$
is  generated   by   spaces $\hor^1(P)$   and
$r_{\omega}(\cal{R})$ $($where $\omega\in \rc(P)${}$)$.
In particular, if
$\J(P)=\{0\}$  then  higher-order horizontal  forms  are
algebraically expressible through the first-order ones.
\end{lem}
\begin{pf}
The algebra
$$\Omega^+(P)= \sideset{}{^\oplus}\sum_{k \geq 1}
\Omega^k(P)$$
is generated  by
the space $\Omega^1(P)$ of 1-forms.  This  fact,  together  with  the
*-bimodule splitting $\Omega^1(P)=\hor^1(P)\oplus \ver^1(P)$
determined by an arbitrary $\omega\in\rc(P)$,
can be used to prove that each element $w\in\Omega^n(P)$  is
expressible in the form
\begin{equation}
 w= \sum_{k=0}^n \biggl[\sum_i \varphi_{ik}\omega^{\wedge}
(\vartheta_{ik})+\sum_j \psi_{jk}\omega^{\otimes}(\eta_{jk})
\biggr],\label{dec1}
\end{equation}
where $\vartheta_{ik}$ and
$\eta_{jk}$ are linearly independent elements in the spaces
$\Gamma_{\inv}^{\wedge k}$ and $S_{\inv}^{\wedge k}$ respectively,
while $\varphi_{ik}$, $\psi_{jk}$ are
horizontal $(n-k)$-forms, expressible as sums of  products
of $n-k$ factors from $\hor^1(P)$. Now using the facts that
$S_{\inv}^{\wedge}$ is  generated by $S_{\inv}^{\wedge 2}$
and  that  $\omega^{\otimes}(Q)$  is  horizontal  for  each
$Q\in S_{\inv}^{\wedge 2}$ (because of $\omega^{\otimes}(Q)\in
r_{\omega}(\cal{R})$) we  can  prove,  inductively
applying \eqref{defreg1} and using the definition of $r_{\omega}$
and identity  \eqref{dxFr}, that the elements
$\omega^{\otimes}(\eta_{jk})$ are expressible as sums of
products of  the  form
$r_{\omega}(a_1)\cdots r_{\omega}(a_m)\omega^{\wedge}
(\vartheta_{ql})$  with $l+2m=k$, and a possibly extended set of
$\vartheta_{ki}$. Inserting this in \eqref{dec1} we
conclude that
\begin{equation}\label{dec2}
w=\sum_{k=0}^n\biggl[\sum_i\tilde{\varphi}_{ik}\omega^{\wedge}
(\vartheta_{ik})\biggr],
\end{equation}
where $\tilde{\varphi}_{ik}$ are horizontal $(n-k)$-forms expressible
as  sums  of
products of elements from $\hor^1(P)$ and $\im(r_{\omega})$. If
$w\in \hor^n(P)$  then
$\tilde{\varphi}_{ik}=0$ for each $1\leq k\leq n$,  according  to
Theorem~\ref{pro:decW} ({\it i}\/).
Hence, higher-order horizontal forms are algebraically expressible
via the first-order ones, and the 2-forms from the space
$r_{\omega}(\cal{R})$.
\end{pf}

In the general case ($P$ is an arbitrary bundle and
$\omega$ is an arbitrary connection) the algebra
$\hor^+(P)$ is generated by spaces $\hor^1(P)$, $r_\omega(\cal{R})$ and
horizontal forms obtained by iteratively acting by $\lh^\omega$ on the
elements from $\hor^1(P)$ and $r_\omega(\cal{R})$. It is interesting to
observe that in the general case
$$ r_\omega(a)\varphi=\sum_k\varphi_k r_\omega(a c_k)+
\lh^\omega\bigl(\pi(a^{(1)}),\lh^\omega(\pi(a^{(2)}),\varphi)\bigr). $$

     If the ideal $\J(P)$ is non-trivial (if every $\omega\in
\rc(P)$ is not multiplicative) then we can ``renormalize'' the
calculus, passing to
the  factoralgebra  $\Omega^\sstar (P)=\Omega(P)/\J(P)$
and  projecting
the   whole formalism from $\Omega(P)$ on $\Omega^\sstar(P)$.
It is worth noticing  that  such  a
factorization preserves the first-order calculus. In the framework
of   this   projected   calculus   regular   connections    become
multiplicative.  More  precisely,   let
$\frak{h}^\sstar(P)\subseteq\Omega^\sstar(P)$
be   the corresponding   algebra   of    horizontal    forms,
and let $\Pi\colon \Omega(P)\rightarrow\Omega^\sstar(P)$
be the projection map.
\begin{lem}
We have
\begin{equation}
\frak{h}^\sstar(P)=\Pi\bigl(\hor(P)\bigr).\label{projhor}
\end{equation}
In particular, if $\omega$ is a  regular  connection  relative  to
$\Omega(P)$ then $\Pi\omega$ is regular in terms of $\Omega^\sstar(P)$.
\end{lem}
\begin{pf} It is  evident  that  $\Pi(\hor(P))\subseteq \frak{h}^\sstar(P)$.
Let  $m^\star_{\omega}$ be the factorization map corresponding to
the calculus $\Omega^\sstar(P)$ and  to the connection $\Pi\omega$. We have
then
$$m^\star_{\omega}\bigl[\bigl(\Pi\restr \hor(P)\bigr)\otimes \id
\bigr]=\Pi m_{\omega}.$$
In particular, \eqref{projhor} holds.
\end{pf}

 With the help of the identification $m_{\omega}$ the corresponding {\it
horizontal
projection
operator\/} $h_{\omega}\colon \Omega(P)\rightarrow \hor(P)$
can be defined as follows
\begin{equation}
h_{\omega}=(\id \otimes p_0)m_{\omega}^{-1}.\label{defhproj}
\end{equation}
Clearly,  $h_{\omega}$ projects $\Omega(P)$ onto $\hor(P)$.

 The domain of the previously introduced covariant derivative
$D_{\omega}$
will be now extended from $\hor(P)$ to the whole algebra $\Omega(P)$.
Let  a
map $D_{\omega}\colon \Omega(P)\rightarrow \hor(P)$
be defined as follows
\begin{equation}
D_{\omega}=h_{\omega}d.\label{extD}
\end{equation}
This is a straightforward generalization of the corresponding
classical definition.
     The main properties  of $h_{\omega}$   and  $D_{\omega}$   are
collected in the following theorem.
\begin{thm}
\bla{i} The diagrams
\begin{equation}
\begin{CD}
\Omega(P) @>{\mbox{$F^\wedge$}}>> \Omega(P)\otimes\cal{A}\\
@V{\mbox{$h_\omega$}}VV  @VV{\mbox{$h_\omega\otimes \id $}}V\\
\hor(P) @>>{\mbox{$F^\wedge$}}> \hor(P)\otimes\cal{A}
\end{CD}\qquad\quad
\begin{CD}
\Omega(P) @>{\mbox{$F^\wedge$}}>> \Omega(P)\otimes\cal{A}\\
@V{\mbox{$D_\omega$}}VV  @VV{\mbox{$D_\omega\otimes \id $}}V\\
\hor(P) @>>{\mbox{$F^\wedge$}}> \hor(P)\otimes\cal{A}
\end{CD}\label{Dcov}
\end{equation}
are commutative.

\smallskip
     \bla{ii} The map $D_{\omega}$ extends  the previously
defined covariant derivative.

\smallskip
\bla{iii} If  $\omega$  is  regular  and  if  $\J(P)=\{0\}$
then $h_{\omega}$ is a *-homomorphism and
\begin{align}
D_{\omega}(wu)&=D_{\omega}(w)h_{\omega}(u)+(-1)^{\partial
w}h_{\omega}(w)D_{\omega} (u)\label{Dm}\\
D_{\omega}(w^* )&=D_{\omega}(w)^*,\label{*D}
\end{align}
for each $w,u\in\Omega(P)$.
\end{thm}
\begin{pf}
 The statement ({\it i}\/) follows from the construction of
$h_{\omega}$ and $D_{\omega}$,
and properties \eqref{dwF} and
\eqref{mcov}.
The statement ({\it iii}\/)
is a consequence of property ({\it iii\/}) in Theorem~\ref{pro:decW},
and elementary properties of $d\colon\Omega(P)\rightarrow\Omega(P)$.
Finally,  in  the  first  definition \eqref{Dhor}
of the covariant derivative, the differential of a horizontal
form  $\varphi$  is written as
$$
d\varphi=m_{\omega}\biggl[D_{\omega}(\varphi)\otimes1+(-1)
^{\partial\varphi}\sum_k\varphi_k \otimes \pi(c_k)\biggr],
$$
which implies that the new definition includes the previous
one.
\end{pf}

 According to \eqref{Dcov}, compositions of
 pseudotensorial  forms  with
$D_{\omega}$ and $h_{\omega}$ are
tensorial. In particular, it is possible to define, via  these
compositions,    the    covariant    derivative and the horizontal
projection as maps
$D_{\omega},h_{\omega}\colon \psi (P)\rightarrow\tau (P)$.

     The following proposition gives a more geometrical  description  of
the curvature map, establishing a close analogy with classical geometry.
\begin{pro}

We have
\begin{equation}
R_{\omega} =D_{\omega}\omega\label{RD}
\end{equation}
for each $\omega\in \con(P)$.
\end{pro}
\begin{pf} From definitions of $m_{\omega}$  and $R_{\omega}$
we find
$$R_{\omega}(\vartheta)\otimes 1=m_{\omega}^{-1}d\omega(\vartheta)-
1\otimes\delta(\vartheta), $$
for each $\vartheta\in\Gamma_{\inv}$. Hence,
$D_{\omega}\omega=h_{\omega}d\omega=R_{\omega}.$
\end{pf}

Let us assume that $\omega\in\rc(P)$ and $\J(P)=\{0\}$. The  above
proposition implies
\begin{equation}\label{R-reg}
\begin{aligned}
R_\omega(\vartheta)\varphi&=\sum_k\varphi_kR_\omega(\vartheta\circ
c_k)\\
\varphi R_\omega(\vartheta)&=\sum_k R_\omega\bigl(\vartheta\circ
\k^{-1}(c_k)\bigr)
\varphi_k
\end{aligned}
\end{equation}
for each $\varphi\in\hor(P)$ and $\vartheta\in\Gamma_{\inv}$.
Evidently,  the   above   commutation   relations   are   mutually
equivalent. To obtain the first it is  sufficient  to  act  by
$D_\omega$ on \eqref{defreg1}, and apply \eqref{Dm} and \eqref{RD}.
\section{Characteristic Classes}

     In this section a quantum generalization of classical
Weil's theory
of characteristic classes  will  be  presented.  Conceptually,  we
follow the exposition  of \cite{KN}.  As  in  the  classical
case, the main result will  be  a  construction  of  an  invariant
homomorphism defined on an algebra playing the role of  ``invariant
polynomials'' over the ``Lie algebra''  of  $G$,  with  values  in  the
algebra of cohomology classes of $M$.
We shall assume that the bundle $P$ admits regular connections,
and that $\J(P)=\{0\}$.

 For each $k\geq 0$ let $\cal{I}^k \subseteq\Gamma_{\inv}^{\otimes k}$
be the subspace of $\adj^{\otimes k}$-invariant
elements, and let $\cal{I}$ be  the
direct  sum  of  these spaces. Clearly,
$\cal{I}$ is a unital *-subalgebra of $\Gamma^{\otimes}_{\inv}$.

     Let us consider a  connection  $\omega$.
There  exists  the  unique
unital homomorphism $R_{\omega}^{\otimes}\colon \Gamma_{\inv}^{\otimes}
\rightarrow \Omega(P)$ extending the curvature $R_{\omega}$.
The map $R_{\omega}^{\otimes}$ is *-preserving, horizontally valued,
intertwines
$\adj^{\otimes}$
and $F^{\wedge}$, and multiplies degrees by $2$.
Here, we are interested for the  values  of  the  restriction
map $R_{\omega}^{\otimes}\restr\cal{I}$.
\begin{pro}If $\vartheta\in \cal{I}^k$ then the form
$R_{\omega}^{\otimes}(\vartheta)$
belongs to $\Omega^{2k}(M)$. Moreover,
if $\omega\in\rc(P)$ then $R_{\omega}^{\otimes}
(\vartheta)$ is closed.
\end{pro}
\begin{pf}
Equation \eqref{dxFR} implies
$$ \widehat{F}R_{\omega}^{\otimes}(\vartheta)=
\sum_k R_{\omega}^\otimes(\vartheta_k)\otimes c_k,$$
for each $\vartheta\in\Gamma_{\inv}^{\otimes}$, where
$\Sum_k\vartheta_k\otimes c_k=\adj^{\otimes}(\vartheta)$.
Now, the first statement follows from the
assumption $\vartheta\in \cal{I}^k$, and from the definition of
$\Omega(M)$. The second statement
follows from \eqref{Dprod}, ({\it iv}\/)--Proposition~\ref{pro:propD},
and from Bianchi
identity $D_{\omega}R_{\omega}=0$,
which holds for regular multiplicative connections.
\end{pf}

     Now we prove that the cohomological class of
$R_{\omega}^{\otimes}(\vartheta)$ in $\Omega(M)$ is
independent of $\omega\in\rc(P)$. Let $\tau$ be another regular connection,
and let
\begin{equation}
\omega_t =\omega +t\varphi
\end{equation}
where $\varphi=\tau-\omega$ and $t\in [0,1]$, be
the segment  in  the space $\rc(P)$ determined by $\omega$ and
$\tau$.
\begin{lem}
We have
\begin{equation}\label{dtR}
(d/dt)R_{\omega_t}=D_{\omega_t}(\varphi).
\end{equation}
\end{lem}
\begin{pf}
Applying Lemma~\ref{lem:dxFq} ({\it ii}\/), Proposition~\ref{pro:Dten}
and \eqref{defR} we obtain
\begin{equation*}
\begin{split}
(d/dt)R_{\omega_t}&=(d/dt)[d\omega+td\varphi-
\langle \omega+t\varphi,\omega+t\varphi\rangle ]\\
&=d\varphi-\langle \varphi,\omega\rangle -\langle
\omega,\varphi\rangle -2t\langle \varphi,\varphi\rangle \\
&=d\varphi-\langle \omega+t\varphi,\varphi\rangle
-\langle \varphi,\omega+t\varphi\rangle
=D_{\omega_t}(\varphi).\qed
\end{split}
\end{equation*}
\renewcommand{\qed}{}
\end{pf}

 Let  us  consider an element $\vartheta\in\Gamma_{\inv}^{\otimes
k}$ and let
$\vartheta=\Sum_i c_i\vartheta_1^i\otimes\cdots\otimes
\vartheta^i_k$, where $c_i$ are complex numbers and
$\vartheta^i_j\in\Gamma_{\inv}$. Applying \eqref{dtR}
and property \eqref{Dprod}, the definition of
$R_{\omega}^{\otimes}$
and the Bianchi identity we obtain
\begin{equation*}
\begin{split}
(d/dt)R^{\otimes}_{\omega_t}(\vartheta)&=\sum_i c_i\bigl[
D_{\omega_t}(\varphi)(\vartheta_1^i)\ldots R_{\omega_t}
(\vartheta_k^i)+\dots +R_{\omega_t}(\vartheta_1^i)\dots
D_{\omega_t}(\varphi)(\vartheta_k^i)\bigr]\\
&=\sum_i c_iD_{\omega_t}\bigl[\varphi(\vartheta_1^i)
\ldots R_{\omega_t}
(\vartheta_k^i)+\dots+R_{\omega_t}(\vartheta_1^i)\dots
\varphi(\vartheta_k^i)\bigr].
\end{split}
\end{equation*}

Using the tensoriality property of $\varphi$ and $R_{\omega_t}$
we see that  if
$\vartheta\in \cal{I}^k$  then the form
$$
\psi_t(\vartheta)=
\sum_i c_i\bigl[\varphi(\vartheta_1^i)
\ldots R_{\omega_t}
(\vartheta_k^i)+\cdots+R_{\omega_t}(\vartheta_1^i)\ldots
\varphi(\vartheta_k^i)\bigr]
$$
belongs to $\Omega^{2k}(M)$. Hence
\begin{equation}
(d/dt)R_{\omega_t}^{\otimes}(\vartheta)=d\psi_t (\vartheta),
\end{equation}
according to \eqref{d-Dw}.
Integrating the above equality from $0$ to $1$ we obtain
\begin{equation}
R_{\tau}(\vartheta)=R_{\omega}(\vartheta)+
d\Bigl(\int_0^1\!\psi_t(\vartheta)\,dt\Bigr).\label{Rint}
\end{equation}

 Let $H(M)$ be the graded *-algebra of cohomology classes associated
to $\Omega(M)$. We have proved the following theorem.
\begin{thm}
\bla{i} The cohomological class of
$R^{\otimes}_{\omega}(\vartheta)$ in
$\Omega(M)$ is independent of the choice  of  a  regular  connection
$\omega$, for each $\vartheta\in \cal{I}$.

\smallskip
\bla{ii} The map $W\colon \cal{I}\rightarrow H(M)$ given by
\begin{equation}
W(\vartheta)=[R_{\omega}^{\otimes}(\vartheta)]
\end{equation}
is a unital *-homomorphism. $\qed$
\end{thm}

     The homomorphism $W$ plays the role of the Weil  homomorphism
in classical differential geometry \cite{KN}.

     In fact, in classical  geometry  the  domain  of  the  Weil
homomorphism is restricted to the algebra of {\it symmetric}
invariant
elements  of   the   corresponding   tensor   algebra   (invariant
polynomials). However, besides simplifying the domain of $W$ such  a
restriction  gives  nothing  new:  the   image   of   the   Weil
homomorphism will be the same.

 A  similar  situation  holds  in  the  noncommutative  case.  Let
$\cal{S}$ be the graded *-algebra obtained from
$\Gamma^{\otimes}_{\inv}$    by  factorizing
through the ideal $\cal{J}$ generated  by  the space
$\im(I- \sigma)\subseteq\Gamma_{\inv}^{\otimes 2}$.

The algebra $\cal{S}$ plays the role of polynoms over  the  ``Lie
algebra'' of $G$. The adjoint action
$\adj^\otimes$ naturally induces the action
$\adj_{\cal{S}}\colon\cal{S}\rightarrow\cal{S}\otimes\cal{A}$.
Let $\cal{I}_{\cal{S}}\subseteq \cal{S}$
be  the  subalgebra  consisting of elements  invariant   under
$\adj_{\cal{S}}$. Clearly,
$$\cal{I}_{\cal{S}}=\cal{I}/(\cal{I}\cap\cal{J}).$$
\begin{lem}
If $\omega$ is regular then
\begin{equation}
R_{\omega}^{\otimes}\sigma(\vartheta)=R_{\omega}^{\otimes}
(\vartheta)\label{Rflip}
\end{equation}
for each $\vartheta\in\Gamma_{\inv}^{\otimes 2}$.
\end{lem}
\begin{pf}
Applying \eqref{R-reg} we find
\begin{equation}
m_\Omega\Bigl\{R_\omega\otimes\varphi\Bigr\}=m_\Omega\Bigl\{
\varphi\otimes R_\omega\Bigr\}\sigma
\end{equation}
for each $\varphi\in\tau(P)$. In particular, \eqref{Rflip} holds.
\end{pf}

 The  above  statement  implies  that  both  maps  $W$  and
$R_{\omega}^{\otimes}$   are
factorizable through the ideal $\cal{J}$ (so that $\cal{I}_{\cal{S}}$
is the natural domain for the Weil homomorphism).

Let us briefly consider some types of quantum principal bundles,
particularly interesting from the point of
view of the theory of characteristic classes.
\par\smallskip
{\it Bundles over Classical Smooth Manifolds}
\par\smallskip
Let us assume that $M$ is
a classical (compact) smooth manifold. According to \cite{D}
(locally trivial) quantum principal bundles $P$ over $M$ are
classified by
classical $G_{cl}$-bundles $P_{cl}$ (over the same manifold),
where $G_{cl}$ is
the classical part of $G$. Let us assume that $\Gamma$ is the
minimal admissible (bicovariant *-) calculus over $G$ (in the sense of
\cite{D}). Let $\Omega(P)$ be the graded-differential *-algebra
canonically associated to $P$ and $\Gamma$. Then characteristic
classes of $P$ are naturally interpretable
as (classical) characteristic classes of $P_{cl}$ (however, the
converse is generally not true).
\par\smallskip
{\it Quantum Line Bundles}
\par\smallskip
If $G=U(1)$ and if the calculus on $G$ is
1-dimensional then
there is essentially the only one
characteristic class given by
the curvature form. It corresponds to the Euler class in the classical
theory. The connection actually defines a ``global angular form''
on the bundle.

The $\circ$-structure is characterized by
$\zeta\circ z=\lambda\zeta$, where $\lambda\in\Re\setminus\{0\}$ and
$z\colon G\rightarrow\Bbb{C}$ is the canonical generator of $\cal{A}$.
The higher-order components are
trivial, $S_{\inv}^\wedge$ is generated by $\zeta\otimes\zeta$ (if
$\lambda\neq-1$) and $\cal{R}$
is generated by $z^{-1}+z/\lambda-(1+1/\lambda)$. Hence a
connection $\omega$ is
multiplicative iff $\omega(\zeta)^2=0$. The regularity condition
can be written as
$$\omega(\zeta)\varphi=
(-1)^{\partial\varphi}\lambda^k\varphi\omega(\zeta),$$
if $F^\wedge(\varphi)=\varphi\otimes z^k$.
\par\smallskip
{\it Special Differential Structures}
\par\smallskip
Let us assume that the calculus on $G$ is such that elements
$\pi(u_{ij})$ linearly generate $\Gamma_{\inv}$.
Here $u_{ij}$ are matrix elements of the
fundamental representation $u\colon\Bbb{C}^n\rightarrow\Bbb{C}^n
\otimes\cal{A}$ of $G$. Let us consider an arbitrary element
$\Delta\in\cal{I}_{\cal{S}}$. Then it is possible to construct a series
of characteristic classes, associated to coefficients
of the polynomial $p(\lambda,\Delta)=\xi_\lambda(\Delta)$, where
$\xi_\lambda\colon\cal{S}\rightarrow\cal{S}$ are automorphisms  specified
by
$$\xi_\lambda\pi(u_{ij})=\lambda\delta_{ij}1-\pi(u_{ij}),$$
if the above formula is compatible with the ideal $\cal{J}$ (in any case
$\xi_\lambda$ are acting in the tensor algebra $\Gamma_{\inv}^\otimes$).
These automorphisms intertwine the adjoint action $\adj_{\cal{S}}$
and hence preserve the algebra $\cal{I}_{\cal{S}}$.
In particular, if $\Delta$ is a quantum determinant (appropriately defined)
then coefficients of $p(\lambda,\Delta)$ define counterparts of classical
Chern classes.

Finally, if the structure group $G$ is classical
and if the calculus on $G$ is classical then the construction
of characteristic classes becomes the same as in the
classical case. In particular, regular connections graded-commute with
horizontal differential forms. However this gives a relatively
strong constraint for differential calculus on the bundle.

\section{Examples, Remarks $\&$ Some Additional Constructions}
\subsection{Infinitesimal Gauge Transformations A}
The *-$\cal V$-module  $\cal{E}=\tau^0 (P)$  of  tensorial
 0-forms  is
definable independently of the choice of a  differential  calculus
on  the  bundle  $P$.  The  elements  of  this  space  are   quantum
counterparts of infinitesimal gauge  transformations  (vertical
equivariant vector fields  on  the  bundle).  The
space $\cal{E}$ will be here analyzed from this point of view.

 Explicitly,   $\cal{E}$    is    consisting    of    linear    maps
$\zeta\colon \Gamma_{\inv}\rightarrow\cal{B}$ such that
the diagram
\begin{equation}\label{infg}
\begin{CD}
\Gamma_{\inv} @>{\mbox{$\zeta$}}>> \cal{B}\\
@V{\mbox{$\adj $}}VV   @VV{\mbox{$F$}}V\\
\Gamma_{\inv}\otimes\cal{A} @>>{\mbox{$\zeta\otimes \id $}}>
\cal{B}\otimes \cal{A}
\end{CD}
\end{equation}
is commutative.

 Let us observe that $\cal{E}$ is closed under  operations  $\langle,
\rangle $  and
$[,]$.   In   classical   geometry,  we have $[,]=-2\langle,\rangle$,
and $[,]\colon \cal{E}\times\cal{E}\rightarrow\cal{E}$  coincides   with
the standard commutator of vector fields (up to the sign).

     We are going to construct quantum analogs  of  contraction
operators associated to vector fields  representing  infinitesimal
gauge transformations.

 For each $\zeta\in\cal{E}$ let us consider  a  map
$\iota_{\zeta}\colon \Omega(P)\rightarrow\Omega(P)$ defined by
\begin{equation}
 \iota_\zeta(w)=-(-1)^{\partial w} \sum_k u_k \zeta(\eta_k )
\label{jota}
\end{equation}
where $\Sum_k u_k\otimes\eta_k =(\id \otimes\pi_{\inv}p_1)
\widehat{F}(w)$.
\begin{lem}
The diagram
\begin{equation}\label{jwF}
\begin{CD}
\Omega(P) @>{\mbox{$F^{\wedge}$}}>>
\Omega(P)\otimes\cal{A}\\
@V{\mbox{$\iota_\zeta$}}VV @VV{\mbox{$\iota_\zeta\otimes \id $}}V\\
\Omega(P) @>>{\mbox{$F^\wedge$}}>  \Omega(P)\otimes\cal{A}
\end{CD}
\end{equation}
is commutative.
\end{lem}
\begin{pf}
A direct computation gives
\begin{equation*}
\begin{split}
F^{\wedge}\iota_{\zeta}(w)&=-(-1)^{\partial w} \sum_i (\id \otimes
p_0)\widehat{F}(w_i)\e(a_i)F\zeta p_1(\vartheta_i)\\
&=-(-1)^{\partial w} \sum_i(w_i\otimes a_i)
(\zeta\otimes \id) \adj p_1(\vartheta_i)\\
&=-(-1)^{\partial w} (m_{\Omega}\otimes \id) (\id \otimes\zeta\pi_{\inv}p_1
\otimes p_0)\biggl(\sum_iw_i\otimes\widehat{\phi}(a_i \vartheta_i )\biggr)\\
&=-(-1)^{\partial w} (m_{\Omega}\otimes \id) (\id \otimes\zeta\pi_{\inv}p_1
\otimes p_0)(\widehat{F}\otimes \id) \widehat{F}(w)\\
&=(\iota_{\zeta}\otimes \id) F^{\wedge}(w),
\end{split}
\end{equation*}
where $\widehat{F}(w)=\Sum_iw_i \otimes a_i \vartheta_i$,
with $w_i \in\Omega(P)$, $a_i \in\cal{A}$ and $\vartheta_i
\in\Gamma^{\wedge}_{\inv}$.
\end{pf}

     The definition of $\iota_{\zeta}$  implies
\begin{equation}\label{h-linj}
\iota_\zeta(\varphi w)=(-1)^{\partial\varphi}\varphi\iota_\zeta(w),
\end{equation}
for each $\varphi\in\hor(P)$ and $w\in\Omega(P)$. In particular,
\begin{equation}
\iota_{\zeta}\bigl(\hor(P)\bigr)=\{0\},
\end{equation}
for each $\zeta\in\cal{E}$.

Let  us  consider linear maps
$\ell_{\zeta}\colon \Omega(P)\rightarrow\Omega(P)
$ given by
\begin{equation}\label{ld}
\ell_{\zeta} =d\iota_{\zeta} +\iota_{\zeta}d.
\end{equation}
These maps play the role of the corresponding Lie derivatives.
\begin{lem}\bla{i} The diagram
\begin{equation}\label{lwF}
\begin{CD}
\Omega(P) @>{\mbox{$F^{\wedge}$}}>> \Omega(P)\otimes\cal{A}\\
@V{\mbox{$\ell_\zeta$}}VV      @VV{\mbox{$\ell_\zeta\otimes \id $}}V\\
\Omega(P) @>>{\mbox{$F^\wedge$}}>  \Omega(P)\otimes\cal{A}
\end{CD}
\end{equation}
is commutative.

\smallskip
\bla{ii} The following equality holds
\begin{equation}
\ell_{\zeta}(\varphi)=\sum_k\varphi_k \zeta\pi(c_k),\label{lz}
\end{equation}
where $\varphi\in \hor(P)$ and
$\Sum_k\varphi_k \otimes c_k =F^{\wedge} (\varphi)$. In particular,
\begin{equation}
\ell_{\zeta}\bigl(\hor(P)\bigr)\subseteq \hor(P).
\end{equation}
\end{lem}
\begin{pf}
Diagram \eqref{lwF} follows  from \eqref{jwF}, \eqref{ld}
and \eqref{dwF}.
Identity \eqref{lz} directly follows from definitions of
$\iota_{\zeta}$   and $\ell_{\zeta}$:
\begin{multline*}
\ell_{\zeta}\varphi=\iota_{\zeta} d\varphi=(-1)^{\partial\varphi}
m_{\Omega} (\id \otimes\zeta\pi_{\inv}p_1)
\biggl[\sum_kd\varphi_k \otimes c_k +(-1)^{\partial\varphi}
\varphi_k \otimes dc_k \biggr] \\
=m_{\Omega}(\id \otimes\zeta\pi_{\inv})\biggl[\sum_k\varphi_k
\otimes dc_k
\biggr]=\sum_k\varphi_k \zeta\pi(c_k).\qed
\end{multline*}
\renewcommand{\qed}{}
\end{pf}

 Covariance properties \eqref{jwF} and \eqref{lwF}
imply  that
$\psi(P)$  is
invariant   under   compositions   with    $\iota_{\zeta}$     and
$\ell_{\zeta}$. In particular,
$\ell_{\zeta}\tau(P)\subseteq\tau(P)$ for each $\zeta\in\cal{E}$.
\begin{lem}\label{lem:lbra2}
We have
\begin{equation}
\ell_{\zeta} \varphi=[\varphi,\zeta],\label{lbra2}
\end{equation}
for each $\zeta\in\cal{E}$ and $\varphi\in\tau(P)$.
\end{lem}
\begin{pf}
It follows from \eqref{lz}, definitions of $c^\top$ and $[,\,]$,
and the tensoriality of $\varphi$.
\end{pf}

 Let us  compute  actions of $\iota_{\zeta}$ and $\ell_{\zeta}$
on connection forms.
\begin{lem}
We have
\begin{align}
\iota_{\zeta} \omega&=\zeta\label{jcon}\\
\ell_{\zeta}\omega&=d\zeta+[\omega,\zeta], \label{lcon}
\end{align}
for each $\zeta\in\cal{E}$ and $\omega\in \con(P)$.
\end{lem}
\begin{pf}
Both identities directly follow from property \eqref{dxFcon}
and
from definitions of $\iota_{\zeta}$  and $\ell_{\zeta}$.
\end{pf}

 In the framework of the theory presented in the previous paper, a
very  important  class  of  infinitesimal  gauge   transformations
naturally  appear.  These  transformations  can  be  described  as
infinitesimal generators (in the standard sense) of the  group  of
vertical automorphisms of the  bundle  $P$.  They  form  a  subspace
$\cal{G}\subseteq\cal{E}$. A more detailed geometrical  analysis
shows  that  the
elements  from  $\cal{G}$
are  naturally  identificable   with   standard
infinitesimal gauge transformations of the classical part $P_{cl}$
 of
the bundle $P$ (vertical  automorphisms  of  $P$  are  in  a  natural
bijection with standard gauge transformations of  $P_{cl}$).  Moreover,
the space $\cal{G}$ is closed under brackets $[,]$  and, in terms of the
mentioned identification $-[,]$ becomes the standard Lie  bracket  of
vector fields. The elements from $\cal{G}$ naturally act as derivations on
$\Omega(P)$. However, the action of the derivation generated  by  an
element  $\zeta\in\cal{G}$
generally  {\it differs}  from  the  action  of  the
corresponding Lie derivative $\ell_\zeta$ introduced in this subsection.
This is visible from \eqref{lcon}, which can be rewritten in the form
\begin{equation*}
\ell_{\zeta} \omega=D_{\omega} \zeta+[\omega,\zeta]+[\zeta,\omega].
\end{equation*}

 The last two  summands  generally  give  a  nontrivial  non-horizontal
contribution, even  in  the  case  $\zeta\in\cal{G}$.
Only  in  classical
geometry we have $[\omega,\zeta]+[\zeta,\omega]=0$  (more  generally
$[\varphi,\eta]=-(-1)^{\partial\varphi\partial\eta}
[\eta,\varphi]$  for   each
$\varphi,\eta\in\psi(P)$), due to the antisymmetricity  of $c^\top$,  and
the graded-commutativity of $\Omega(P)$.

\subsection{Infinitesimal Gauge Transformations B}
 Motivated by the  above  remarks,  a  slightly  different
approach to defining quantum analogs of the Lie  derivative  and
the contraction operator will be now presented.

 The main property of this approach is that in the special case of
bundles over smooth manifolds the
Lie derivative of an arbitrary element  $\zeta\in\cal{G}$
 coincides  with
the derivation generated by $\zeta$.

     We shall also introduce a general quantum counterpart of  the
space $\cal{G}$, and briefly analyze its properties.

\begin{lem}\label{lem:kj}
\bla{i} For each $\zeta\in\cal{E}$
there exists the unique $\varsigma_{\zeta}^\star\colon \vh(P)
\rightarrow\vh(P)$ such that
\begin{align}
\varsigma_{\zeta}^\star&\bigl(\hor(P)\bigr)=\{0\}\label{khor}\\
\varsigma_{\zeta}^\star&(w\vartheta)=\varsigma_{\zeta}^\star(w)\vartheta+
(-1)^{\partial w}w\zeta(\vartheta),\label{kprod}
\end{align}
for each $w\in \vh(P)$  and
$\vartheta\in\Gamma_{\inv}$.

\smallskip
\bla{ii} Similarly, for each $\zeta\in\cal{E}$ there exists the unique
$\iota_{\zeta}^\star\colon \vh(P)\rightarrow \vh(P)$
such that
\begin{align}
\iota_{\zeta}^\star(\vartheta\eta)&=
-\vartheta\iota_{\zeta}^\star(\eta)
+\sum_k\eta_k\zeta(\vartheta\circ c_k)\label{j2}\\
 \iota_{\zeta}^\star(\varphi\eta)&=(-1)^{\partial\varphi}\varphi
\iota^\star_{\zeta}(\eta),\label{j3}
\end{align}
for each $\varphi\in \hor(P)$, $\eta\in\Gamma_{\inv}^{\wedge}$ and
$\vartheta\in\Gamma_{\inv}$, where $\Sum_k\eta_k\otimes c_k
=\adj^{\wedge}(\eta)$. In particular,
\begin{equation}
\iota_{\zeta}^\star\bigl(\hor(P)\bigr)=\{0\}\label{j1}.
\end{equation}

 \bla{iii} The following identities hold
\begin{align}
F_{v\!h}\iota_{\zeta}^\star&=(\iota_{\zeta}^\star\otimes \id) F_{v\!h}
\label{wfj*}\\
F_{v\!h}\varsigma_{\zeta}^\star&=(\varsigma_{\zeta}^\star\otimes \id)
F_{v\!h}\label{wfk*}\\
\varsigma_{\zeta}^\star(w\vartheta)&=\varsigma_{\zeta}^\star(w)\vartheta+
(-1)^{\partial w}w\varsigma^\star_{\zeta}(\vartheta), \label{k*prod}
\end{align}
where $w\in \vh(P)$
and $\vartheta\in\Gamma_{\inv}^{\wedge}$.
\end{lem}
\begin{pf}
We shall prove ({\it i}\/) and properties \eqref{wfk*}
and \eqref{k*prod}. The  statements
about the map $\iota_{\zeta}^\star$  follow in a similar way.
It  is  clear  that
conditions \eqref{khor} and \eqref{kprod}
uniquely fix the values of $\varsigma_{\zeta}^\star$,
if it  exists.
Also, \eqref{k*prod} directly follows from \eqref{kprod}.
In order to establish the
existence of $\varsigma_{\zeta}^\star$, it is sufficient to check that
\eqref{kprod} is not in  a contradiction with the quadratic constraint
generating the ideal $S_{\inv}^{\wedge}$.

     For each $a\in\cal{R}$ we have
\begin{equation*}
\begin{split}
     \left\{\pi(a^{(1)})\pi(a^{(2)})\right\}&
\longrightarrow\left(\zeta\pi(a^{(1)})\pi(a^{(2)})-\pi(a^{(1)
}   )\zeta\pi(a^{(2)})\right)\\
&=\bigl(\zeta\pi(a^{(1)})\bigr)\pi(a^{(2)})\\
&\phantom{=}-\bigl(\zeta\pi(a^{(3)})\bigr)\pi(a
^{(1)})\circ \bigl\{\k(a^{(2)})a^{(4)}\bigr\}\\
&=\bigl(\zeta\pi(a^{(1)})\bigr)\pi(a^{(2)})\\
&\phantom{=}-\bigl(\zeta\pi(a^{(3)})\bigr)
\pi\left[\bigl(a^{(1)}-\e(a^{(1)})1\bigr)\k(a^{(2)})a^{(4)}\right]\\
&=\zeta\pi(a^{(2)})\pi\bigl(\k(a^{(1)})a^{(3)}\bigr)=0,
\end{split}
\end{equation*}
and hence $\varsigma_{\zeta}^\star$  exists.

 Finally, let us check \eqref{wfk*}. This equality is
 satisfied trivially
on elements from $\hor(P)$. The definition of $\zeta$ implies  that
it is satisfied on elements  from  $\Gamma_{\inv}$.
Finally, inductively applying \eqref{kprod} we conclude that
\eqref{wfk*} holds  on  the  whole  algebra $\vh(P)$.
\end{pf}

     Let us assume that the bundle $P$ admits  regular  connections,
as well as that $\J(P)=\{0\}$.
\begin{lem}
We have
\begin{equation}
\iota_{\zeta}=m_{\omega}\iota_{\zeta}^\star m_{\omega}^{-1}
\label{simj}
\end{equation}
for each $\zeta\in\cal{E}$ and $\omega\in\rc(P)$.
\end{lem}
\begin{pf}
Let  us  fix  a  connection $\omega\in\rc(P)$.
For each $\zeta\in\cal{E}$ let
$\iota'_{\zeta}\colon \vh(P)\rightarrow\vh(P)$ be a map defined by
$$\iota'_{\zeta}=m_{\omega}^{-1}\iota_{\zeta}m_{\omega}.$$
If $\varphi\in \hor(P)$  and
$\eta\in\Gamma_{\inv}^{\wedge}$ then
\begin{equation*}
\iota'_\zeta(\varphi\eta)=m_{\omega}^{-1}\iota_{\zeta}
\bigl(\varphi\omega^{\wedge}(\eta)\bigr)
=(-1)^{\partial\varphi}m_{\omega}^{-1}\bigl(\varphi\iota_{\zeta}
\omega^{\wedge}(\eta)\bigr)=(-1)^{\partial\varphi}\varphi
\iota'_{\zeta}\eta
\end{equation*}
according to \eqref{h-linj}.
Further, if $\vartheta\in\Gamma_{\inv}$ then
\begin{equation*}
\begin{split}
\iota'_\zeta(\vartheta\eta)=m_{\omega}^{-1}\iota_{\zeta}
\bigl(\omega(\vartheta)
\omega^{\wedge}(\eta)\bigr)&=(-1)^{\partial\eta}
m_{\omega}^{-1}m_{\Omega}(\id \otimes\zeta\pi_{\inv}p_1
)\widehat{F}\bigl[\omega(\vartheta)\omega^{\wedge}(\eta)\bigr]\\
&=m_{\omega}^{-1}m_{\Omega}
(\id \otimes\zeta)\Bigl(\sum_k
\omega^{\wedge}(\eta_k )
\otimes\vartheta\circ c_k \Bigr)\\
&\phantom{=}+(-1)^{\partial\eta}
m_\omega^{-1} \Bigl\{\omega(\vartheta)m_{\Omega}\Bigl[
(\id \otimes\zeta\pi_{\inv}p_1)
\widehat{F}\omega^{\wedge}(\eta)\Bigr]\Bigr\}\\
&=\Sum_k\eta_k \zeta(\vartheta\circ c_k)-m_{\omega}^{-1}
\bigl(\omega(\vartheta)\iota_{\zeta}\omega^{\wedge}(\eta)\bigr)\\
&=\sum_k\eta_k \zeta(\vartheta\circ c_k)-\vartheta
\iota'_{\zeta}(\eta).
\end{split}
\end{equation*}

 Here, $\Sum_k\eta_k \otimes c_k =\adj^{\wedge}(\eta)$
and we have  used  elementary
properties of entities figuring in the game. Applying  ({\it ii}\/)
Lemma~\ref{lem:kj}
we find that $\iota'_{\zeta}=\iota_{\zeta}^\star$. Hence \eqref{simj}
holds.
\end{pf}

 For    each    $\omega\in\rc(P)$    and $\zeta\in\cal{E}$ let
$\varsigma_{\zeta,\omega}  \colon \Omega(P)\rightarrow\Omega(P)$  be  a
map  introduced via the diagram
\begin{equation}\label{defkW}
\begin{CD}
\vh(P) @>{\mbox{$\varsigma^\star_{\zeta}$}}>>\vh(P)\\
@V{\mbox{$m_\omega$}}VV  @VV{\mbox{$m_\omega$}}V\\
\Omega(P) @>>{\mbox{$\varsigma_{\zeta,\omega}$}}> \Omega(P)
\end{CD}
\end{equation}
and let $\ell_{\zeta,\omega}\colon \Omega(P)\rightarrow\Omega(P)$
be a map given by
\begin{equation} \label{deflW}
\ell_{\zeta,\omega}=d\varsigma_{\zeta,\omega}+\varsigma_{\zeta,\omega}
d.
\end{equation}

     It is evident that $\ell_{\zeta,\omega}$ and
$\varsigma_{\zeta,\omega}$    are right-covariant maps,  in
the sense that
\begin{align}
F^{\wedge}\ell_{\zeta,\omega}&=(\ell_{\zeta,\omega}\otimes \id) F^{\wedge}
\label{lWcov}\\
F^{\wedge}\varsigma_{\zeta,\omega}&=(\varsigma_{\zeta,\omega}
\otimes \id) F^{\wedge}.\label{kWcov}
\end{align}

     The maps $\ell_{\zeta,\omega}$ and $\varsigma_{\zeta,\omega}$,
are also interpretable as quantum
counterparts  of  the  Lie  derivative  and  the   contraction
operator respectively. In contrast  to  the  classical  case,  these
maps  are generally connection-dependent. However,
\begin{lem}\label{lem:k1j}
If $w\in\Omega_1(P)$ then
\begin{equation}
\varsigma_{\zeta,\omega}(w)=\iota_{\zeta}(w),\label{k1j}
\end{equation}
for each $\zeta\in\cal{E}$ and $\omega\in\rc(P)$. In  particular,
operators  $\ell_{\zeta}$ and $\ell_{\zeta,\omega}$
possess the same restrictions on $\hor(P)$.
\end{lem}
\begin{pf}
It follows from the fact that $\iota^\star_{\zeta}$  and
$\varsigma^\star_{\zeta}$ coincide  on the spaces $\hor(P)$ and
$\hor(P)\otimes\Gamma_{\inv}$.
\end{pf}

Covariance properties \eqref{lWcov}--\eqref{kWcov}
enable us to  define  actions  of $\ell_{\zeta,\omega}$
and $\varsigma_{\zeta,\omega}$ in the space $\psi(P)$.

 The $\omega$-dependence of constructed operators becomes explicitly
visible if we consider the action of $\ell_{\zeta,\omega}$
on connection forms.
\begin{lem} We have
\begin{equation}
\ell_{\zeta,\omega}\tau=D_{\tau}\zeta+
[\tau-\omega,\zeta]+[\zeta,\tau-\omega]\label{wdep}
\end{equation}
for each  $\zeta\in\cal{E}$, $\omega\in\rc(P)$  and
$\tau\in \con(P)$.  In
particular $\ell_{\zeta,\omega}\tau$ is always tensorial.
\end{lem}
\begin{pf}
Using Lemmas~\ref{lem:lbra2} and~\ref{lem:k1j}, and properties
\eqref{jcon} and \eqref{deflW}
we obtain
$$
\ell_{\zeta,\omega}\tau=\ell_{\zeta,\omega}\omega
+[\tau-\omega,\zeta]=d\zeta+[\tau-\omega,\zeta]+
\varsigma_{\zeta,\omega}d\omega. $$
On the other hand, \eqref{defR} together with the regularity of
$\omega$, tensoriality of $\zeta$ and the definition of
$\varsigma_{\zeta,\omega}$ gives
\begin{equation*}
\begin{split}
\varsigma_{\zeta,\omega}
d\omega(\vartheta)&=\varsigma_{\zeta,\omega}
(\langle \omega,\omega\rangle +R_{\omega})(\vartheta)=
\varsigma_{\zeta,\omega}
\langle \omega,\omega\rangle (\vartheta)\\
&=-\zeta\pi(a^{(1)})\omega\pi(a^{(2)})
+\omega\pi(a^{(1)})\zeta\pi(a^{(2)})\\
&=-\zeta\pi(a^{(1)})\omega\pi(a^{(2)})+
\zeta\pi(a^{(3)})\omega\left[\pi(a^{(1)})\circ \bigl(
\k(a^{(2)})a^{(4)}\bigr)\right]\\
&=-\zeta\pi(a^{(2)})\omega\pi\bigl(\k(a^{(1)})a^{(3)}\bigr)
=-[\zeta,\omega](\vartheta),
\end{split}
\end{equation*}
where $a\in\cal{A}$ satisfies \eqref{delpi}. Consequently,
$$
\ell_{\zeta,\omega}\tau=D_{\omega}\zeta+
[\tau-\omega,\zeta]=D_{\tau}
\zeta+[\tau-\omega,\zeta]+[\zeta,\tau-
\omega].\qed$$
\renewcommand{\qed}{}
\end{pf}

 Finally, let $\cal{G}\subseteq \cal{E}$ be the
space of elements $\zeta$ satisfying
\begin{equation}
\zeta(\vartheta)\varphi=\sum_k\varphi_k   \zeta(\vartheta\circ
c_k)\label{regg}
\end{equation}
for each $\varphi\in \hor(P)$ and $\vartheta\in\Gamma_{\inv}$.
Let us assume  that
$\cal{G}$ is nontrivial.
\begin{pro}\bla{i} The space
$\cal{G}$ is closed under the action
of brackets $[,]$.
We have
\begin{equation}
\zeta\pi(a^{(1)})\xi\pi(a^{(2)})-\xi\pi(a^{(1)})\zeta
\pi(a^{(2)})=[\zeta,\xi]\pi(a)\label{lieg}
\end{equation}
for each $\zeta,\xi\in\cal{G}$ and
$a\in\cal{A}.$
In particular, brackets $[,]$ determine a Lie algebra structure
on $\cal{G}$.

\smallskip
 \bla{ii} Operators $\varsigma_{\zeta,\omega}$ and $\ell_{\zeta,\omega}$
are   $\omega$-independent, if $\zeta\in\cal{G}$.

\smallskip
    \bla{iii} The following identities hold
\begin{gather}
\ell^\sstar_{\zeta}(wu)=\ell^\sstar_{\zeta}(w)u+w\ell^\sstar_{\zeta}(u)
\label{lwu}\\
\varsigma_{\zeta}(wu)=\varsigma_{\zeta}(w)u+(-1)^{\partial w}w\varsigma
_{\zeta}(u) \label{kwu}\\
\varsigma_{\zeta}\varsigma_{\xi}+\varsigma_{\xi}\varsigma_{\zeta}=0
\label{antcom}\\
\ell^\sstar_{\zeta}\ell^\sstar_{\xi}-\ell^\sstar_{\xi}
\ell^\sstar_{\zeta}=-\ell^\sstar_{[\zeta,\xi]}\label{coml}\\
\ell^\sstar_{\zeta}\varsigma_{\xi}-\varsigma_{\xi}
\ell^\sstar_{\zeta}=-\varsigma_{[\zeta,\xi]}.
\label{comlk}
\end{gather}
Here $\zeta,\xi\in\cal{G}$ and $w,u\in\Omega(P)$ while
$\varsigma_{\zeta}=\varsigma_{\zeta,\omega}$ and
$\ell^\sstar_\zeta=\ell_{\zeta,\omega}$.
\end{pro}
\begin{pf}
We compute
\begin{multline*}
\xi\pi(a^{(1)})\zeta\pi(a^{(2)})=\zeta\pi(a^{(3)})\xi\bigl[\pi(a
^{(1)})\circ\bigl(\k(a^{(2)})a^{(4)}\bigr)\bigr]\\
=\zeta\pi(a^{(1)})\xi\pi(a^{(2)})-[\zeta,\xi]\pi(a).
\end{multline*}

     Let us check that $\cal{G}$ is closed under the brackets
$[,]$.  Using properties \eqref{lwF}--\eqref{lz} and
\eqref{regg}--\eqref{lieg} we find
\begin{equation*}
\begin{split}
([\zeta,\xi]\pi(a))\varphi&=\sum_k \varphi_k
\zeta\bigl(\pi(a^{(1)})\circ c_k^{(1)}\bigr)\xi
\bigl(\pi(a^{(2)})\circ c_k^{(2)}\bigr)\\
&\phantom{=}-\sum_k\varphi_k\xi\bigl(\pi(a^{(1)})\circ c_k^{(1)}\bigr)
\zeta
\bigl(\pi(a^{(2)})\circ c_k^{(2)} \bigr)\\
&=\sum_k\left(\varphi_k\zeta\pi(a^{(1)}c_k^{(1)})\xi\pi
(a^{(2)}c_k^{(2)})-
\varphi_k\xi\pi(a^{(1)}c_k^{(1)})\zeta\pi(a^{(2)}c_k
^{(2)})\right)\\
&\phantom{=}-\sum_k \ell_\zeta(\varphi_k)\xi\pi(ac_k)-
\sum_k\varphi_k \zeta\pi(ac^{(1)}_k)\xi\pi(c_k^{(2)})\\
&\phantom{=}+\sum_k\ell_\xi(\varphi_k)\zeta\pi(ac_k
)+\sum_k\varphi_k \xi\pi(ac_k^{(1)})\zeta\pi(c_k^{(2)})\\
&=\sum_k \varphi_k[\zeta,\xi]
\pi(ac_k)-\sum_k\xi\pi(a)\ell_\zeta(\varphi)
-\sum_k\zeta\pi(a)\varphi_k\xi\pi(c_k)\\
&\phantom{=}+\sum_k\zeta\pi(a)\ell_\xi(\varphi)+
\sum_k\xi\pi(a)\varphi_k\zeta\pi(c_k)\\
&=\sum_k\varphi_k[\zeta,\xi]\pi(ac_k)=\sum_k\varphi_k[\zeta,\xi]
\bigl(\pi(a)\circ c_k\bigr),
\end{split}
\end{equation*}
where $a\in \ker(\e)$ and $\varphi\in \hor(P)$.

Now we shall prove identities \eqref{lwu}--\eqref{comlk}.
Let  us  observe  that \eqref{lwu} directly follows from \eqref{kwu}
and \eqref{deflW}. On the other  hand
the fact that $\varsigma_{\zeta,\omega}$ is an antiderivation together
with  Lemma~\ref{lem:k1j}
shows     that     $\varsigma_{\zeta,\omega}$     (and      therefore
$\ell_{\zeta,\omega}$)  is $\omega$-independent.
Evidently, \eqref{kwu} is    equivalent    to
the    fact    that  $\varsigma_{\zeta}^\star$   is  an
antiderivation on $\vh(P)$.
Having in mind  identity
\eqref{k*prod} and property \eqref{khor} it
is sufficient to check that
$$\varsigma_{\zeta}^\star(\vartheta\varphi)=
\varsigma_{\zeta}^\star(\vartheta)
\varphi, $$ for each $\vartheta\in\Gamma_{\inv}^\wedge$ and
$\varphi\in \hor(P)$. However,  this
easily  follows  from  property  \eqref{regg} and  the  definition   of
$\varsigma_{\zeta}^\star$.

 For each $\zeta,\xi\in \cal{G}$ the anticommutator of
$\varsigma_{\zeta}^\star$   and $\varsigma^\star_{\xi}$   is  an
antiderivation  on  $\vh(P)$.
This   anticommutator
vanishes  on  $\hor(P)$  and   $\Gamma_{\inv}$.
Therefore   it   vanishes identically.

 Having  in  mind  that  $\ell_{\zeta}^\sstar$ are  derivations  on
$\Omega(P)$
commuting with the differential, and that derivations form  a  Lie
algebra, it is sufficient to check that \eqref{coml}
holds  on  elements
$b\in\cal{B}$. We have
\begin{equation*}
\begin{split}
(\ell^\sstar_{\zeta}\ell^\sstar_{\xi}-\ell^\sstar_{\xi}
\ell^\sstar_{\zeta})(b)&=\ell^\sstar_{\zeta}\sum_kb_k
\xi\pi(a_k)-\ell^\sstar_{\xi}\sum_kb_k\zeta\pi(a_k)\\
&=\sum_k \left[b_k \xi\pi(a^{(1)}_k)\zeta\pi(a_k^{(2)})-
b_k\zeta\pi(a^{(1)}_k)\xi\pi(a_k^{(2)})\right]\\
&=\sum_kb_k [\xi,\zeta]\pi(a_k)=\ell^\sstar_{[\xi,\zeta]}(b),
\end{split}
\end{equation*}
where $F(b)=\Sum_kb_k \otimes a_k$.
 Similarly, it is sufficient to check that  \eqref{comlk}
holds  on elements of the form $\omega(\vartheta)$. We have
$$\left(\ell^\sstar_{\zeta}\varsigma_{\xi}-\varsigma_{\xi}
\ell^\sstar_{\zeta}\right)\omega=
\ell^\sstar_{\zeta}\xi-\varsigma_{\xi}D_{\omega}\zeta=[\xi,\zeta]=
-\varsigma_{[\zeta,\xi]}\omega,$$
which completes the proof.
\end{pf}
\subsection{Some Interrelations}
A similar  approach  to  general  quantum  principal  bundles  is
presented in \cite{BM}. Let us briefly consider interrelations
between the \cite{BM} and the theory developed here.
At the level of spaces both formulations coincide (modulo the
*-structure).  The  main  difference  appears  at  the  level   of
differential calculus on the bundle.
As first, it  is assumed in \cite{BM} that
the higher-order differential calculus uniquely follows  from  the
first-order  one,  beeing  based  on  (an  appropriate)  universal
envelope of  the  first-order  differential  structure.  Secondly,
the total ``pull back''  of  the  right  action  does  not  figure
in \cite{BM} and horizontal forms are defined in a
different way.

More precisely, let $P=(\cal B,i,F)$ be a quantum principal
$G$-bundle over $M$ and let $\Omega(P)$ be  an  arbitrary
graded-differential *-algebra satisfying properties ({\it diff1/2\/}).
Let $\Omega_{\hr}\subseteq \hor(P)$ be a (*-) subalgebra generated by
$\cal B$ and $d\bigl(i(\cal  V)\bigr)$.  This  algebra  is
a  counterpart  of horizontal forms introduced in \cite{BM}.
One  of  the main  conditions postulated in \cite{BM} (in the
definition of differential calculus) is that the sequence
\begin{equation*}
0\rightarrow\Omega_{\hr}^1
\hookrightarrow\Omega^1(P) @>{\pi_v}>> \ver^1(P)\rightarrow
0\end{equation*}
is exact. According to Lemma~\ref{lem:seq} this is equivalent to
\begin{equation}\label{bm-md}
\Omega^1_{\hr}=\hor^1(P).
\end{equation}
This condition can  be  understood
as a (relatively strong) condition {\it for the bundle},
if applied to the universal case, for example. Namely, {\it a  trivial}
differential
calculus on the bundle can be always  constructed  by  taking  the
universal differential  envelope  of  $\cal  B$  (conditions  ({\it
diff1/2\/}) hold). However  \eqref{bm-md}  does  not  generally  hold
(although it holds in various interesting special cases).
In particular, not all {\it quantum homogeneous spaces}
(endowed with universal differential calculus) can be
included in the theory presented in \cite{BM}.

It is worth noticing that if \eqref{bm-md} holds and if the bundle
admits regular and multiplicative connections then $\Omega_{\hr}=
\hor(P)$. Also, two horizontal algebras coincide if $\Omega(M)$
is generated (as a differential algebra) by $i(\cal{V})$. This
follows from the fact that
each $\varphi\in \hor(P)$ can be written in the form
\begin{equation*}
\varphi=\sum_i w_ib_i
\end{equation*}
where $b_i\in\cal B$ and $w_i\in\Omega(M)$
(as explained in Appendix B).

The definition of connection forms given in \cite{BM} is (modulo
the *-structure and  differencies  between  differential  calculi)
equivalent to the definition proposed in this work.
On the other hand  the  embedded differential map $\delta$ does
not figure in the formalism of connections. In particular,
operators of covariant derivative and curvature described in
\cite{BM} are generally different from $D_\omega$ and $R_\omega$
constructed here.

\subsection{Trivial Bundles}
For given quantum space $M$ and compact matrix quantum group
$G$ let us  define
a  *-algebra  $\cal{B}$  and  maps
$F\colon \cal{B}\rightarrow\cal{B}\otimes\cal{A}$ and
$i\colon \cal{V}\rightarrow\cal{B}$ by
\begin{gather*}
\cal{B}=\cal{V}\otimes\cal{A}\\
i(\,)=(\,)\otimes 1 \\
\id\otimes \phi=F.
\end{gather*}

 The triplet $P=(\cal{B},i,F)$ is a {\it trivial} quantum principal bundle
over $M$. Geometrically $P=M\times G$.

The algebra of verticalized forms is isomorphic  to  the
tensor product
\begin{equation*}
\ver(P)=\cal{V}\otimes\Gamma^{\wedge},\label{triv4}
\end{equation*}
with $d_v\leftrightarrow\id \otimes d$.

 Let $\Omega(M)$ be an arbitrary graded-differential *-algebra
generated by $\cal{V}=\Omega^0(M)$,
representing  a  differential  calculus on $M$. Then  it  is
natural  to  define  the  algebra  $\Omega(P)$
(representing a differential calculus on the bundle $P$) as the graded
tensor product
\begin{equation*}
 \Omega(P)=\Omega(M)\grten\Gamma^{\wedge}.
\end{equation*}
We have then
\begin{equation*}
\widehat{F}=\id \otimes\widehat{\phi}.
\end{equation*}
Horizontal forms constitute a *-subalgebra
\begin{equation*}
\hor(P)=\Omega(M)\otimes\cal{A}.
\end{equation*}

     We shall  now  analyze  in  more  details  the  structure  of
tensorial forms, connection  forms,  and  operators  of  covariant
derivative and curvature, in the special case of trivial  bundles.

Let $\cal{X}$ be the graded *-$\Omega(M)$ module of linear maps
$L\colon \Gamma_{\inv}\rightarrow\Omega(M)$.
\begin{lem}
For each $\varphi\in\tau(P)$ there exists
the unique $L^\varphi\in\cal{X}$ such that
\begin{equation}
\varphi(\vartheta)=(L^\varphi\otimes \id) \adj (\vartheta)
\label{triv9}
\end{equation}
for  each  $\vartheta\in\Gamma_{\inv}$.
The  above   formula
establishes an isomorphism between  graded
*-$\Omega(M)$-modules.
\end{lem}

\begin{pf}
For a given $L\in\cal{X}$ let $\varphi_L\colon \Gamma_{\inv}
\rightarrow\Omega(P)$ be a map
determined by equality $\varphi_L(\vartheta)=(L\otimes
\id) \adj(\vartheta) $. Evidently, the image of
$\varphi_L$  is contained in $\hor(P)$ and
$$
F^{\wedge}\varphi_L=(L\otimes\phi)\adj =\bigl[(L\otimes \id)
\adj \otimes \id \bigr]\adj =
(\varphi_L\otimes \id) \adj . $$
Hence $\varphi_L\in\tau(P)$.
It is clear that the map $L\mapsto\varphi_L$ is a
monomorphism of *-$\Omega(M)$-modules
(because $L=(\id \otimes \e)\varphi_L$).

 Let us consider an  arbitrary  tensorial  form  $\varphi$.
Acting by $\id\otimes\e\otimes \id $ on the tensoriality identity
for $\varphi$ we find that \eqref{triv9} holds, with
$L^\varphi=(\id \otimes \e)\varphi$.
\end{pf}

     The following lemma gives a similar description of connection
forms.
\begin{lem}\bla{i} The formula
\begin{equation}
\omega(\vartheta)=(A^\omega\otimes \id) \adj (\vartheta)+
1\otimes\vartheta \label{triv10}
\end{equation}
establishes a bijective affine correspondence between  connections
on $P$ and hermitian elements of $\cal{X}^1$.

\smallskip
    \bla{ii} A connection $\omega$ is regular iff
\begin{gather}
A^\omega(\vartheta)\zeta=(-1)^{\partial\zeta}\zeta
A^\omega(\vartheta)
\label{triv11}\\
A^\omega(\vartheta\circ a)=\e(a)A^\omega(\vartheta)\label{triv12}
\end{gather}
for each $a\in\cal{A}$, $\vartheta\in\Gamma_{\inv}$
and $\zeta\in\Omega(M)$.
\end{lem}
\begin{pf}
The formula
$\omega_0(\vartheta)=1\otimes\vartheta$
determines a canonical ``flat'' connection on $P$. The  statement
({\it i}\/) follows from the  previous  lemma  and  the  fact  that
hermitian elements of $\tau^1(P)$ form the vector space associated to
$\con(P)$.

     Let us assume that $\omega\in\rc(P)$. In other words
$$ \omega(\vartheta)(\zeta\otimes a)=(-1)^{\partial \zeta}
(\zeta\otimes a^{(1)})\omega(\vartheta\circ a^{(2)})
$$
for each $\zeta\in\Omega(M)$, $a\in\cal{A}$ and
$\vartheta\in\Gamma_{\inv}$. This is
equivalent to
\begin{equation}
\sum_k A^\omega(\vartheta_k)\zeta\otimes c_ka=
(-1)^{\partial \zeta}\sum_k \zeta A^\omega(\vartheta_k \circ a^{(1)})
\otimes c_ka^{(2)},\label{triv13}
\end{equation}
where $\Sum_k\vartheta_k \otimes c_k =\adj (\vartheta).$
Acting by $\id\otimes \e$ on both sides on the above equality we obtain
\begin{equation}
 A^\omega(\vartheta)\zeta\e(a)=(-1)^{\partial
\zeta}\zeta A^\omega(\vartheta\circ a),
\label{triv14}
\end{equation}
which is equivalent to conditions listed in ({\it ii}).  Conversely
\eqref{triv14} imply \eqref{triv13}, evidently.
\end{pf}

     The   bijection $\omega\leftrightarrow A^\omega$
generalizes    the  classical  correspondence
between connections and their gauge potentials. In the  previous
paper, a similar correspondence  was  established,  at  the  local
level. This was  possible  because  of  the local  triviality  of
considered bundles. However, in the general quantum context it  is
not  possible  to  speak  about  local  domains  on   the   base
space, and hence it is not possible to speak about locally trivial
bundles.

\begin{lem} We have
\begin{equation}
R_{\omega}(\vartheta)=(F^\omega\otimes\id)\adj(\vartheta) \label{triv15}
\end{equation}
where $F^\omega\in \cal{X}^2$  is a hermitian element given by
\begin{equation}
F^\omega=dA^\omega-\langle A^\omega,A^\omega\rangle .\label{triv16}
\end{equation}

     Further, if $\varphi\in\hor(P)$ then
\begin{equation}
D_{\omega}(\varphi)\leftrightarrow q^{\omega,\varphi}
\label{triv17}
\end{equation}
where
\begin{equation}
q^{\omega,\varphi}=dL^\varphi-
(-1)^{\partial\varphi}[L^\varphi,A^\omega].\label{triv18}
\end{equation}
\end{lem}
\begin{pf}
Inserting \eqref{triv9} and \eqref{triv10} in \eqref{Dten}
we obtain
\begin{equation*}
\begin{split}
(D_{\omega}\varphi)(\vartheta)&=\sum_kdL^\varphi(\vartheta_k)\otimes
c_k +(-1)^{\partial \varphi}\sum_k L^\varphi(\vartheta_k)\otimes dc_k \\
&\phantom{=}-(-1)^{\partial \varphi}\sum_k
\left[L^\varphi(\vartheta_k)\otimes
c_k^{(1)}\pi(c_k^{(2)})+
L^\varphi(\vartheta_k)A^\omega\pi(c_k^{(1)})\otimes
c_k^{(2)}\right]\\
&=\bigl[\bigl(dL^\varphi-(-1)^{\partial \varphi}[L^\varphi,A^\omega]\bigr)
\otimes \id\bigr]\adj (\vartheta).
\end{split}
\end{equation*}
Similarly, inserting \eqref{triv10} in \eqref{defR} we obtain
\begin{equation*}
\begin{split}
R_{\omega}(\vartheta)&=\sum_kdA^\omega(\vartheta_k)\otimes
c_k
-\sum_kA^\omega(\vartheta_k)\otimes dc_k +1\otimes d\vartheta\\
&\phantom{=}+1\otimes\pi(a^{(1)})\pi(a^{(2)})+
A^\omega\pi(a^{(2)})\otimes \k(a^{(1)})a^{(3)}\pi(a^{(4)})\\
&\phantom{=}-A^\omega\pi(a^{(3)})\otimes\pi(a^{(1)})\k(a^{(2)})a^{(4)}\\
&\phantom{=}+A^\omega\pi(a^{(2)})A^\omega\pi(a^{(3)})\otimes
\k(a^{(1)})a^{(4)}
=\bigl[\bigl(dA^\omega-\langle A^\omega,A^\omega
\rangle \bigr)\otimes \id \bigr]
\adj (\vartheta).
\end{split}
\end{equation*}
Here, $a\in \ker(\e)$ is such that \eqref{delpi} holds.
\end{pf}

Infinitesimal gauge transformations are  in  a  natural  bijection
with linear maps $\gamma\colon \Gamma_{\inv}\rightarrow\cal{V}$.
The elements of  $\cal{G}$ correspond  to
functions $\gamma$ satisfying
\begin{gather*}
 \gamma(\vartheta\circ a)=\e(a)\gamma(\vartheta)\\
w\gamma(\vartheta)=\gamma(\vartheta)w
\end{gather*}
for each $\vartheta\in\Gamma_{\inv}$, $a\in\cal{A}$ and
$w\in\Omega(M)$.

 Let  us  assume  that  $\Gamma$ is the minimal admissible
(bicovariant *-) calculus over $G$ (in the sense of \cite{D}). The
following  natural  identifications hold
\begin{align*}
\cal{G}&\leftrightarrow {\cal Z}^0(M)\otimes \lie(G_{cl})\\
\rc(P)&\leftrightarrow {\cal Z}^1(M)\otimes \lie(G_{cl}),
\end{align*}
where $G_{cl}$ is the classical part of $G$ and $\cal{Z}(M)$ is the
(graded) centre of $\Omega(M)$.

\subsection{Quantum Homogeneous Spaces}

Let $H$ be a compact matrix quantum group and let
$G$ be a (compact) subgroup of $H$.
At the formal level, this  presumes  a
specification of a *-epimorphism (the corresponding
``restriction map'') $j\colon \cal B\rightarrow\cal{A}$ such that
\begin{gather*}
(j\otimes j)\phi'=\phi j\\
\e j=\e'\\
\k j=j\k'.
\end{gather*}
Here $\cal B$ is the functional Hopf *-algebra for $H$.
In what follows
entities endowed with the prime will reffer to $H$.

     The *-homomorphism
$F\colon \cal B\rightarrow\cal B\otimes\cal A$ given by
\begin{equation*}
F=(\id \otimes j)\phi'
\end{equation*}
is interpretable as the right action of $G$ on $H$.
Let $M$  be  the
corresponding ``orbit space''. At the formal level, $M$
is represented by the fixed-point *-subalgebra
$\cal{V}\subseteq\cal B$.  Let  $i\colon \cal V\hookrightarrow\cal B$
be the inclusion map.
\begin{lem}
The  triplet  $P=(\cal{B},i,F)$
is a quantum principal $G$-bundle over $M$.
\end{lem}
\begin{pf}
It is evident that conditions ({\it qpb1/2}\/) of Definition 3.1 are
satisfied. We have
$$1\otimes j(b)=\k(b^{(1)})F(b^{(2)}), $$
for each $b\in\cal B$. Hence ({\it qpb4}\/) holds.
\end{pf}

Because of the inclusion $\phi'(\cal{V})\subseteq\cal{B}
\otimes\cal{V}$
there exists a natural left action of $H$ on $M$, defined by
$\phi'i\colon \cal{V}\rightarrow\cal{B}\otimes \cal V$.   This
action is ``transitive'' in the  sense  that  only  scalar  elements
of $\cal{V}$ are invariant.
In this sense $M$ is understandable as a quantum homogeneous $H$-space.

Now a construction of a differential calculus on $H$ will be
presented,  which  explicitly  takes  care  about  the
``fibered'' geometrical framework.

 Let $\Psi$ be a left-covariant first-order *-calculus over
$H$ and   let $\cal{R}'\subseteq \ker(\e')$   be    the
corresponding right
$\cal{B}$-ideal. Let us assume that
\begin{align}
j(\cal{R}')&\subseteq\cal{R}\label{r-incl}\\
(\id \otimes j)\ad^\prime({\cal R}^\prime)&\subseteq {\cal R}^\prime
\otimes\cal A\label{r-bicov}
\end{align}
where $\cal{R}\subseteq \ker(\e)$ is the right $\cal{A}$-ideal
which  determines the   calculus
$\Gamma$    over $G$. Condition \eqref{r-incl} ensures the existence of
the projection map $\rho\colon\Psi_{\inv}\rightarrow\Gamma_{\inv}$, which
is determined by the formula
\begin{equation}
\rho\pi^\prime=\pi j.
\end{equation}
The meaning of condition \eqref{r-bicov} is that the calculus $\Psi$ is
right-covariant, relative to $G$. Consequently, there exists the
corresponding adjoint action
$\chi\colon\Psi_{\inv}\rightarrow\Psi_{\inv}\otimes\cal{A}$. This map is
explicitly given by
 \begin{equation}
\chi\pi^\prime=(\pi^\prime\otimes j)\ad^\prime.
\end{equation}

 Maps $\rho$ and $\chi$ are hermitian and
\begin{gather}
\rho(\vartheta\circ a)=\rho(\vartheta)\circ j(a)\\
\adj \rho=(\rho\otimes \id) \chi.
\end{gather}
In  particular,  the  space  $\L=\ker(\rho)$   is   a   *-   and
$\chi$-invariant submodule of $\Psi_{\inv}$.

Let us now assume that the full calculus on the bundle $P$
($\Leftrightarrow$ the fibered $H$) is described by  a
graded-differential  *-algebra  $\Omega(P)$ built   over
$\Psi$ which is such  that  the  map  $F$ (and therefore
$\chi$) can  be  extended  to a differential  algebra
homomorphism
$\widehat{F}\colon \Omega(P)\rightarrow\Omega(P)
\grten\Gamma^\wedge$
(that is, property ({\it diff2\/}) holds).
Let us consider a
*-invariant
and $\chi$-invariant complement $\L^\perp$ of $\L$
in $\Psi_{\inv}$.

\begin{lem}
A linear map $\omega\colon \Gamma_{\inv}\rightarrow\Omega(P)$
given by
\begin{equation}
\omega(\vartheta)=\left(\rho\restr \L^\perp\right)^{-1}(
\vartheta)
\end{equation}
is a connection on $P$.
\end{lem}
\begin{pf}  By  construction,  it  follows  that  $\omega$  is   a
hermitian   pseudotensorial   1-form.   Condition
\eqref{pivcon}
directly follows from the observation that
\begin{equation*}
\pi_v(\vartheta)=1\otimes \rho(\vartheta)
\end{equation*}
for each $\vartheta\in \Psi_{\inv}.$
\end{pf}

Let us assume that the subspace $\L^\perp$ is also  a  submodule
of $\Psi_{\inv}$. In other words
\begin{equation}
\Proj(\vartheta\circ b)=\Proj(\vartheta)\circ b,
\end{equation}
where $\Proj\colon \Psi_{\inv}\rightarrow \L$
is the projection  map, corresponding to
the splitting
$$\Psi_{\inv}=\L\oplus \L^\perp.$$
In what follows we shall identify the spaces $\L^\perp$ and
$\Gamma_{\inv}$, via the map $\rho$. The right
$\cal B$-module structure on $\L^\perp$  can  be
naturally ``projected'' to  the right $\cal A$-module
structure on this space,
so that $\rho\leftrightarrow\L^\perp$ becomes  a  right
$\cal A$-module isomorphism,
because of $$\L^\perp\circ \ker(j)=\{0\}.$$

Further, let  us  assume  that  $\Omega(P)$   is
left-covariant and   let
$\Lf\subseteq\Omega(P)$  be  a  *-subalgebra  consisting  of
left-invariant elements.
Finally, let us assume that $\omega$ constructed in the above lemma
is a regular and
multiplicative connection. This
assumption implies certain specific  algebraic  relations  between
elements of $\Lf$.

It is  clear  that elements $\eta\in  \L$  are
horizontal. Hence the following relations hold
\begin{equation}\label{L-L}
\vartheta\eta+\sum_k \eta_k(\vartheta\circ a_k)=0
\end{equation}
where $\vartheta\in \L^\perp$ and $\Sum_k\eta_k\otimes a_k=
\chi(\eta)$.

The  action  of  the
covariant derivative on elements from $\cal B$
and $\L$ is described by

\begin{lem}
The following identities hold
\begin{gather}
\proj(b)=\k^\prime(b^{(1)})D_\omega(b^{(2)})
\label{homD2}\\
D_\omega(b)=b^{(1)}\proj(b^{(2)})\label{homD1}\\
D_\omega\proj(b)=
-R_\omega\pi j(b)-\proj(b^{(1)})\proj(b^{(2)}),\label{homDR}
\end{gather}
where $\proj=\Proj\pi'$.
\end{lem}
\begin{pf} Evidently, \eqref{homD2} and \eqref{homD1} are mutually
equivalent. Equation  \eqref{homD1}  directly  follows  from
\eqref{pi}, \eqref{Dhor}
and from the definition of $\omega$. Acting by $D_\omega$ on  both
sides of \eqref{homD2} and applying \eqref{DD2}--\eqref{Dprod}
we obtain
\begin{equation*}
\begin{split}
D_\omega\proj(b)&=D_\omega \k^\prime(b^{(1)})D_\omega
(b^{(2)})+\k^\prime(b^{(1)})D_\omega^2(b^{(2)})\\
&=\bigl\{D_\omega \k^\prime(b^{(1)})\bigr\}b^{(2)}\proj
(b^{(3)})
-\k^\prime(b^{(1)})b^{(2)}R_\omega\pi j(b^{(3)})\\
&=-\proj(b^{(1)})\proj(b^{(2)})-
R_\omega\pi j(b). \qed
\end{split}
\end{equation*}
\renewcommand{\qed}{}
\end{pf}

In fact, formula \eqref{homDR} {\it defines} $D_\omega$  and
$R_\omega.$ If $\pi^\prime(b)\in \L$ then
$$D_\omega\bigl(\pi^\prime(
b)\bigr)=-\proj(b^{(1)})\proj(b^{(2)})$$  and
similarly $$R_\omega\bigl(\pi^\prime(b)\bigr)=-
\proj(b^{(1)})\proj(b^{(2)}),$$ if
$\proj(b)=0.$  The  above  two formulas  are   in   fact
equivalent to  \eqref{homDR}.  Both  give  the  same  consistency
condition for $\Lf$.
If $b\in\cal{R}^\prime $ then
\begin{equation}
\proj(b^{(1)})\proj(b^{(2)})=0.\label{
cons}
\end{equation}
The  above  constraint  generates  a  further  constraint  in   the
third-order level, because it must be compatible with $D_\omega$
satisfying the graded Leibniz rule. Explicitly,
$$R_\omega\pi j(b^{(1)})\proj(b^{(2)})
-\proj(b^{(1)})R_\omega\pi j(b^{(2)})=0
$$
for each $b\in\cal R^\prime$. More generally, it follows that
\begin{equation}
R_\omega(\vartheta)\eta=\sum_k\eta_k R_\omega(\vartheta
\circ c_k)
\end{equation}
where $\vartheta\in\Gamma_{\inv}$ and $\eta\in \L$, while
$\Sum_k\eta_k\otimes c_k=\chi(\eta)$.

It is worth noticing that $\rho$ is extendible to a homomorphism
$\rho^\wedge\colon\Lf\rightarrow\Gamma_{\inv}^\wedge$ of
graded-differential algebras. In terms of the canonical identification
$\Omega(P)\leftrightarrow\cal{B}\otimes\Lf$ of spaces, the
verticalization homomorphism is given by
$\pi_v\leftrightarrow\id\otimes\rho^\wedge$.

Motivated by the  derived  expressions  and  constraints we shall
now construct ``the universal''
higher-order calculus on the
bundle, admitting regular and multiplicative connections
of the described geometrical nature. The starting  point  will  be  a
left-covariant
*-calculus  $\Psi$  over  $H$,   endowed   with   a
splitting of the form
$\Psi_{\inv}\cong \L\oplus\Gamma_{\inv}$.  We
shall   assume   that   this   splitting   possesses    all above
introduced properties.
Let $\cal{K}_1$ be the ideal in the tensor  algebra
$\L^\otimes$
generated by elements of the form
\begin{equation}
w=\proj(b^{(1)})\otimes\proj(b^{(2)})
\label{w-LL}
\end{equation}
where $b\in\cal R^\prime$. The formulas
\begin{align*}
D\bigl(\pi^\prime(b)\bigr)&=-
\proj(b^{(1)})\proj(b^{(2)})\\
R\bigl(\pi^\prime(q)\bigr)&=-\proj(q^{(1)})
\proj (q^{(2)})
\end{align*}
consistently define linear maps
$D\colon \L\rightarrow\L^\otimes/\cal{K}_1$ and $R\colon
\Gamma_{\inv}\rightarrow
\L^\otimes/\cal{K}_1$. Here, $\pi j(b)=0$ and
$\proj(q)=0$. We have
\begin{equation*}
D\proj(b)=
-R\pi j(b)-\proj(b^{(1)})
\proj(b^{(2)}).
\end{equation*}
Let $\cal{K}_2$ be the ideal in $\L^\otimes/\cal{K}_1$
generated by relations  of the form
$$ R(\vartheta)\eta=\sum_k\eta_k R(\vartheta\circ c_k)$$
where $\chi(\eta)=\Sum_k\eta_k\otimes c_k$.

The map $D$ can be uniquely extended to a
first-order derivation $D\colon \L^\sstar\rightarrow\L^\sstar$, where
$\L^\sstar=\bigl[(\L^\otimes/\cal{K}_1)/\cal{K}_2\bigr]$. Indeed, it
is sufficient to check that the graded Leibniz rule for $D$ is not in a
contradiction with relations generating $\cal{K}_1$ and $\cal{K}_2$.
This follows from the above derived equations.

Both ideals
$\cal{K}_1$ and $\cal{K}_2$ are right and *-invariant, in a natural manner.
In other
words $\L^\sstar$ is a *-algebra, endowed with the right action
$\chi\colon \L^\sstar\rightarrow\L^\sstar\otimes\cal A$.
Let us assume that $R$ is factorized through the ideal $\cal{K}_2$.
By construction, $D$ and $R$ are hermitian and right-covariant maps. In
particular, it follows that $DR=0$.

The $\circ$-structure on $\L^\otimes$ can be naturally ``projected''
to $\L^\sstar$, through ideals $\cal{K}_1$ and $\cal{K}_2$. The following
identities hold
\begin{gather}
D(\vartheta\circ b)=D(\vartheta)\circ b-
\proj(b^{(1)})(\vartheta\circ b^{(2)})+(-1)^{\partial
\vartheta}(\vartheta\circ b^{(1)})\proj(b^{(2)})\\
R\bigl(\vartheta\circ j(b)\bigr)=R(\vartheta)\circ b.
\end{gather}
Finally, let us consider a graded *-algebra defined as
$$\hor_P=\cal B\otimes\L^\sstar$$
at the level of graded vector spaces, while the
product  and  the *-structure are given by
\begin{gather*}
(q\otimes\eta)(b\otimes\vartheta)=qb^{(1)}\otimes(\eta\circ
b^{(2)})\vartheta\\
(b\otimes\vartheta)^*=b^{(1)*}\otimes (\vartheta^*\circ b^{(2)*}).
\end{gather*}
Evidently,  $\cal  B$  and  $\L^\sstar$  are  *-subalgebras   of
$\hor_P$. The formulas
\begin{gather*}
D(b\otimes\vartheta)=b^{(1)}\otimes\proj(b^{(2)})
\vartheta+b\otimes D(\vartheta)\\
F^{\sstar}(b\otimes\vartheta)=F(b)\chi(\vartheta)
\end{gather*}
define extensions $D\colon \hor_P\rightarrow \hor_P$ and
$F^{\sstar}\colon\hor_P\rightarrow \hor_P\otimes\cal A$ of the previously
introduced maps.  By construction, $F^{\sstar}$ defines the action of $G$
by ``automorphisms'' of $\hor_P$. Further,
\begin{lem}
The map $D$ is a hermitian right-covariant first-order
antiderivation
on $\hor_P$. The following identities hold
\begin{gather*}
D^2(\varphi)=-\sum_k\varphi_kR\pi(c_k)\\
R(\vartheta)\varphi=\sum_k\varphi_kR(\vartheta\circ c_k),
\end{gather*}
where
$F^\sstar(\varphi)=\Sum_k\varphi_k\otimes c_k$.
\end{lem}
\begin{pf} We compute
\begin{equation*}
\begin{split}
D(\vartheta b)&=D(b^{(1)})\vartheta\circ b^{(2)}+
b^{(1)}D(\vartheta\circ b^{(2)})\\
&=b^{(1)}\proj(b^{(2)})\vartheta\circ b^{(3)}+
b^{(1)}D(\vartheta)\circ b^{(2)}\\
&\phantom{=}-b^{(1)}\proj(b^{(2)})\vartheta\circ
b^{(3)}+(-1)^{\partial\vartheta}b^{(1)}(\vartheta\circ b^{(2)})
\proj(b^{(3)})\\
&=D(\vartheta)b+(-1)^{\partial\vartheta}\vartheta D(b).
\end{split}
\end{equation*}
This implies that $D$ is  a  (first-order)  antiderivation
on $\hor_P$. Furthermore, it is sufficient to check that the relation
between $D^2$ and $R$ holds on elements from $\L$ and $\cal{B}$. We have
\begin{multline*}
D^2(b)=D\bigl(b^{(1)}\proj(b^{(2)})\bigr)=
b^{(1)}\proj (b^{(2)})\proj(b^{(3)})-b^{(1)}R\pi j(b^{(2)})-
b^{(1)}\proj (b^{(2)})\proj(b^{(3)})\\ =-b^{(1)}R\pi j(b^{(2)})
\end{multline*}
and similarly
\begin{equation*}
\begin{split}
D^2\proj(b)&=-D\bigl(R\pi j(b)+\proj(b^{(1)})\proj(b^{(2)})\bigr)\\
&=R\pi j(b^{(1)})\proj(b^{(2)})-\proj(b^{(1)})R\pi j(b^{(2)})\\
&=\proj(b^{(3)})R\bigl[\pi j(b^{(1)})\circ
j\bigl(\k(b^{(2)})b^{(4)}\bigr)\bigr]-\proj(b^{(1)})R\pi j(b^{(2)})\\
&=-\proj(b^{(2)}) R\pi j\bigl(\k(b^{(1)})b^{(3)}\bigr).
\end{split}
\end{equation*}
Finally, we have to check the commutation relation between $R$ and
$\cal{B}$. Direct transformations give
\begin{multline*}
R\pi j(q)b=-\proj(q^{(1)})\proj(q^{(2)})b=-b^{(1)}\bigl(
\proj(q^{(1)})\circ b^{(2)}\bigr)\bigl(
\proj(q^{(2)})\circ b^{(3)}\bigr)\\
=-b^{(1)}\proj(q^{(1)}b^{(2)})\proj(q^{(2)}b^{(3)})=-b^{(1)}R\pi
j(qb^{(2)}),
\end{multline*}
where $\proj(q)=0$. This completes the proof.
\end{pf}

Now the construction of the  full  differential  calculus  on  the
bundle $P$ can be completed applying ideas of Subsection 6.
The initial splitting is naturally understandable as a regular and
multiplicative connection $\omega$ on $P$. Maps $D$  and  $R$  are
interpretable  as  the  corresponding   operators   of   covariant
derivative  and  curvature.  The  full  algebra   $\Omega(P)$   is
left-covariant   (over   $H$).   The associated first-order calculus
coincides with $\Psi$. The   corresponding
differential  *-subalgebra  $\Lf$  of  left-invariant
elements can be independently described as follows. At  the  level
of (graded) vector spaces
$$\Lf=\L^\sstar\otimes\Gamma_{\inv}^\wedge.$$
The differential *-algebra structure is specified by
\begin{gather*}
(\eta\otimes\xi)(\vartheta\otimes\zeta)=(-1)^{\partial\xi\partial\vartheta}
\sum_k\eta\vartheta_k\otimes(\xi\circ c_k)\zeta\\
(\vartheta\otimes\zeta)^*=\sum_k\vartheta_k^*\otimes(\zeta^*\circ
c_k^*)\\
d^\wedge(\vartheta\otimes\zeta)=D(\vartheta)\otimes\zeta
+(-1)^{\partial\vartheta}\sum_k
\vartheta_k\otimes\pi(c_k)\zeta+(-1)^{\partial\vartheta}\vartheta
d^\wedge(\zeta).
\end{gather*}
Here $d^\wedge\colon \Lf\rightarrow\Lf$ is the
corresponding differential and $\Sum_k\vartheta_k\otimes
c_k=\chi(\vartheta)$. Finally,
$$d^\wedge(\zeta)=d(\zeta)+R(\zeta)$$
for $\zeta\in\Gamma_{\inv}$. The map $\rho^\wedge$ is given by
$\rho^\wedge\leftrightarrow\epsilon_{\L}\otimes\id$, where
$\epsilon_{\L}$ is a character on $\L^\sstar$ specified by
$\epsilon_{\L}(\L)=\{0\}$.

As a concrete example, let us consider the quantum  Hopf  fibering
\cite{P1}. This bundle is decribed by the quantum group
$H=S_\mu U(2)$, and its subgroup $G=H_{cl}=U(1)$. The base  manifold
is the quantum 2-sphere \cite{P1}.

By definition \cite{W-su2}, $\cal B$  is  the
*-algebra generated
by elements $\alpha$ and $\gamma$ and the following relations
\begin{gather*}
\alpha\alpha^*+\mu^2\gamma\gamma^*=1\qquad\alpha^*\alpha+\gamma^*\gamma=1\\
\alpha\gamma=\mu\gamma\alpha\qquad
\alpha\gamma^*=\mu\gamma^*\alpha\qquad
\gamma\gamma^*=\gamma^*\gamma.
\end{gather*}
The fundamental representation is given by
\begin{equation*}
u=(u^\dagger)^{-1}=\begin{pmatrix}
\alpha&-\mu\gamma^*\\
\gamma&\phantom{-\mu}\alpha^*
\end{pmatrix}
\end{equation*}
where $\mu\in(-1,1)\setminus\{0\}$.

Let us first assume  that the first-order differential structure
$\Psi$ over
$H$ coincides with the $3D$-calculus \cite{W-su2},
based on a  right  ${\cal  B}$-ideal
$$
\cal{R}'=\gen\Bigl\{\,\gamma^2,\,
\gamma\gamma^*,\,\gamma^{*2},\,
\alpha\gamma-\gamma,\,\alpha\gamma^*-\gamma^*,\,
\mu^2\alpha+\alpha^*-(1+\mu^2)1\,\Bigr\}.
$$
The  space  $\Psi_{\inv}$  is
$3$-dimensional and spanned by elements
$$\eta=\pi^\prime(\alpha-\alpha^*)\qquad
\eta_+=\pi^\prime(\gamma)\qquad
\eta_-=\pi^\prime(\gamma^*).$$
The corresponding right ${\cal B}$-module structure $\circ$
is specified by
\begin{equation*}
\begin{aligned}
\mu^2\eta\circ\alpha&=\eta\\
\mu\eta_\pm\circ\alpha&=\eta_\pm
\end{aligned}\qquad
\begin{aligned}
\eta\circ\alpha^*&=\mu^2\eta\\
\eta_\pm\circ\alpha^*&=\mu\eta_\pm
\end{aligned}
\end{equation*}
with $\Psi_{\inv}\circ\gamma=\Psi_{\inv}\circ\gamma^*=\{0\}$.

The *-algebra $\cal A$ of polynomial functions on $G$ is generated
by the  canonical  unitary  element  $z=j(\alpha)$  (and  we  have
$j(\gamma)=j(\gamma^*)=0$). Let us assume that  $\Gamma$  is  a
left-covariant calculus over $G$ based on the right
$\cal   A$-ideal   $\cal   R=j({\cal   R}^\prime)$.   The    space
$\Gamma_{\inv}$    is     1-dimensional,     and     spanned     by
$\zeta=\rho(\eta)=\pi(z-z^*)$. We have
$\rho(\eta_+)=\rho(\eta_-)=0$  and  $\Gamma^{\wedge  k}=\{0\}$
for $k\geq 2$. However, the calculus based on $\Gamma^\wedge$
differs from the classical  differential  calculus,  because
the right $\cal A$-module structure on $\Gamma_{\inv}$ is given by
\begin{equation*}
\mu^2\zeta\circ z=\zeta\qquad
\zeta\circ z^*=\mu^2\zeta.
\end{equation*}

The space $\L$ is spanned by $\eta_+$ and $\eta_-$. Let us define a
splitting
$\Psi_{\inv}\cong \L\oplus\Gamma_{\inv}$ ($\Leftrightarrow$
 the space $\L^\perp$) by   identifying
$\eta$   and   $\zeta$.  The corresponding relations   determining
the algebra $\L^\sstar$ are
\begin{gather*}
\eta_+^2=\eta_-^2=0\\
\eta_+\eta_-=-\mu^2\eta_-\eta_+
\end{gather*}
while the third-order relations  are trivialized.
Additional    relations    determining    the algebra
$\Lf$ are
\begin{equation*}
\eta\eta_+=-\frac{1}{\mu^4}\eta_+\eta\qquad
\eta^2=0\qquad
\eta\eta_-=-\mu^4\eta_-\eta.
\end{equation*}
It is worth noticing that  $\Omega(P)=\Psi^\wedge$  and  hence
the  higher-order  calculus  on  the  bundle  coincides  with  the
calculus constructed in \cite{W-su2}.
Modulo differences between general formulations,
the canonical regular connection $\omega$ (associated to the fixed
splitting  of  $\Psi_{\inv}$)  coincides  with   a   connection
constructed in \cite{BM}.
The  curvature  of $\omega$ is
given by
\begin{equation}
R_\omega(\zeta)=d\eta=\mu(1+\mu^2)\eta_-\eta_+.
\end{equation}

Next, let us consider the case of the $4D_+$-calculus over  $S_\mu
U(2).$   This   calculus   $\Psi$   is   bicovariant   and
*-covariant.  By  definition \cite{W3} the
corresponding   ${\cal
R}^\prime$ is generated (as a right ideal) by multiplets
\begin{gather*}
\mbox{\it 1}=\Bigl\{\,a\bigl(\mu^2\alpha+\alpha^*-(1+\mu^2)1\bigr)\,\Bigr\}
\qquad\mbox{\it
3}=\Bigl\{\,a\gamma,\,a(\alpha-\alpha^*),\,a\gamma^*\,\Bigr\}\\
\mbox{\it 5}=\Bigl\{\,
\gamma^2,\,\gamma(\alpha-\alpha^*),\,\mu^2\alpha^{*2}-
(1+\mu^2)(\alpha\alpha^*-\gamma\gamma^*)+\alpha^2,\,
\gamma^*(\alpha-\alpha^*),\,\gamma^{*2}\,\Bigr\}
\end{gather*}
where $a=\mu^2\alpha+\alpha^*-(\mu^3+1/\mu)1$. The  space
$\Psi_{\inv}$ is 4-dimensional. A natural  basis  is  given  by
elements
\begin{gather*}
\tau=\pi^\prime(\mu^2\alpha+\alpha^*)\\
\eta_+=\pi^\prime(\gamma)\qquad\eta=\pi^\prime(\alpha-\alpha^*)
\qquad\eta_-=\pi^\prime(\gamma^*).
\end{gather*}
Elements $\eta_+$, $\eta_-$ and $\eta$ form a triplet (relative
to the adjoint action). The element $\tau$ is
$\adj '$-invariant. We have (\cite{D}--Section 6)
\begin{gather*}
\eta_+\circ\gamma^*=\eta_-\circ\gamma=-\frac{(1+\mu)(1-\mu^2)}{\mu(1+
\mu^2)(1-\mu^3)}\tau-\frac{1-\mu^2}{\mu(1+\mu^2)}\eta\\
\eta_+\circ\gamma=\eta_-\circ\gamma^*=0\\
\begin{aligned}
\eta_+\circ\alpha&=\eta_+=\eta_+\circ\alpha^*\\
\eta\circ\gamma&=-\frac{1-\mu^2}{\mu}\eta_+
\end{aligned}\qquad\qquad
\begin{aligned}
\eta_-\circ\alpha&=\eta_-=\eta_-\circ\alpha^*\\
\eta\circ\gamma^*&=-\frac{1-\mu^2}{\mu}\eta_-
\end{aligned}\\
\begin{aligned}
-\eta\circ\alpha^*&=\frac{(1+\mu)(1-\mu^2)}{\mu(1+\mu^2)
(1-\mu^3)}\tau-\frac{2\mu}{1+\mu^2}\eta\\
\eta\circ\alpha&=\frac{\mu(1+\mu)(1-\mu^2)}{(1+\mu^2)
(1-\mu^3)}\tau+\frac{2\mu}{1+\mu^2}\eta
\end{aligned}\\
\begin{aligned}
\tau\circ\gamma&=\frac{(1-\mu)(1-\mu^3)}{\mu}\eta_+\\
\tau\circ\gamma^*&=\frac{(1-\mu)(1-\mu^3)}{\mu}\eta_-
\end{aligned}\qquad
\begin{aligned}
\tau\circ\alpha^*&=\frac{1+\mu^4}{\mu(1+\mu^2)}\tau-
\frac{\mu(1-\mu)(1-\mu^3)}{1+\mu^2}\eta\\
\tau\circ\alpha&=\frac{1+\mu^4}{\mu(1+\mu^2)}\tau+
\frac{(1-\mu)(1-\mu^3)}{\mu(1+\mu^2)}\eta.
\end{aligned}
\end{gather*}

It turns  out  that  the  ideal  $\cal  R=j({\cal  R}^\prime)$  is
generated by the element $z+\mu z^*-(1+\mu)1$. The  space  $\L$  is
spanned by elements $\eta_+$, $\eta_-$ and
$\xi=\tau+\bigl((1-\mu^3)/(1+\mu)\bigr)\eta.$ It is  worth  noticing  that
$\xi$ is (the unique) common characteristic vector for all
operators of the form $\circ
b$ (where $b\in\cal B$). Explicitly,
$$\xi\circ\alpha=\frac{1}{\mu}\xi\qquad\xi\circ\alpha^*=\mu\xi
\qquad\xi\circ\gamma=\xi\circ\gamma^*=0.$$

However, the  space  $\L$  does  not  possess  a  $\circ$-invariant
complement. A natural choice of the  complement  $\L^\perp$  is  to
consider the subspace spanned by $\eta$ (following
a weak analogy with the previous example). The calculus on $G=U(1)$
is   non-standard.   The   higher-order   calculus   is    trivial.
The corresponding right $\cal A$-module structure is given by
\begin{equation*}
\zeta\circ z=\mu\zeta\qquad
\zeta\circ z^*=\frac{1}{\mu}\zeta
\end{equation*}
where $\zeta=\rho(\eta).$
Let us  assume  that  the
higher-order calculus on the bundle is based  on  the bicovariant
exterior  algebra  $\Psi^\vee$  \cite{W3}.
Let  $\omega$  be   the
connection corresponding to the above splitting.  This  connection
is not multiplicative.

The most general form of a  connection (coming  from a
complement $\L^\perp$) is
\begin{equation}
\omega(\zeta)=\eta+t\xi
\end{equation}
where $t\in \Re.$
All these connections are non-regular. A connection  $\omega$
is multiplicative iff $(\eta+t\xi)^2=0$, which  is  equivalent  to
$t=-(1+\mu)/(1-\mu^3)$,   in   other   words   the   corresponding
complement is spanned by $\tau$.
However
it  is worth noticing that
if the higher-order calculus  is  described  by  $\Psi^\wedge$
then $\omega$ is not multiplicative, because
$\tau^2\neq 0$ in this case.

The curvature is given by
\begin{equation}
R_{\omega}=d(\eta+t\xi)=\biggl(
\frac{\mu t}{1-\mu^2}+\frac{\mu}{(1-\mu)(1-\mu^3)}\biggr)
(\tau\eta+\eta\tau).
\end{equation}
We see that for the above mentioned special value of
$t$ the curvature vanishes.

Finally, let us assume that  the  calculus  on  $G$  is {\it
classical }
(based  on  standard  differential  forms).  Let  us  assume  that
$\Psi$ is a bicovariant *-calculus. This implies
$(X\otimes \id) \ad^\prime({\cal R}^\prime)=\{0\}$, where $X$  is  the
canonical generator of $\lie(G)$. In other words, $\Psi$  is
{\it admissible} (in the sense of the previous paper).
Let us assume  that  $\Psi$  is  the  {\it minimal}
admissible
(bicovariant  *-)  calculus.  The  space  $\Psi_{\inv}$  can   be
naturally  identified  with  the  algebra  $\cal  D$ consisting of
all elements $b\in\cal B$ invariant under the left action  of  $G$
on $H$ (polynomial functions on a
quantum 2-sphere). In terms of this identification
\begin{equation*}
\adj '\,\leftrightarrow\,(\phi^\prime\restr\cal D)
\colon \cal D\rightarrow\cal D\otimes\cal B\qquad
(\,)\circ b \,\leftrightarrow\, \k(b^{(1)})(\,)b^{(2)} .
\end{equation*}
Let  $\tau$  be  the  element  corresponding to $1\in\cal D$
and let us assume that $\L^\perp$ is spanned by $\tau$.
Explicitly,
\begin{equation}
\tau=\frac{2}{\mu^2-1}\pi^\prime(\mu^2\alpha+\alpha^*).
\end{equation}
The element $\tau$ is right-invariant ($\L^\perp$ consists precisely
of right-invariant elements of $\Psi_{\inv}$
and each integer-valued
irreducible  multiplet   appears   without   degeneracy   in   the
decomposition   of   $\adj '$   into    irreducible
components), and
\begin{equation}\label{tau-mod}
\tau\circ b=\e^\prime(b)\tau.
\end{equation}
In particular
\begin{equation}\label{tau-sgm}
\begin{split}
\sigma(\tau\otimes\vartheta)&=\vartheta\otimes\tau\\
\sigma(\vartheta\otimes\tau)&=\tau\otimes\vartheta.
\end{split}
\end{equation}
Let us assume that the higher-order  calculus  on  the  bundle  is
described by the exterior algebra $\Psi^\vee$. We have
\begin{gather*}
\tau^2=0\\
d\tau=0
\end{gather*}
Identities
\eqref{tau-mod} and \eqref{tau-sgm} imply
\begin{equation*}
\tau w=(-1)^{\partial w}w\tau
\end{equation*}
for   each   $w\in   \Psi^\vee$.   The   connection   $\omega$
corresponding  to  $\L^\perp$  is   regular   and   multiplicative.
Moreover, $\omega$ is {\it flat} in the sense that $R_\omega=0$.

\subsection{A Constructive Approach to Differential Calculus}
Every  regular  connection  $\omega$  induces  the
isomorphism
$m_{\omega}\colon \vh(P)\leftrightarrow
\Omega(P)$
(if $\J(P)=\{0\}$). Moreover,  if
two algebras are identified  with  the  help of $m_{\omega}$
then  it  is
possible to express the differential  structure  on
$\Omega(P)$  in
terms of the algebra structure on
$\vh(P)$, and  the following maps:
\begin{gather*}
D_{\omega}\colon \hor(P)\rightarrow \hor(P)\qquad
R_\omega\colon\Gamma_{\inv}\rightarrow\hor(P)\\
d\colon \Gamma_{\inv}^{\wedge}\rightarrow\Gamma_{\inv}^{\wedge}\qquad
F^{\wedge}\colon \hor(P)\rightarrow \hor(P)\otimes\cal{A}.
\end{gather*}
Explicitly,
\begin{equation}
d^\wedge(\varphi\otimes\vartheta)=D_{\omega}(\varphi)\otimes\vartheta
+(-1)^{\partial\varphi}\sum_k
\varphi_k \otimes\pi(c_k)\vartheta+(-1)^{\partial\varphi}
\varphi d^\wedge(\vartheta),
\end{equation}
where $d^\wedge$ is the corresponding differential
on $\vh(P)$ and $d^\wedge\restr\Gamma^\wedge_{\inv}$ is fixed by
$$ d^\wedge(\vartheta)=d(\vartheta)+R_\omega(\vartheta),$$
for $\vartheta\in\Gamma_{\inv}$, and extended on
$\Gamma_{\inv}^\wedge$ by the graded Leibniz rule.

     It is worth noticing that  the  curvature $R_{\omega}$
is  completely determined by the covariant derivative,  as  easily
follows from \eqref{DD2} and property ({\it qpb4\/}) from Section~3. Indeed,
we have
\begin{equation}
R_\omega\pi(a)=-\sum_k q_kD^2_\omega(b_k)
\end{equation}
where $q_k,b_k\in\cal B$ are such that \eqref{free} holds.

 In this subsection an ``opposite''  construction will be presented,
which builds the algebra $\Omega(P)$ and a regular multiplicative
connection $\omega $
starting from a *-algebra playing the role  of  horizontal  forms,
and three operators imitating the covariant derivative, curvature
and the right action.

     Let $P=(\cal{B},i,F)$ be a quantum principal $G$-bundle  over
$M$. Let
$$\hor_P=\sideset{}{^\oplus}\sum_{k\geq 0} \hor_P^k$$
be a graded *-algebra such  that $\hor^0_P=\cal{B}$.
Further,  let  us assume that a grade-preserving
*-homomorphism
$F^\sstar\colon \hor_P\rightarrow \hor_P\otimes\cal{A}$
extending      the      map
$F$ is given such that
\begin{align}
(F^\sstar\otimes \id) F^\sstar&=(\id \otimes\phi)F^\sstar\\
(\id \otimes \e)F^\sstar&=\id .
\end{align}
Let us assume that a linear first-order
map $D\colon \hor_P\rightarrow \hor_P$
is given such that the following properties hold
\begin{gather}
D(\varphi\psi)=D(\varphi)\psi+(-1)^{\partial\varphi}\varphi D
(\psi)\label{E4}\\
D*=*D\\
F^{\sstar}D=(D\otimes
\id) F^{\sstar}.\label{r-covD}
\end{gather}

Finally, let us assume that there exists a linear map
$R\colon \Gamma_{\inv}\rightarrow \hor_P$ such that
\begin{equation}\label{E7}
D^2(\varphi)=-\sum_k\varphi_k R\pi(c_k),
\end{equation}
for each $\varphi\in \hor_P$, where $F^{\sstar}(\varphi)=\Sum_k
\varphi_k\otimes c_k$.

 Evidently, $D$ plays the role of the  covariant
derivative.
The  map  $R$
is  determined uniquely by the above
condition.
It plays the role of the curvature map. Explicitly
\begin{equation}\label{E12}
R\pi(a)=-\sum_kq_kD^2(b_k)
\end{equation}
where $q_k,b_k\in{\cal B}$ are such that \eqref{free} holds.

We pass to a
``reconstruction'' of the calculus $\Omega(P)$.
Let us  first  analyze the map $R$ in more details.
\begin{pro} The following identities hold
\begin{gather}
F^{\sstar}R=(R\otimes \id) \adj \label{E9}\\
DR=0\label{E8}\\
R*=*R\label{E10}\\
R(\vartheta)\varphi=\sum_k\varphi_k R(\vartheta\circ c_k).
\label{E11}
\end{gather}
\end{pro}
\begin{pf}
Acting by $D\otimes \id $ on equality \eqref{free}
and using \eqref{E4} we obtain
$$0=\sum_kD(q_k)F(b_k)+\sum_kq_k(D\otimes \id) F(b_k). $$
This, together with \eqref{r-covD} and \eqref{E7}, implies
$$0=\sum_kD(q_k)D^2(b_k)+\sum_kq_kD^3(b_k)=D\Bigl(\sum_kq_kD^2(b_k)\Bigr)
=-DR\pi(a). $$
Hence, \eqref{E8} holds. Further, \eqref{free} implies
$$ 1\otimes1\otimes a=\sum_{kj}F(q_k)F(b_{kj})\otimes a_{kj}$$
and hence
\begin{equation}
1\otimes \k(a^{(1)})\otimes a^{(2)}=\sum_{ijk}q_{ki}b_{kj}
\otimes z_{ki}\otimes a_{kj},\label{E13}
\end{equation}
where  $F(b_k)=\Sum_j  b_{kj}\otimes  a_{kj}$  and $F(q_k)=\Sum_i
q_{ki}\otimes z_{ki}$. From \eqref{E13} we find
\begin{gather}
1\otimes \k(a)=\sum_kF(q_k)b_k\label{E14}\\
\e(a)1=\sum_k q_kb_k\label{E15}\\
1\otimes a^{(2)}\otimes \k(a^{(1)})a^{(3)}=\sum_{ijk}q_{ki}F(
b_{kj})\otimes z_{ki}a_{kj}.\label{E16}
\end{gather}

Identity \eqref{E15} also directly follows
from \eqref{free}.  A  direct
computation now gives
$$
R\bigl[\pi(a)^*\bigr]=-R\pi\bigl(\k(a)^*\bigr)=\sum_kb_k^*D^2(q_k^*)=
-\sum_k D^2(b_k^* )q_k^* =\bigl[R\pi(a)\bigr]^*,
$$
and similarly
\begin{multline*}
(R\otimes \id) \adj \pi(a)=R\pi(a^{(2)})\otimes
\k(a^{(1)})a^{(2)}=-
\sum_{ijk}q_{ki}D^2(b_{kj})
\otimes z_{ki}a_{kj}\\
=-F^{\sstar}\biggl(\sum_kq_kD^2(b_k)\biggr)=F^{\sstar}R\pi(a).
\end{multline*}
This proves \eqref{E9} and \eqref{E10}.
Finally, let us prove \eqref{E11}. We have
$$\sum_{nkj} q_kb_{kj}\varphi_n \otimes a_{kj}c_n =
\sum_n \varphi_n \otimes ac_n, $$
for each $\varphi\in \hor_P$. Hence
$$\sum_{nkj} q_kb_{kj}\varphi_n \otimes \pi(a_{kj}c_n) =
\sum_n \varphi_n \otimes \pi(a)\circ c_n, $$
if $a\in \ker(\e)$. The above equality, together with \eqref{E7}
implies
$$
\sum_{nkj} q_kb_{kj}\varphi_nR\pi(a_{kj}c_n) =
-\sum_kq_kD^2(b_k\varphi)=\sum_n \varphi_n R\bigl[\pi(a)\circ c_n\bigr].$$
On the other hand, \eqref{E12} and \eqref{E15} imply
$$
\bigl[R\pi(a)\bigr]\varphi=-\sum_kq_kD^2(b_k)\varphi=
-\sum_kq_kD^2(b_k\varphi). $$
Hence, property \eqref{E11} holds.
\end{pf}

 Let us consider the graded space
$\Omega(P)=\hor_P\otimes\Gamma_{\inv}^{\wedge}$, endowed
with the following *-algebra structure
\begin{align}
(\psi\otimes\eta)(\varphi\otimes\vartheta)&=(-1)^{\partial
\varphi\partial\eta}
\sum_k \psi\varphi_k \otimes(\eta\circ c_k)\vartheta\\
(\varphi\otimes\vartheta)^*&=\sum_k\varphi_k^*
\otimes(\vartheta^* \circ c^*_k),
\end{align}
where $F^{\sstar}(\varphi)=\Sum_k\varphi_k\otimes c_k$.
Algebras $\hor_P$ and $\Gamma^{\wedge}_{\inv}$ are understandable as
*-sub\-algeb\-ras of $\Omega(P)$, in a natural manner.
\begin{lem}\label{lem:R+d}
There   exists   the   unique antiderivation
$d^{\wedge}\colon \Gamma^{\wedge}_{\inv}\rightarrow \Omega(P)$ satisfying
\begin{equation}
d^{\wedge}(\vartheta)=R(\vartheta)+d(\vartheta)\label{E19}
\end{equation}
for each $\vartheta\in\Gamma_{\inv}$.
\end{lem}
\begin{pf}
The graded Leibniz rule implies
that the values of $d^{\wedge}$
on higher-order forms are completely determined by the
restriction $d^{\wedge}\restr\Gamma_{\inv}$ (and we have $d^{\wedge}1=0$.)
Hence, $d^{\wedge}$ is unique, if exists.

 Let  us  prove  that $d^{\wedge}$
can  be  consistently  constructed   by
extending, with the help of the graded Leibniz rule,
a linear map (acting on  $\Gamma_{\inv}$)
given by \eqref{E19}.
The extension exists iff  contradictions  do  not
appear at the level of second-order constraints defining the algebra
$\Gamma_{\inv}^{\wedge}$.
Simple transformations give
\begin{equation*}
\begin{split}
\left\{\pi(a^{(1)})\pi(a^{(2)})\right\}&\rightarrow
R\pi(a^{(1)})\pi(a^{(2)})\\
&\phantom{\rightarrow}-\pi(a^{(1)})\pi(a^{(2)})\pi(a^{(3)})\\
&\phantom{\rightarrow} +
\pi(a^{(1)})\pi(a^{(2)})\pi(a^{(3)})
-\pi(a^{(1)})R\pi(a^{(2)})\\
&=R\pi(a^{(1)})\pi(a^{(2)})\\
&\phantom{\rightarrow}-R\pi(a^{(3)})\bigl[\pi(a^{(1)})\circ
\bigl(\k(a^{(2)})a^{(4)}\bigr)\bigr]\\
&=R\pi(a^{(2)})\pi\bigl(\k(a^{(1)})a^{(3)}\bigr)=0,
\end{split}
\end{equation*}
for each $a\in\cal{R}$. Hence, $d^{\wedge}$ exists.
\end{pf}

The formula
\begin{equation}
d^{\wedge}
(\varphi\otimes\vartheta)=D(\varphi)\otimes\vartheta+
(-1)^{\partial\varphi}\varphi d^{\wedge}(\vartheta)+
(-1)^{\partial\varphi}\sum_k\varphi_k\otimes\pi(c_k)
\vartheta
\end{equation}
defines a linear first-order map
$d^{\wedge}\colon \Omega(P)\rightarrow\Omega(P)$ extending $d^{\wedge}$
introduced in the previous lemma.
\begin{pro}
The following identities hold
\begin{gather}
d^{\wedge} *=*d^{\wedge}\\
(d^{\wedge})^2 =0\\
d^{\wedge}(wu)=
d^{\wedge}
(w)u+(-1)^{\partial w}wd^{\wedge}(u).\label{E25}
\end{gather}
\end{pro}
\begin{pf}
A  direct  calculation   gives for
$w\in\Gamma_{\inv}$ and $u=\varphi\otimes\vartheta$,
\begin{equation*}
\begin{split}
 d^{\wedge}(wu)&=(-1)^{\partial\varphi}d^{\wedge}\biggl(\sum_k
\varphi_k\otimes(w\circ c_k)\vartheta\biggr)
=\sum_k\varphi_k\otimes d(w\circ c_k)\vartheta\\
&\phantom{=}+\sum_k\varphi_kR(w\circ c_k)\otimes\vartheta
-\biggl(\sum_k\varphi_k\otimes
(w\circ c_k)\biggr)d^\wedge(\vartheta)\\
&\phantom{=}+(-1)^{\partial\varphi}\sum_k\bigl(
D(\varphi_k)\otimes(w
\circ c_k)\vartheta\bigr) +\sum_k\bigl(
\varphi_k\otimes\pi(c_k^{(1)})(w\circ c_k^{(2)}
)\vartheta\bigr)\\
&=-w\bigl[D(\varphi)\otimes\vartheta\bigr]-(-1)^{\partial\varphi}w
\varphi d^{\wedge}(\vartheta)+d^{\wedge}(w)u\\
&\phantom{=}-(-1)^{\partial\varphi}\sum_kw\bigl(\varphi_k\otimes
\pi(c_k)\vartheta\bigr)=d^{\wedge}(w)u-wd^{\wedge}(u).
\end{split}
\end{equation*}
A similar computation shows that
\eqref{E25} holds   for $w\in \hor_P$,
\begin{equation*}
\begin{split}
d^\wedge(wu)&=(-1)^{\partial\varphi+\partial w}\sum_{lk}w_l\varphi_k
\otimes\pi(d_lc_k)\vartheta
+(-1)^{\partial\varphi+\partial w}w\varphi d^\wedge(\vartheta)
+D(w\varphi)\otimes\vartheta\\
&=D(w)\varphi\otimes\vartheta+(-1)^{\partial w}wD(\varphi)\otimes\vartheta+
(-1)^{\partial w+\partial\varphi}w\varphi d^\wedge(\vartheta)\\
&\phantom{=}+(-1)^{\partial\varphi +\partial w}\sum_k w\varphi_k\otimes
\pi(c_k)\vartheta+(-1)^{\partial\varphi +\partial w}\sum_{kl}
w_l\varphi_k\otimes\bigl(\pi(d_l)\circ c_k\bigr)\vartheta\\
&=d^\wedge(w)u+(-1)^{\partial w}wd^\wedge(u),
\end{split}
\end{equation*}
where $\Sum_l w_l\otimes d_l=F^{\sstar}(w)$. It follows that \eqref{E25}
holds  for arbitrary $w,u\in\Omega(P)$.

 Let us prove that the square of $d^{\wedge}$ vanishes. We have
\begin{equation*}
\begin{split}
(d^{\wedge})^2\pi(a)&=d^{\wedge}\bigl(-\pi(a^{(1)})\pi(a^{(2)}
)+R\pi(a)\bigr)=DR\pi(a)\\
&\phantom{=}+\pi(a^{(1)})\pi(a^{(2)})\pi(a^{(3)})-
\pi(a^{(1)})\pi(a^{(2)})\pi(a^{(3)})\\
&\phantom{=}-\bigl[R\pi(a^{(1)})\bigr]\pi(a^{(2)})+
\pi(a^{(1)})R\pi(a^{(2)})\\
&\phantom{=}+R\pi(a^{(2)})\pi\bigl(\k(a^{(1)})a^{(3)}\bigr)\\
&=R\pi(a^{(3)})\bigl[\pi(a^{(1)})\circ \bigl(\k(a^{(2)})a^{(4)}\bigr)\bigr]-
\bigl(R\pi(a^{(1)}
)\bigr)\pi(a^{(2)})\\
&\phantom{=}+R\pi(a^{(2)})\pi\bigl(\k(a^{(1)})a^{(3)}\bigr)=0
\end{split}
\end{equation*}
for each
$a\in\cal{A}$. Further,
\begin{multline*}
 (d^{\wedge})^2\varphi=d^{\wedge}\biggl(D\varphi+(-1)^{\partial
\varphi}\sum_k\varphi_k \pi(c_k
)\biggr)=D^2\varphi-(-1)^{\partial\varphi}
\sum_kD(\varphi_k)\pi(c_k)\\
\mbox{}+(-1)^{\partial\varphi}\sum_kD(\varphi_k)\pi(c_k
)+\sum_k\Bigl\{\varphi_k\pi(c_k^{(1)})\pi(c_k^{(2)})+
\varphi_kR\pi(c_k)+\varphi_kd\pi(c_k)\Bigr\}=0,
\end{multline*}
for each $\varphi\in \hor_P$. Having in mind that  spaces
$\hor_P$ and $\Gamma_{\inv}$
generate $\Omega(P)$, and using the fact that the square of an
antiderivation is a derivation we conclude that
$(d^{\wedge})^2=0$.

     Finally,
\begin{equation*}
\begin{split}
d^{\wedge}(\varphi^*)&=D(\varphi^*)+(-1)^{\partial\varphi}
\sum_k\varphi_k^*\pi(c_k^*)\\
&=(D\varphi)^* +(-1)^{\partial\varphi}\sum_k\varphi_k^*
\pi(c_k^{(2)})^* \circ c_k^{(1)*}\\
&=\bigl(D\varphi+(-1)^{\partial
\varphi}\sum_k\varphi_k\pi(c_k)\bigr)^*=d^{\wedge}(\varphi)^*
\end{split}
\end{equation*}
for each $\varphi\in \hor_P$. It follows that $d^{\wedge}$
is a hermitian map.
\end{pf}
 We are  going  to  prove  that $\Omega(P)$   satisfies  condition
({\it diff2\/}). The formula
\begin{equation}
\widehat{F}(\varphi\otimes\vartheta)=F^{\sstar}
(\varphi)\widehat{\adj} (\vartheta)
\end{equation}
determines a linear map
$\widehat{F}\colon \Omega(P)\rightarrow\Omega(P)
\grten\Gamma^{\wedge}$
extending $F^{\sstar}$  and $\widehat{\adj} $.
\begin{pro}
The  map  $\widehat{F}$ is  a  homomorphism  of differential *-algebras.
\end{pro}
\begin{pf}
For each $\varphi\in \hor_P$ and
$\vartheta\in\Gamma_{\inv}$ we have
\begin{equation*}
\begin{split}
\widehat{F}(\vartheta\varphi)&=(-1)^{\partial\varphi}
\sum_k\widehat{F}\bigl(\varphi_k
(\vartheta\circ c_k
)\bigr)=(-1)^{\partial\varphi}
\sum_k\bigl(\varphi_k\otimes c_k^{(1)}(\vartheta\circ
c_k^{(2)})\bigr)\\
&\phantom{=}+ (-1)^{\partial\varphi}
\sum_{kl}(\varphi_k\otimes c_k^{(1)})
(\vartheta_l\circ c_k^{(3)}
\otimes \k(c_k^{(2)})a_lc_k^{(4)})\\
&=(-1)^{\partial\varphi}\sum_k\varphi_k \otimes
(\vartheta c_k)+ (-1)^{\partial\varphi}
\sum_{kl}\varphi_k(\vartheta_l\circ c_k^{(1)})\otimes
a_lc_k^{(2)}\\
&=\bigl(1\otimes\vartheta+\adj (\vartheta)\bigr)
\sum_k\varphi_k \otimes c_k =\widehat{F}(\vartheta)
\widehat{F}(\varphi).
\end{split}
\end{equation*}
Here, $\Sum_k\varphi_k\otimes c_k=F^{\sstar}(\varphi)$        and
$\adj (\vartheta)=\Sum_l\vartheta_l\otimes a_l$.
Using  the   facts   that
$F^{\sstar}$  and $\widehat{\adj} $  are  multiplicative,  and
that  $\hor_P$ and $\Gamma_{\inv}$
generate $\Omega(P)$,   we   conclude    that $\widehat{F}$ is
multiplicative, too.

 Let us prove that $\widehat{F}$ intertwines differentials. We
have
\begin{equation*}
\begin{split}
\widehat{F}d^{\wedge}(\varphi)&=\widehat{F}
\biggl(D(\varphi)+(-1)^{\partial\varphi}
\sum_k\varphi_k \pi(c_k)\biggr)\\
&=\sum_kD\varphi_k\otimes c_k +
(-1)^{\partial\varphi}\sum_k\varphi_k\otimes c_k^{(1)}
\pi(c_k^{(2)})\\
&\phantom{=}+(-1)^{\partial\varphi}\sum_k
\varphi_k\pi(c_k^{(3)})\otimes c_k^{(1)}
\k(c_k^{(2)})c_k^{(4)}\\
&=\sum_kD\varphi_k\otimes c_k
+(-1)^{\partial\varphi}\sum_k\varphi_k \otimes dc_k\\
&\phantom{=}+(-1)^{\partial\varphi}\sum_k
\varphi_k    \pi(c^{(1)}_k)\otimes
c_k^{(2)}\\
&=\sum_kd^{\wedge}(\varphi_k)\otimes c_k+(-1)^{\partial\varphi}
\sum_k\varphi_k\otimes dc_k \\
&=\bigl(d^{\wedge}\otimes \id +(-1)^{\partial\varphi}
\id\otimes  d\bigr)\widehat{F}(\varphi).
\end{split}
\end{equation*}
Further, if $\vartheta\in\Gamma_{\inv}$ then
\begin{equation*}
\begin{split}
\widehat{F}d^{\wedge}(\vartheta)&=
F^{\sstar} R(\vartheta)+\widehat{\adj} (d\vartheta)\\
&=(R\otimes \id) \adj (\vartheta)
+1\otimes d\vartheta+(d\otimes
\id) \adj (\vartheta)-(\id \otimes
d)\adj (\vartheta)\\
&=(d^{\wedge}\otimes
\id) \adj (\vartheta)-(\id \otimes
d)\adj (\vartheta)+1\otimes
d\vartheta\\
&=\bigl[d^{\wedge}\otimes
\id+(-1)^{\partial *}\id\otimes  d\bigr]
\bigl[1\otimes\vartheta+\adj (\vartheta)\bigr]\\
&=\bigl[d^{\wedge}\otimes \id +
(-1)^{\partial *}\id\otimes  d\bigr]\widehat{\adj}
(\vartheta).
\end{split}
\end{equation*}
Hence,
$\widehat{F}$ preserves differential structures. Finally,
\begin{equation*}
\begin{split}
\widehat{F}\bigl[(\varphi\otimes\vartheta)^*\bigr]&=\sum_k\widehat{F}
\bigl(\varphi_k^*
\otimes(\vartheta^*\circ c_k^*)\bigr)
=\sum_{ki}\varphi_k^* (\vartheta_i^*\circ c_k^{(1)*})
\otimes w_i^*c_k^{(2)*}\\
&=\sum_{ki}(\varphi_k\vartheta_i)^*\otimes w^*_ic_k^*
=\bigl[\widehat{F}(\varphi\otimes\vartheta)\bigr]^*
\end{split}
\end{equation*}
for each $\varphi\in \hor_P$ and
$\vartheta\in\Gamma_{\inv}^{\wedge}$,
where $\widehat{\adj} (\vartheta)=\Sum_i\vartheta_i\otimes w_i$.
Thus, $\widehat{F}$ is a hermitian map.
\end{pf}

 Clearly,   $\cal{B}=\Omega^0(P)$.  The   algebra   of
horizontal forms corresponding to $\Omega(P)$
coincides   with   the    initial    one.    In    other    words,
we can write $\hor_P=\hor(P)$. Further, $F^{\sstar}$ coincides with (the
restriction of) the corresponding right action $F^\wedge$.

 Let  us  consider  a  map $\omega\colon \Gamma_{\inv}
\rightarrow\Omega(P)$  given   by
$\omega(\vartheta)=1\otimes\vartheta$.
\begin{pro}
\bla{i} The map $\omega$ is a regular multiplicative connection
on  $P$. In particular, $\J(P)=\{0\}$.

\smallskip
\bla{ii} We have $R=R_{\omega}$  and $D=D_{\omega}$.
\end{pro}
\begin{pf}
It is evident that $\omega$ is a hermitian map. According to
the definition of $\widehat{F}$, we have
$$
\widehat{F}\omega(\vartheta)=\widehat{\adj} (\vartheta)=
\adj (\vartheta)+1\otimes
\vartheta=(\omega\otimes
\id) \adj (\vartheta)+1\otimes\vartheta.
$$
Hence,
$\omega$ is a  connection  on $P$.
Multiplicativity  of $\omega$ directly follows from the fact
that $\Gamma_{\inv}^{\wedge}$    is a  subalgebra  of
$\Omega(P)$. In particular,
$\omega^{\wedge}(\vartheta)=1\otimes\vartheta$ for
each  $\vartheta\in\Gamma_{\inv}^{\wedge}$. Regularity
follows  from   the
definition of the product in $\Omega(P)$. Finally, ({\it ii}\/) follows from
definitions of $R_{\omega}$, $D_{\omega}$ and $d^{\wedge}.$
\end{pf}

     The corresponding factorization map $m_{\omega}$ reduces to
the identity.
It is worth noticing that the construction of the algebra
$\ver(P)$ of verticalized forms can be understood as a trivial
special  case
of the construction presented in this subsection.  Indeed,  if  we
define  $\hor_P=\cal{B}$ (with $\hor_P^k=\{0\}$  for  $k\geq1$),
$D=0$ (and  hence $R=0$) and $F^{\sstar}=F$ then
$\Omega(P)=\ver(P)$ and $\widehat{F}=\widehat{F}_v$. The algebra
$(\vh(P),d_{v\!h})$ can be viewed in a similar way.

\appendix
\section{On Bicovariant Exterior Algebras}
In the presented theory we have assumed that the higher-order calculus
on the structure quantum group is described by  the  corresponding
universal envelope.  This  assumption  is  conceptually  the  most
natural.   However   all   the   formalism    can    be repeated
(straightforwardly,  or  with  natural   modifications)    if  the
higher-order calculus  on  $G$  is  described  by  an  appropriate
non-universal differential structure.

Here, it will be assumed that the higher-order calculus on $G$  is
based on the corresponding bicovariant exterior algebra \cite{W3}.
The appendix is devoted to the analysis of some aspects of
these structures, interesting from the point of view of
differential calculus on quantum principal bundles.

Let  $\Gamma$  be  a  bicovariant   first-order   differential
calculus over $G$ and let  us  consider  the  canonical  flip-over
automorphism $\sigma\colon \Gamma\otimes_{\cal A}\Gamma
\rightarrow\Gamma\otimes_{\cal   A}\Gamma$
(its  ``left-invariant'' restriction is given   by \eqref{flip}).

The    corresponding    {\it    exterior    algebra}     \cite{W3}
$\Gamma^\vee$
can be constructed by factorising
$\Gamma^\otimes$ through the ideal
$$S^\vee=\ker(A).$$
Here $ A=\Sum_{n}^\oplus A_n $ is the ``total antisymetrizer'' map, with
$A_n\colon \Gamma^{\otimes n}\rightarrow\Gamma^{\otimes n}$
given by
\begin{equation*}
A_n=\sum_{\pi\in S_n} (-1)^{\pi}\sigma_\pi
\end{equation*}
where  $\sigma_\pi$  is  the  operator   obtained   by   replacing
transpositions $i\leftrightarrow\,  i+1$  figuring  in  a  minimal
decomposition of $\pi$, by the corresponding $\sigma$-twists
(this definition  is  consistent, due to the braid equation  for
$\sigma$). By definition,  $A$ acts as the identity
transformation  on $\cal A$ and $\Gamma$.
The following
decomposition holds
\begin{equation}\label{antsym}
A_{k+l}=(A_k\otimes A_l)A_{kl}
\end{equation}
where
$$A_{kl}=\sum_{\pi\in S_{kl}} (-1)^{\pi}\sigma_{\pi^{-1}}$$
and $S_{kl}\subseteq  S_{k+l}$  is  the  subset  consisting  of
permutations preserving the order of the first $k$  and  the  last
$l$ factors.

The differential map $d\colon \cal   A\rightarrow\Gamma$ can be naturally
extended to the differential  on  the  whole  $\Gamma^\vee$.
By  universality,  there  exists  the  unique  graded-differential
homomorphism  ${\between}\colon \Gamma^\wedge\rightarrow\Gamma^\vee$
reducing to identity maps on $\cal A$ and $\Gamma$.
If
$\Gamma$ is *-covariant then the *-involutions on $\Gamma$  and
$\cal A$ can be extended to the *-structure on $\Gamma^\vee$ (so
that ${\between}$ is a hermitian map).

\begin{pro}\bla{i} There  exists  the  unique  differential   algebra
homomorphism
$\phi^\vee\colon \Gamma^\vee\rightarrow\Gamma^\vee\grten
\Gamma^\vee$ extending the map $\phi$.

\smallskip
\bla{ii} There exists  the  unique
graded-antimultiplicative extension $\k^\vee\colon \Gamma^\vee
\rightarrow \Gamma^\vee$
of the antipode $\k$, satisfying
\begin{equation}
\k^\vee d=d\, \k^\vee.\label{wk-d}
\end{equation}

\bla{iii} The following identities hold
\begin{gather}
(\phi^\vee\otimes
\id) \phi^\vee=(\id \otimes\phi^\vee)\phi^\vee
\label{w-coas}\\
m^\vee(\k^\vee\otimes \id) \phi^\vee
=m^\vee(\id  \otimes \k^\vee)\phi^\vee=1\e^\vee\label{w-k}\\
(\e^\vee\otimes \id) \phi^\vee=(\id \otimes \e^\vee)
\phi^\vee=\id \label{w-e}
\end{gather}
where $m^\vee$  is  the  product  map  in  $\Gamma^\vee$  and
$\e^\vee$ is a linear functional specified by
$$\e^\vee(\vartheta)=\e p_0(\vartheta).$$

\smallskip
\bla{iv} If $\Gamma$ is *-covariant then
\begin{gather}
(*\k^\vee)^2(\vartheta)=\vartheta \label{*wk}\\
\phi^\vee*=(*\otimes *)\phi^\vee\label{*wphi}.
\end{gather}
\end{pro}
\begin{pf} Uniqueness of maps $\phi^\vee$ and $\k^\vee$ is  a
consequence of the fact that $\Gamma^\vee$ is  generated,  as  a
differential algebra, by $\cal A$.
Let $\kG\colon \Gamma\rightarrow\Gamma$ be the  canonical  extension  of
the antipode map \cite{D}--Appendix B.  There  exists  the  unique
graded-antimultiplicative extension $\k^\otimes\colon \Gamma^\otimes
\rightarrow\Gamma^\otimes$ of $\k$ and $\kG$. We have
\begin{equation}\label{k-sigma}
\sigma \k^{\otimes}(\vartheta)=\k^{\otimes} \sigma(\vartheta)
\end{equation}
for each $\vartheta\in\Gamma\otimes_{\cal{A}}\Gamma$.
This directly follows from the fact that $\kG$ maps left-invariant
to right-invariant elements, and conversely. Hence,
\begin{equation}\label{k-inv}
\k^{\otimes}_n\sigma_\pi=\sigma_{j\pi j}\k^{\otimes}_n
\end{equation}
for   each   $\pi\in   S_n$,   where   $j$    is    the ``total
inverse''
permutation. This implies that operators
$\k^{\otimes n}$ commute with $A_n$.
Therefore  $\k^\otimes$  can  be  factorized  through   the   ideal
$S^\vee$. In such a way we  obtain a
graded-antimultiplicative
map $\k^\vee\colon \Gamma^\vee\rightarrow \Gamma^\vee$ satisfying
\eqref{wk-d}.

Let us now consider a $\cal{A}\otimes\cal{A}$-module $\Psi$ given by
$$ \Psi=\left(\cal A\otimes\Gamma\right)\oplus\left(\Gamma
\otimes\cal A\right).$$
It is easy to see that  $\Psi$  is  a  bicovariant  bimodule  over
the group $G\times G$. In particular, the corresponding right  and
left actions of $G\times G$ on $\Psi$ are given by
\begin{gather*}
\rig_\Psi\bigl((a\otimes\vartheta)\oplus(\eta\otimes b)\bigr)
=\bigl(\phi(a)\rig_\Gamma(\vartheta)\bigr)\oplus
\bigl(\rig_\Gamma(\eta)\phi(b)\bigr)\\
\ell_\Psi\bigl((a\otimes\vartheta)\oplus(\eta\otimes b)\bigr)
=\bigl(\phi(a)\ell_\Gamma(\vartheta)\bigr)\oplus
\bigl(\ell_\Gamma(\eta)\phi(b)\bigr),
\end{gather*}
where on the right-hand side the tensor multiplication is assumed.
The following natural isomorphism holds
$$ \Psi_{\inv}\cong\Gamma_{\inv}\oplus\Gamma_{\inv}.$$
Further, $\Psi$ is a first-order calculus over $G\times G$,  in  a
natural    manner.    The    corresponding    differential     map
$D\colon \cal{A}\otimes\cal{A}\rightarrow\Psi$ is given by
$D=d\otimes \id +\id\otimes  d$.
In terms of the above identification the corresponding right $\cal A
\otimes\cal A$-module structure on $\Psi_{\inv}$ is given by
\begin{equation}
(\vartheta\oplus\eta)\circ(a\otimes b)=
\e(a)(\vartheta\circ b)\oplus(\eta\circ a)\e(b)
\end{equation}
and  the  action  of  the  corresponding  flip-over  operator
$\Sigma\colon \Psi^{\otimes 2}_{\inv}\rightarrow\Psi^{\otimes 2}_{\inv}$
is determined by the block-matrix
\begin{equation}
\Sigma=\begin{pmatrix}\sigma&0&0&0\\
0&0&\tau&0\\
0&\tau&0&0\\
0&0&0&\sigma\end{pmatrix}
\end{equation}
where
$\tau\colon \Gamma_{\inv}^{\otimes2}\rightarrow\Gamma_{\inv}^{\otimes2}$
is the standard transposition. This implies
$$ \Psi^\vee\cong\Gamma^\vee\grten\Gamma^\vee.$$

Let $\fG\colon \Gamma\rightarrow\Psi$ be a map given by
\begin{equation*}
\fG=\ell_\Gamma\oplus\rig_\Gamma.
\end{equation*}
The following identities hold
\begin{align}
D\phi&=\fG d\label{D-d}\\
\fG(\vartheta a)&=\fG(\vartheta)\phi(a)\label{phi*mult1}\\
\fG(a\vartheta)&=\phi(a)\fG(\vartheta).\label{phi*mult2}
\end{align}

Equalities  \eqref{phi*mult1}--\eqref{phi*mult2}  imply  that
$\phi$ and $\fG$ can be consistently extended to a homomorphism
$\phi^\otimes\colon \Gamma^\otimes\rightarrow\Psi^\otimes$.
The following inclusion holds
$$\phi^\otimes\bigl(\ker(A)\bigr)\subseteq
\ker(A^\Sigma),$$
where $A^\Sigma$ is the antisymmetrizer corresponding to $\Sigma$.
Hence $\phi^\otimes$  can  be  projected
through  ideals $\ker(A)$  and
$\ker(A^\Sigma).$  In  such  a  way  we  obtain  the
homomorphism
$\phi^\vee\colon \Gamma^\vee\rightarrow\Gamma^\vee\grten
\Gamma^\vee$. Equality \eqref{D-d}
implies   that   $\phi^\vee$   intertwines    the    corresponding
differentials.

Properties \eqref{w-coas}--\eqref{w-e} as well as \eqref{*wk}--\eqref{*wphi}
simply follow from analogous properties  for  $\phi$
and $\k$. It is worth noticing that $\k^\vee$ and $\phi^\vee$ can be
obtained  by  projecting  $\k^\wedge$  and  $\widehat{\phi}$   from
$\Gamma^\wedge$ to $\Gamma^\vee$.
\end{pf}

Let us now turn to the conceptual framework of the previous paper,
and assume that $M$ is a classical compact  smooth  manifold,  and
$P$ a principal $G$-bundle over $M$. Further, let us  assume  that
$\Gamma$ is the  minimal  admissible  left-covariant  calculus
over $G$ (this calculus is bicovariant and *-covariant, too).
Let
$\tau=\left(\pi_U\right)_{U\in\cal U}$ be an arbitrary
trivialization  system
for $P$, and let ${\cal C}_\tau$ be the corresponding  $G$-cocycle,
consisting of ``transition functions'' $\psi_{U\!V}\colon
S(U\cap V)\otimes{\cal A}\rightarrow S(U\cap V)\otimes{\cal A}$.
The restrictions $\varphi_{U\!V}=\psi_{U\!V}\restr\cal{A}$ are
(*-homomorphisms) explicitly given by
$$\varphi_{U\!V}=(g_{V\!U}\otimes \id) \phi$$ where $g_{V\!U}\colon
U\cap V\rightarrow  G_{cl}$  represent  the  classical   $G_{cl}$-cocycle
corresponding  to  ${\cal   C}_\tau$   (understood   here   as
*-homomorphisms $g_{V\!U}\colon \cal A\rightarrow S(U\cap V)$).

There exists the unique map
$\tG_{U\!V}\colon \Gamma\rightarrow\bigl[\Omega^1(U\cap
V)\otimes\cal A\bigr]\oplus\bigl[ S(U\cap V)\otimes \Gamma\bigr]$
satisfying
\begin{align*}
\tG_{U\!V}(a\xi)&=\varphi_{U\!V}(a)\tG_{U\!V}(\xi)\\
\tG_{U\!V}(da)&=\left[d\otimes \id +\id\otimes
d\right]\varphi_{U\!V}(a),
\end{align*}
for each $a\in\cal A$ and $\xi\in\Gamma$. This implies also
$$ \tG_{U\!V}(\xi a)=\tG_{U\!V}(\xi)\varphi_{U\!V}(a) $$
and hence there exists (the unique) homomorphism
$\varphi_{U\!V}^\otimes\colon \Gamma^\otimes\rightarrow
\Omega(U\cap V)\grten
\Gamma^\otimes$     extending     both     $\varphi_{U\!V}$      and
$\tG_{U\!V}$.

All antisymetrizing operators are left and  right  covariant  in  a
natural  manner.  In  particular   they   are   reduced   in   the
corresponding spaces of left-invariant elements. In  what  follows
we shall denote by the same symbols their restrictions in
$\Gamma_{\inv}^\otimes$ (if there is no ambiguity from the context).
\begin{pro} We have
\begin{equation}\label{inv-A}
\varphi_{U\!V}^\otimes\left[S^\vee\right]\subseteq
\Omega(U\cap V)\grten S^\vee
\end{equation}
for each $(U,V)\in N^2(\cal U)$.
\end{pro}
\begin{pf}

The ideal $S^\vee$ is bicovariant. In
particular, it has the form
$$S^\vee\cong{\cal A}\otimes S^\vee_{\inv}$$
where $S^\vee_{\inv}$ is the
left-invariant part of $S^\vee$. The following equality holds
\begin{equation}\label{tran1}
\tG_{U\!V}(\vartheta)=1\otimes\vartheta+
\sum_k\partial^{U\!V}(\vartheta_k)\otimes c_k
\end{equation}
for each $\vartheta\in\Gamma_{\inv}$. Here, $\Sum_k\vartheta_k
\otimes c_k=\adj (\vartheta)$ and $\partial^{U\!V}\colon
\Gamma_{\inv}\rightarrow\Omega(U\cap V)$ is the map specified by
\begin{equation}\label{tran2}
\partial^{U\!V}\pi(a)=g_{V\!U}(a^{(1)})d\bigl(g_{U\!V}(a^{(2)})\bigr).
\end{equation}
Let us observe that
\begin{equation}\label{tran3}
\partial^{U\!V}(\eta)\partial^{U\!V}(\zeta)=
-\sum_k \partial^{U\!V}(\zeta_k)\partial^{U\!V}(\eta_k)
\end{equation}
for each $\eta,\zeta\in\Gamma_{\inv}$, where
$\Sum_k\zeta_k\otimes\eta_k=\sigma(\eta\otimes\zeta)$.
Indeed this follows from \eqref{flip}, and from the identity
\begin{equation}\label{tran4}
\partial^{U\!V}(\vartheta\circ a)=\e(a)\partial^{U\!V}(\vartheta).
\end{equation}

Now for an arbitrary $\vartheta\in\Gamma_{\inv}^{\otimes n}$
let us consider the elements
$\vartheta_k=(\id \otimes p_k)\varphi_{U\!V}^\otimes(\vartheta)$
where $0\leq k\leq n$. It follows from \eqref{flip} and \eqref{tran1}
that these  elements  have the form
\begin{equation}
\vartheta_k=(\nabla^{U\!V}_{l}\otimes \id^k)A_{lk}(\vartheta)
\end{equation}
with $l=n-k$. Here we have identified $\Gamma^\otimes$ and
$\cal A\otimes\Gamma^\otimes_{\inv}$, and $\nabla_l^{U\!V}$
are components of the unital multiplicative extension
of the map $\nabla^{U\!V}=(\partial^{U\!V}\otimes \id) \adj .$
According to \eqref{tran3} the above expression can  be  rewritten
as
$$ \vartheta_k=\frac{1}{l!}(\nabla^{U\!V}_{l}
A_{l}\otimes
\id^k)A_{lk}(\vartheta).
$$
In particular
\begin{equation}
(\id \otimes A_k)
(\vartheta_k)=\frac{1}{l!}
(\nabla^{U\!V}_{l}\otimes \id^k)A_n
(\vartheta),
\end{equation}
according to \eqref{antsym}. Hence, \eqref{inv-A} holds.
\end{pf}

The map $\varphi_{U\!V}^\otimes$ can be factorized through the ideal
$S^\vee$. In such a way we obtain a homomorphism
$\varphi_{U\!V}^\vee\colon \Gamma^\vee\rightarrow\Omega(U\cap V)\widehat{
\otimes}\Gamma^\vee$  of  graded-differential *-algebras.  Now  the
formula
\begin{equation}
\psi_{U\!V}^\vee(\alpha\otimes\vartheta)=\alpha\varphi_{U\!V}^\vee(
\vartheta)
\end{equation}
defines a graded-differential *-automorphism of $\Omega(U
\cap V)\grten\Gamma^\vee$ (extending both
$\psi_{U\!V}$ and $\varphi_{U\!V}^\vee$).
These maps are $\Omega(U\cup V)$-linear and satisfy the following
cocycle conditions
$$\psi_{U\!V}^\vee\psi_{V\!W}^\vee(\varphi)=\psi_{U\!W}^\vee(\varphi)$$
for each $(U,V,W)\in N^3(\cal U)$ and
$\varphi\in\Omega_c(U\cap V\cap W)\grten\Gamma^\vee$.

Applying the above results, and using similar constructions as in
the   previous   paper,   it   is   possible   to   construct    a
graded-differential *-algebra $\Lambda(P)$, representing the
corresponding ``differential forms'' on the bundle. This algebra  is
locally trivial,  in  the  sense  that  any  local  trivialization
$(U,\pi_U)$  of  the  bundle  can  be  ``extended''   to   a   local
representation of the form
$\Omega(U)\grten\Gamma^\vee\leftrightarrow\Lambda(P\vert_U)$.

By construction, the right action  $F$  can  be (uniquely)
extended  to  the
homomorphism
$F^\vee
\colon \Lambda(P)\rightarrow\Lambda(P)\grten\Gamma^\vee$
of  graded-differential   *-algebras.   Moreover,   all   entities
naturally appearing in the differential calculus on $P$
constructed from the universal envelope $\Gamma^\wedge$ have
counterparts  in  the  calculus  based  on  $\Lambda(P)$,  and  all
algebraic relations are preserved (because the  formalism  can
be obtained by ``projecting'' the first calculus on $\Lambda(P))$.

     In a certain sense,  $\Gamma^\vee$  is  {\it the  minimal}
bicovariant graded algebra (built over $\Gamma$) compatible  with
all  possible ``transition
functions'' for the bundle $P$.   Namely,   let   us    assume    that
$\cal{N}\subseteq  \Gamma^\otimes$   is   a   bicovariant   graded-ideal
satisfying $\cal{N}^k=\{0\}$ for $k\in\{0,1\}$ and
\begin{equation}\label{rtrv}
\varphi^\otimes_{U\!V}(\cal{N})\subseteq\Omega(U\cap V)\grten \cal{N}
\end{equation}
for each trivialization system  $\tau=\left(\pi_U\right)_{U\in\cal
U}$ and each $(U,V)\in N^2(\cal U)$. This  ensures  the  possibility  to
construct the corresponding global algebra  for  the  bundle  $P$,
which locally will be of the form $\Omega(U)\grten\bigl[
\Gamma^\otimes/\cal{N}\bigr].$
\begin{lem}\label{lem:min}
Under the above assumptions we have
\begin{equation}
\cal{N}\subseteq S^\vee.
\end{equation}
\end{lem}
\begin{pf} We shall  prove  inductively  that
\begin{equation}\label{inv-ext}
\cal{N}_{\inv}^k\subseteq S^{\vee k}_{\inv},
\end{equation}
for each $k\geq 2$. Let us assume that
$\Sum_i\vartheta_i\otimes\eta_i=\psi\in \cal{N}^2_{\inv}$. Applying
\eqref{rtrv} and \eqref{tran1} we obtain
\begin{equation*}
0=\sum_{ikl}\partial^{U\!V}(\vartheta_{ik})\partial^{U\!V}(\eta_{il
})\otimes c_{ik}d_{il}+
\sum_{ik}\partial^{U\!V}(\vartheta_{ik})\otimes
c_{ik}\eta_i-\sum_{il}\partial^{U\!V}(\eta_{il})\otimes
\vartheta_i d_{il}
\end{equation*}
where $\adj (\vartheta_i)=\Sum_k\vartheta_{ik}\otimes
c_{ik}$ and
$\adj (\eta_i)=\Sum_l\eta_{il}\otimes d_{il}.$
In particular
$$0=\sum_{ik}\partial^{U\!V}(\vartheta_{ik})\otimes
c_{ik}\eta_i-\sum_{il}\partial^{U\!V}(\eta_{il})\otimes
\vartheta_id_{il}.
$$
In other words, modulo the identification
$\Gamma\leftrightarrow\cal A\otimes \Gamma_{\inv}$, we have
\begin{equation*}
(\nabla^{U\!V}\otimes \id) (I -\sigma)(\psi)=0.
\end{equation*}
Having in mind that the family of maps $\nabla^{U\!V}$ distinguishes
elements of $\Gamma_{\inv}$ (a consequence  of  the  minimality  of
$\Gamma$)  we  conclude  that  $\psi=\sigma(\psi)$.  In  other
words
$\cal{N}_{\inv}^2\subseteq \ker(I -\sigma)=S^{\vee 2}_{\inv}$.

Let us assume that  \eqref{inv-ext}  holds  for  some  $k\geq  2$.
Then
$$
0=(\id \otimes A_kp_k)\varphi^\otimes_{U\!V}(\psi)=(\nabla^{U\!V}
\otimes A_k)A_{1k}(\psi)
$$
for each $\psi\in \cal{N}_{\inv}^{k+1}.$ Because of the
arbitrariness of $\tau$ it follows that
$$A_{k+1}(\psi)=(\id \otimes A_k)A_{1k}(\psi)=0,$$
in other words $\cal{N}_{\inv}^{k+1}\subseteq (S^{\vee}_{\inv})^{k+1}$.
\end{pf}

We  pass  to  the  study  of   the   problem   of   passing   from
$\Gamma^\wedge$  to  $\Gamma^\vee$,  in  the  framework of the
general theory.

Let $P=(\cal B,i,F)$ be a  quantum  principal  $G$-bundle  over  a
quantum space $M$.

\begin{lem}
Modulo the natural identification $\Gamma^{\wedge
}\leftrightarrow{\cal A}\otimes\Gamma^{\wedge}_{\inv}$ we have
\begin{equation}
(\id \otimes p_k)\widehat{F}\omega^{\otimes n}=
(F^\wedge\omega^{\otimes l}\otimes [\,]_k^\wedge)
A_{lk}
\end{equation}
where $n=k+l$, and $\omega$ is an arbitrary connection on $P$.
\end{lem}
\begin{pf}
Essentially the same reasoning as in the  proof  of
Proposition A.2, applying identity \eqref{dxFcon} instead of
\eqref{tran1}.
\end{pf}

 The above equation implies
\begin{equation}\label{FwA}
\widehat{F}\omega^\otimes\left(S^\vee_{\inv}\right)
\subseteq
\omega^\otimes\left(S^\vee_{\inv}\right)\grten
\Gamma^\wedge
+\Omega(P)\grten[S^\vee]^\wedge,
\end{equation}
for each $\omega\in\con(P)$.
Let us assume that the calculus  (described  by
$\Omega(P)$) admits regular  connections,  and  that  $\J(P)=\{0\}$
(multiplicativity of regular connections).

Let $\JJ\subseteq\Omega(P)$ be the space linearly  generated  by
elements of  the  form  $\varphi\omega^\otimes(\vartheta)$  where
$\vartheta\in S^\vee_{\inv}$ and $\varphi\in\hor(P)$, while
$\omega\in\rc(P)$.

\begin{lem} \bla{i} The space $\JJ$ is a (two-sided)
graded-differential
*-ideal in $\Omega(P)$,
independent  of  the  choice  of  a  regular
connection $\omega$.

\smallskip
\bla{ii} We have
\begin{equation}\label{FIG}
\widehat{F}\bigl(\JJ\bigr)\subseteq \JJ\grten\Gamma^\wedge+
\Omega(P)\grten[S^\vee]^\wedge.
\end{equation}
\end{lem}
\begin{pf}
The space $m_\omega^{-1}\bigl(\JJ\bigr)=
\hor(P)\otimes[S_{\inv}^\vee]^\wedge$ is  a  graded two-sided
*-ideal  in  $\vh(P)$. Hence the space
$\JJ$ is the graded *-ideal in $\Omega(P)$.
Inclusion \eqref{FIG} directly follows from \eqref{FwA}.

Let us prove that $d\bigl(\JJ\bigr)\subseteq \JJ$. It  is
sufficient  to  check
that $d\omega^\otimes(\vartheta)\in \JJ$ for each $\vartheta\in
S^\vee_{\inv}$.  The following equality holds
\begin{equation}
d\omega^{\otimes}(\vartheta)=\omega^{\otimes}\delta^*
(\vartheta)+m_\Omega \left(R_\omega\otimes\omega^{\otimes}
\right)A_{1k-1}(\vartheta)
\end{equation}
for each $k\geq 2$ and $\vartheta\in \Gamma_{\inv}^{\otimes k}$. Here,
$\delta^*\colon \Gamma_{\inv}^\otimes\rightarrow\Gamma_{\inv}^\otimes$
is  the  unique (hermitian) antiderivation
extending a given embedded
differential map $\delta$. In particular,  if
$\vartheta\in S^{\vee k}_{\inv}$
then both summands on the right-hand side of  the  above  equality
belong to $\JJ$, because of
$\delta^*(S^{\vee k}_{\inv})\subseteq
(S_{\inv}^{\vee})^{k+1}$
and $A_{1k-1}(S^{\vee k}_{\inv})\subseteq
\Gamma_{\inv}\otimes (S^{\vee}_{\inv})^{k-1}$.

Finally, let $\zeta$ be an arbitrary (hermitian) tensorial
1-form satisfying \eqref{vreg1}. We have then
$$(\omega+\zeta)^\otimes(\vartheta)=
\sum_{k+l=n}\frac{1}{k!}m_\Omega(\zeta^{\otimes k} A_k\otimes
\omega^{\otimes l})A_{kl}(\vartheta)$$
for each $\vartheta\in\Gamma_{\inv}^{\otimes n}$. This shows that
$\JJ$ is $\omega$-independent.
\end{pf}

Hence, it is possible to pass jointly to the
factoralgebra $\Lambda(P)=\Omega(P)/S^\vee$ (as  a  representative
of the calculus on the bundle), and to  the  exterior
bicovariant algebra $\Gamma^\vee$ (representing  the  calculus  on
$G$).  This  factorization  does  not
change the first-order differential structure.  For  this  reason
the spaces of connection forms associated to both
calculi on $P$ are the same. Further, the
spaces $\hor(P)$ and $\Omega_1(P)$ are preserved.
This  implies  that  regular  connections   relative   to
$\Omega(P)$  and  $\Lambda(P)$  are  the  same.  By  construction,
regular connections are multiplicative, relative  to  $\Lambda(P)$,
too.

Lemma~\ref{lem:min} establishing the
``minimality'' of the exterior  algebra
has a general quantum counterpart.

Let us assume that  the  calculus  on  the  bundle  is  such  that
tensorial 1-forms $\zeta$ satisfying  \eqref{vreg1}  distinguish
elements of $\Gamma_{\inv}$. In the case of bundles over  classical
compact manifolds the algebra $\Omega(P)$ built from the minimal
admissible calculus $\Gamma$ possesses this property.

Let us consider a bicovariant graded-ideal
$\cal{N}\subseteq\Gamma^\otimes$ satisfying $\cal{N}^k=\{0\}$ for
$k\in\{0,1\}$. The space $\cal{N}_\omega=m_\omega\bigl(\hor(P)\otimes
[\cal{N}_{\inv}]^\wedge\bigr)$
is a two-sided ideal in $\Omega(P)$. Let us assume that $\cal{N}_\omega$
is   independent   of   $\omega\in\rc(P)$.
\begin{lem}
Under the above assumptions $\cal{N}\subseteq S^\vee.\qed$
\end{lem}

The construction of ``global'' differential calculus on the bundle
given in Subsection 6 can be applied  to the
exterior algebra case, too.

Starting from  algebras  $\hor_P$  and  $\Gamma_{\inv}^\vee$  it  is
possible  to   construct   a   graded   *-algebra   structure   on
$\Lambda(P)=\hor_P\otimes\Gamma_{\inv}^\vee$.
Using    maps    $R$,     $D$,
and   $d\colon \Gamma_{\inv}^\vee\rightarrow\Gamma_{\inv}^\vee$
it is possible to construct a natural differential $d^\vee$ on
$\Lambda(P)$. All constructions are the same as
in the universal envelope case. The only nontrivial  point  is  to
prove the analog of Lemma~\ref{lem:R+d}.

Let $R^\sstar\colon \Gamma^\otimes_{\inv}\rightarrow \Lambda(P)$
be the unique antiderivation extending $R$. The following equality
holds
\begin{equation}
R^\sstar(\eta)=\bigl(R\otimes [\,]^{\vee n}\bigr)A_{1n}(\eta)
\end{equation}
for each  $\eta\in(\Gamma_{\inv}^{\otimes})^{n+1}$.  This  implies  that
$R^\sstar$ can be factorized through $S^\vee_{\inv}$. In such a way
we obtain the map $R^\vee\colon \Gamma_{\inv}^\vee\rightarrow\Lambda(P)$.
Finally, the  map  $d^\vee\colon \Gamma_{\inv}^\vee\rightarrow\Lambda(P)$
is defined by
$$d^\vee(\vartheta)=d(\vartheta)+R^\vee(\vartheta).$$

\section{Multiple Irreducible Submodules}

     In this appendix the structure of the *-algebra $\cal B$
of  functions on a quantum principal bundle $P=(\cal{B},i,F)$
will be analyzed from the  viewpoint
of the representation theory \cite{W2} of the structure quantum group.

Let $\cal T$ be the set of (equivalence classes of)
irreducible unitary representations of $G$.

     The representation
$F\colon  \cal{B}\rightarrow\cal{B}\otimes\cal{A}$
of $G$  in $\cal{B}$  is  highly  reducible.  For  each
$\alpha\in\cal T$ let
$\cal{B}^\alpha \subseteq\cal{B}$
be  the  multiple  irreducible  subspace
corresponding to $\alpha$. Evidently,
$\cal{B}^\alpha$  is a bimodule over $\cal{V}$. We have
$$
\cal{B}=\sideset{}{^\oplus}\sum_{\alpha\in\cal T}\cal{B}^\alpha.
$$

     For an arbitrary $\alpha\in\cal T$, let
$u^\alpha\colon  \Bbb{C}^n\rightarrow\Bbb{C}^n \otimes\cal{A}$
be a representative  of  this  class (where $n$ is the dimension of
$\alpha$).
Let  $(e_1,\dots,e_n)$  be  the
absolute basis in $\Bbb{C}^n$. We have
$$
u^\alpha(e_i)=\sum_{j=1}^ne_j\otimes u_{ji}^\alpha,
$$
where $u^\alpha_{ij}\in\cal{A}$ are the corresponding
matrix elements.

     Let us consider the space $\Mor(u^\alpha,F)$
consisting of intertwiners (morphisms)
$\varphi\colon  \Bbb{C}^n\rightarrow\cal{B}$ between
representations $u^\alpha$ and $F$. The
space $\Mor(u^\alpha,F)$ is a $\cal{V}$-bimodule, in a natural  manner.
The  maps $\varphi$  take  values  from  $\cal{B}^\alpha$.

Bimodules $\Mor(u^\alpha,F)\otimes\Bbb{C}^n$ and $\cal{B}^\alpha$
are naturally isomorphic, via the correspondence
\begin{equation}
\varphi\otimes x\leftrightarrow\varphi(x).
\end{equation}

     This identification intertwines $F\restr\cal{B}^\alpha$ and
$\id\otimes  u^\alpha$.
Irreducible $G$-multiplets in $\cal{B}^\alpha$
are  of  the  form
$\bigl\{\varphi(e_1),\dots,\varphi(e_n)\bigr\}$, for some intertwiner
$\varphi$.

Let us fix $i,j\in \{1,\dots,n\}$. There  exist  intertwiners
$\mu_1,\dots,\mu_d\in \Mor(u^\alpha,F)$ and
numbers $c_{klpq}\in\Bbb{C}$ (where $k,l\in \{1,\dots,d\}$ and $p,q\in
\{1,\dots,n\}$)
such that
\begin{equation}
\sum_{klpq}c_{klpq}b_{kp}^{\alpha*}F(b_{lq}^\alpha)=
1\otimes u^\alpha_{ij},
\end{equation}
where $b_{kp}^\alpha=\mu_k(e_p)$. Equivalently
\begin{equation}
\sum_{klpq}c_{klpq}F(b_{kp}^{\alpha*})b_{lq}^\alpha=1\otimes
u^{\alpha*}_{ji}.
\end{equation}

     The above equalities imply that the  summation  over  indexes
satisfying $(p,q)\neq (i,j)$ can be dropped. In other words,
\begin{gather}
\sum_{kl}c_{kl}b_{ki}^{\alpha*}F(b_{lj}^\alpha)=1\otimes
u^\alpha_{ij}\label{A4}\\
\sum_{kl}c_{kl}F(b_{ki}^{\alpha*})b_{lj}^\alpha=1\otimes
u^{\alpha*}_{ji}\label{A5}
\end{gather}
where $c_{kl}=c_{klij}$. From the fact that
$\bigl\{b_{k1}^\alpha,\dots,b_{kn}^\alpha\bigr\}$ form  a $G$-multiplet
it follows that the  above  equalities  hold  for
arbitrary $i,j\in \{1,\dots,n\}$.

     Equalities \eqref{A4}--\eqref{A5} imply that the
numbers $c_{kl}$ can
be allways choosen  such  that  the  matrix  $(c_{kl})$ is
hermitian.

Further, without a lack of  generality  we  can  assume  that  the
matrix $(c_{kl})$ is nonsingular. In what follows, it will  be  assumed
that the matrix $(c_{kl})$ is {\it positive}.
Diagonalizing this matrix,  and
redefining in the appropriate way intertwiners $\mu_k$ we obtain
\begin{gather}
\sum_kb_{ki}^{\alpha*}F(b_{kj}^\alpha)=1\otimes
u^\alpha_{ij}\label{A6}\\
\sum_kF(b_{ki}^{\alpha*})b_{kj}^\alpha=1\otimes
u^{\alpha*}_{ji}.\label{A7}
\end{gather}

     Equivalently, the following identities hold
\begin{equation}\label{A8}
\sum_kb_{ki}^{\alpha*}b_{kj}^\alpha=\delta_{ij}.
\end{equation}

     In particular, all irreducible representations appear  in the
decomposition of $F$ into irreducible  components  (all
$\cal{V}$-bimodules
$\cal{B}^\alpha$ are  non-trivial).

     Let us consider, for each $f\in\cal{V}$, the elements
\begin{equation}
\varrho_{kl}(f)=\sum_ib_{ki}^\alpha fb_{li}^{\alpha*}.\label{A9}
\end{equation}
These elements are $F$-invariant and the following identities hold
\begin{align}
\varrho_{kl}(f)^*&=\varrho_{lk}(f^*)\label{A10}\\
\sum_n\varrho_{kn}(f)\varrho_{nl}(g)&=\varrho_{kl}(fg).\label{A11}
\end{align}
Indeed, we have
$F(b_{ki}^\alpha)=\sum_jb_{kj}^\alpha\otimes u^\alpha_{ji}$,
and a direct computation gives
\begin{multline*}
F\bigl(\varrho_{kl}(f)\bigr)=\sum_iF(b_{ki}^\alpha
fb_{li}^{\alpha*})=\sum_{imn}b_{km}^\alpha f
b_{ln}^{\alpha *}\otimes u^\alpha_{mi}u^{\alpha*}_{ni}\\
=\sum_{mn}b_{km}^\alpha fb_{ln}^{\alpha*}\otimes
\delta_{mn}=\varrho_{kl}(f)\otimes 1,\end{multline*}
and similarly
$$
\sum_n \varrho_{kn}(f)\varrho_{nl}(g)=\sum_{nij}b_{ki}^\alpha
fb_{ni}^{\alpha*}
b_{nj}^\alpha gb^{\alpha*}_{lj}=
\sum_{ij}b_{ki}^\alpha fgb_{lj}^{\alpha*}\delta_{ij}
=\sum_ib_{ki}^\alpha fgb_{li}^{\alpha*}=\varrho_{kl}(fg).
$$

In  other  words,  the   maps
$\varrho_{kl}\colon  \cal{V}\rightarrow\cal{V}$
realize a  *-homomorphism  $\varrho\colon
\cal{V}\rightarrow M_d(\cal{V})$,
where $M_d(\cal{V})$ is the *-algebra  of $d\times d$-matrices
over  $\cal{V}.$  In
particular, $\Plft=\varrho(1)$ is a projector in $M_d(\cal{V}).$

     Let us consider the free left $\cal{V}$-module $\cal{V}^d$,
with
the  absolute  basis  $(\varepsilon_1,\dots,\varepsilon_d).$
The elements of the algebra
$M_d(\cal{V})$ are understandable  as
endomorphisms (acting on the right)  of $\cal{V}^d$,
in a natural manner. Explicitly, this realization is  given
by
\begin{equation}
(\varepsilon_k)A=\sum_l A_{kl}\varepsilon_l.
\end{equation}
Let $\cal{E}\subseteq\cal{V}^d$  be the left
$\cal{V}$-submodule determined by the projector $\Plft$.
Evidently, $\cal E$ is
$\varrho$-invariant and $\varrho$, together  with
the left multiplication,  determine  a  (unital)  $\cal{V}$-bimodule
structure on $\cal E$.

     Let $\Jlft\colon  \cal{V}^d \rightarrow \Mor(u^\alpha,F)$  be the left
$\cal{V}$-module homomorphism given by
\begin{equation}
\Jlft(\varepsilon_k)=\mu_k.
\end{equation}

     The following identity holds
\begin{equation}
\Jlft\bigl(\psi \varrho(f)\bigr)=\Jlft(\psi)f,
\end{equation}
in particular $\Jlft(\psi \Plft)=\Jlft(\psi)$.

This implies that the  restriction
$(\Jlft\restr \cal{E})\colon  \cal{E}\rightarrow \Mor(u^\alpha,F)$
is   a    homomorphism    of    unital $\cal{V}$-bimodules.

     Now we shall prove that $\Jlft\restr\cal{E}$ is bijective.
Let us  assume that $\psi\in \ker(\Jlft)$. This implies
$$
\sum_{ki}\psi_kb_{ki}^\alpha b_{li}^{\alpha*}=0,
$$
for each $l\in\{1,\dots,d\}$ where $\psi=\Sum_k\psi_k\varepsilon_k.$  In
other  words,
$\psi\Plft=0$,  which  means  that  $\Jlft\restr \cal{E}$  is
injective. We have
$$\mu=\sum_{k}q_k\mu_k, $$
for each $\mu\in \Mor(u^\alpha,F)$,
where $q_k\in\cal{V}$ are elements given by
$$ q_k =\sum_i \mu(e_i)b_{ki}^{\alpha*}. $$
In other words, elements $\mu_k$  span the left
$\cal{V}$-module $\Mor(u^\alpha,F)$.
This implies that $\Jlft$ is surjective.

     Hence, $\Mor(u^\alpha,F)$  are finite and  projective,  as
left  $\cal{V}$-modules.
This implies that left $\cal{V}$-modules $\cal{B}^\alpha$
are  finite  and  projective.

     Relations \eqref{A8} can be rewritten in the form
\begin{equation}
B^\dagger B=I_n,
\end{equation}
where $B$ is a $d\times n$ matrix with coefficients $b_{ki}^\alpha$ and
$I_n$  is  the  unit matrix in $M_n(\cal B)$.

     Let us assume that the following {\it additional relations} hold
\begin{equation}\label{mod-P}
(ZBC^{-1})^\top B^*=I_n,
\end{equation}
where $C\in M_n(\Bbb{C})$ is the canonical intertwiner
\cite{W2} between $u^\alpha$ and its second  contragradiant
$u^{\alpha}_{cc}$,
and $Z\in M_d(\Bbb{C})$ is a strictly positive matrix.

     Relations of this type naturally  appear  in  a $C^*$-algebraic
version of the theory of quantum principal bundles. The matrix
$Z$ is connected with modular properties of an appropriate invariant
integral on the bundle.

     Let us consider a map
$\lambda\colon  \cal{V}\rightarrow M_d(\cal{V})$ given by
\begin{equation}
\lambda_{kl}(f)=\sum_ib_{ki}^{\alpha*}f[ZBC^{-1}]_{li}.
\end{equation}

     We have then
\begin{gather}
\lambda(f)^\dagger=Z^*\lambda(f^*)(Z^*)^{-1}\\
\lambda(f)\lambda(g)=\lambda(fg).
\end{gather}
Further, the elements of $M_d(\cal{V})$ are naturally identificable
with  endomorphisms of the right $\cal{V}$-module $\cal{V}^d$
(acting  on the left). Let $\cal F\subseteq\cal{V}^d$
be  a  right  $\cal{V}$-submodule
determined by a projection $\Prig=\lambda(1)$. The map $\lambda$
induces a unital left $\cal{V}$-module structure on $\cal F$,   so
that $\cal F$ becomes a
$\cal{V}$-bimodule. Let $\Jrig\colon  \cal{V}^d\rightarrow \Mor(u^\alpha,F)$
be a right $\cal{V}$-module  homomorphism  given  by
$\Jrig(\varepsilon_k)=\nu_k$, where $\nu_k=\Sum_l Z_{kl}\mu_l$. Then the
restriction $(\Jrig\restr\cal{F})\colon  \cal F\rightarrow \Mor(u^\alpha,F)$
is a bimodule isomorphism.

A similar consideration can be applied to all covariant
algebras figuring in the game. In  particular
$$
\hor(P)=\sideset{}{^\oplus}\sum_{\alpha\in\cal T}\cal{H}^\alpha
$$
where $\cal{H}^\alpha$ are corresponding
$\alpha$-multiple irreducible subspaces (relative to
the decompositions of $F^\wedge$).
These spaces are $\Omega(M)$-bimodules.
The following natural decompositions hold
$$
\cal{H}^\alpha=\Mor(u^\alpha,F^\wedge)\otimes\Bbb{C}^n.
$$

Every
$\varphi\in \hor(P)$ can be written in the form
\begin{equation}
\varphi=\sum_k w_k b_k
\end{equation}
where $w_k\in\Omega(M)$ and $b_k\in\cal B$.
Indeed, it is sufficient check the statement for elements of
some irreducible multiplet.

Let us assume that $\bigl\{\varphi_1,
\dots, \varphi_n\bigr\}\in\cal{H}^\alpha$ is an irreducible
$\alpha$-multiplet.
We have then
$$ \varphi_i=\sum_{jk}\varphi_j b^{\alpha*}_{kj}b^\alpha_{ki}$$
for each $i\in\{1,\dots,n\}$. On the other hand the elements $
\Sum_j\varphi_j b^{\alpha*}_{kj}$ belong to $\Omega(M)$.

In particular, if $\hor^+(P)$ is generated by $\hor(P)$ and if every
first-order horizontal form $\varphi$ can be written as
$\varphi=\Sum_kb_kd(g_k)$, where $b_k\in\cal{B}$ and $g_k\in\cal{V}$,
then the spaces $\cal{H}^\alpha$ are linearly spanned by
elements of the form $h^\alpha=b^\alpha d(f_1)\dots d(f_n)$, where
$f_k\in\cal{V}$ and $b^\alpha\in\cal{B}^\alpha$.


\begin{thebibliography}{44}
\bibitem{C} Connes A: {\it Geometrie non commutative}, InterEditions
Paris (1990)
\bibitem{D} {\Dj}ur{\dj}evi\'c M: {\it Geometry of Quantum Principal
Bundles I}, Preprint QmmP 6/92 Belgrade University
\bibitem{BM} Brzezinski T, Majid S: {\it Quantum Group Gauge Theory on
Quantum Spaces}, Preprint DamtP 92-27
\bibitem{W-su2} Woronowicz S L: {\it Twisted SU(2) group. An example of
a non-commutative differential calculus}, RIMS Kyoto University {\bf 23}
117--181 (1987)
\bibitem{W2} Woronowicz S L: {\it Compact Matrix Pseudogroups}, CMP
{\bf 111} 613--665 (1987)
\bibitem{W3} Woronowicz S L: {\it Differential Calculus on Compact
Matrix Pseudogroups/ Quantum Groups}, CMP {\bf 122} 125--170 (1989)
\bibitem{P1} Podles P: {\it Quantum Spheres}, Lett Math Phys {\bf 14}
193--202 (1987)
\bibitem{KN} Kobayashi S, Nomizu K: {\it Foundations
of Differential Geometry}, Interscience Publishers,
New York London (1963)
\end{thebibliography}
\end{document}